\documentclass[leqno]{amsart}
\usepackage[utf8]{inputenc}
\usepackage{newclude}
\usepackage{amsfonts,amssymb,amsmath,amsgen,dsfont}
\allowdisplaybreaks

\usepackage{xcolor,hyperref,soul}
\definecolor{plotblue}{RGB}{0,114,178}
\definecolor{plotvermillion}{RGB}{213,94,0}
\usepackage{tikz}
\usepackage{pgfplots}
\usepackage{pgfplotstable}
\pgfplotsset{compat=1.18}
\usetikzlibrary{arrows.meta}

\def\vortexdegree{1}

\pgfplotstableread{vortex_profiles.dat}\vortexdata

\pgfplotstablegetrowsof{\vortexdata}
\pgfmathtruncatemacro{\vortexlastrow}{\pgfplotsretval-1}

\newcommand{\getvortexamp}[1]{%
  \ifcase\vortexdegree
    \pgfplotstablegetelem{#1}{f0}\of\vortexdata
  \or
    \pgfplotstablegetelem{#1}{f1}\of\vortexdata
  \or
    \pgfplotstablegetelem{#1}{f2}\of\vortexdata
  \or
    \pgfplotstablegetelem{#1}{f3}\of\vortexdata
  \or
    \pgfplotstablegetelem{#1}{f4}\of\vortexdata
  \else
    \PackageError{vortex}
      {vortexdegree must be 0,1,2,3,4}{}
  \fi
}

\usepackage[normalem]{ulem}
\usepackage{cancel}
\usepackage{comment}

\usepackage{enumerate}

\usepackage{relsize} 
\usepackage{bigints}
\usepackage{bm}

\theoremstyle{plain}
\newtheorem{theorem}{Theorem}[section]
\newtheorem{definition}[theorem]{Definition}

\newtheorem{lemma}[theorem]{Lemma}
\newtheorem{corollary}[theorem]{Corollary}
\newtheorem{proposition}[theorem]{Proposition}

\theoremstyle{remark}
\newtheorem{example}[theorem]{Example}
\newtheorem{problem}{Problem}
\newtheorem{remark}[theorem]{Remark}

\def\C{{\mathbb C}}
\def\R{{\mathbb R}}
\def\T{{\mathbb T}}
\def\N{{\mathbb N}}
\def\Z{{\mathbb Z}}

\newcommand*\dd{\mathop{}\!\mathrm{d}}

\def\({\left(}
\def\){\right)}
\def\<{\left\langle}
\def\>{\right\rangle}

\def\1{{\mathbf 1}}

\let\epsilon\varepsilon
\let\geq\geqslant
\let\leq\leqslant

\newcommand{\K}{\mathbb{K}}
\newcommand{\conj}{\overline}

\newcommand{\eul}{\mathrm{e}}
\newcommand{\imu}{\mathrm{i}}

\DeclareMathOperator{\RE}{Re}
\DeclareMathOperator{\IM}{Im}
\DeclareMathOperator{\diver}{div}

\DeclareMathOperator{\dist}{dist}
\DeclareMathOperator{\curl}{curl}

\DeclareMathOperator{\sign}{sign}
\DeclareMathOperator{\Lip}{Lip}

\DeclareFontFamily{U}{mathx}{\hyphenchar\font45}
\DeclareFontShape{U}{mathx}{m}{n}{
      <5> <6> <7> <8> <9> <10>
      <10.95> <12> <14.4> <17.28> <20.74> <24.88>
      mathx10
      }{}
\DeclareSymbolFont{mathx}{U}{mathx}{m}{n}
\DeclareFontSubstitution{U}{mathx}{m}{n}
\DeclareMathAccent{\widecheck}{0}{mathx}{"71}

\newcommand{\ha}[1]{\widehat{#1}}
\newcommand{\supp}{\operatorname{supp}}
\newcommand{\loc}{\mathrm{loc}}

\let\leq\leq
\let\geq\geq

\numberwithin{equation}{section}

\date\today

\title[Boundedness and Continuity of Phase-corrected Gradients]{Boundedness and Continuity of Phase-corrected Gradients on Sobolev Spaces}

\author[L.E. Hientzsch]{Lars Eric Hientzsch}
\address{
Institute for Analysis, Karlsruhe Institute of Technology (KIT), D-76128 Karlsruhe, Germany,
}
\email{lars.hientzsch@kit.edu}

\author[R. Scandone]{Raffaele Scandone}
\address{Universit\`a degli Studi di Napoli “Federico II”, Dipartimento di Matematica e Applicazioni “R. Caccioppoli”, Complesso Monte S. Angelo - Via Cintia, 80126 Napoli, Italy}
\email{raffaele.scandone@unina.it}

\author[R. Wegner]{Robert Wegner}
\address{
Institute for Analysis, Karlsruhe Institute of Technology (KIT), D-76128 Karlsruhe, Germany,
}
\email{robert.wegner@kit.edu}


\begin{document}

\begin{abstract}
	For a function $\psi: \R^d \rightarrow \C$, we define the \emph{phase-corrected gradient} $\nabla_{\!\circ} \psi = |\psi|^{-1} \conj{\psi}\,\nabla\psi$ and study the regularity of the map $\psi \mapsto \nabla_{\!\circ} \psi$ on Sobolev spaces. Since $\RE[\nabla_{\!\circ} \psi] = \nabla |\psi|$, the phase-corrected gradient fits naturally into the framework of nonlinear composition operators. Its imaginary part $\Lambda = \IM[\nabla_{\!\circ} \psi]$, on the other hand, plays a key role in quantum hydrodynamic models. We prove that $\psi \mapsto \nabla_{\!\circ} \psi$ is bounded and continuous from $W^{s,2}(\R^d)$ to $W^{s-1,2}(\R^d)$ for $1 \leq s < \frac32$, and we show that the threshold $s = \frac32$ is sharp. In particular, this answers an open question concerning the continuity of the map $\psi \mapsto |\psi|$ on Sobolev spaces. We also show that our results remain valid for functions that do not necessarily decay at infinity, a relevant setting in the context of Ginzburg--Landau theory. Finally, we extend our analysis to more general phase-corrected gradients of the form $f(|\psi|^{-1} \psi) \nabla \psi$, where $f: \mathbb{S}^1 \rightarrow \C$ is Lipschitz, and discuss implications for certain parabolic PDEs.
\end{abstract}

\maketitle

\tableofcontents

\section{Introduction}
\subsection{Setup of the problem and summary of main results}\label{sec:setup}
Given a complex-valued function $\psi$ on $\R^d$, $d\geq 1$, we investigate the regularity properties in Sobolev spaces of the nonlinear operator
\begin{equation} \label{eqn:phase_corrected}
	\psi\mapsto \nabla_{\!\circ} \psi=\frac{\conj{\psi}\,\nabla\psi}{|\psi|} ,
\end{equation}
which we call the \emph{phase-corrected gradient}. The map \eqref{eqn:phase_corrected} is certainly well-defined when $\psi$ vanishes nowhere: in this case, writing $\psi= a \eul^{\imu \Theta}$ with $a = |\psi|$, we have
\begin{equation} \label{eqn:polar-formal}
	\nabla_{\!\circ}\psi=\nabla a +\imu a \nabla \Theta .
\end{equation}
In particular, taking the real part of $\nabla_{\!\circ} \psi$, we formally obtain the map
\begin{equation}\label{eqn:grad-mod}
	\psi \mapsto \nabla |\psi|,
\end{equation}
which has been widely studied in the context of nonlinear composition operators (see Section \ref{sec:composition} for further details).

Beyond its independent analytical interest, the operator \eqref{eqn:phase_corrected} arises naturally in hydrodynamic models describing quantum fluids, as discussed in detail in Section \ref{sec:qhd}. We anticipate here the main features of this hydrodynamic interpretation, which also leads to the quantities at the core of our analysis. Given a nowhere vanishing function $\psi= a \eul^{\imu \Theta}$, its Madelung transform is defined by
\begin{equation} \label{eqn:formal_madelung}
	\mathcal{M}: \psi \mapsto \left(a=|\psi|, v=\nabla \Theta\right) .
\end{equation}
In quantum hydrodynamics, $\psi$ is a wave function, typically a solution to a nonlinear Schr\"odinger (NLS) equation, and the nonlinear operator $\mathcal{M}$ formally transforms the NLS equation into a compressible Euler-type system with a dispersive correction. In this framework, $a^2$ and $v$ are the mass density and the velocity field of the fluid, respectively. In view of \eqref{eqn:polar-formal}, the phase-corrected gradient encodes these quantities, namely
\begin{equation*}
	\nabla_{\!\circ}\psi=\nabla a+\imu\Lambda ,
\end{equation*}
where $\Lambda := a v$ is the velocity field weighted by the amplitude. The pair $(\nabla a,\Lambda)$ provides natural hydrodynamic variables. Indeed, the momentum density is given by $J=\IM\left[\conj{\psi}\nabla\psi\right]=a\Lambda$, and the gradient energy density $|\nabla\psi|^2$ of NLS splits as $|\nabla\psi|^2=|\nabla a|^2+|\Lambda|^2$, namely as the sum of the quantum (Bohm) energy density $|\nabla a|^2$ and the kinetic energy density $|\Lambda|^2=a^2|v|^2$ of the fluid. If $\psi$ is allowed to vanish, $\nabla a$ is still well-defined as soon as $\psi$ has one locally integrable weak derivative, whereas this is a priori not the case for $\Lambda$, since the phase $\Theta$, and hence the velocity $v$, are not defined on the vacuum set $\{\psi = 0\}$. Following \cite{AntonelliMarcati2026}, we show in Section \ref{sec:realization} that, under the same minimal regularity assumption $\psi \in W^{1,1}_{\loc}(\R^d)$, the phase-corrected gradient, and thus also the hydrodynamic variables $(\nabla a, \Lambda)$, are well-defined also in the presence of vacuum.

Our aim is then to develop a systematic approach yielding sharp boundedness and continuity results for the phase-corrected gradient, which at the same time provides a precise description of the behavior of its components $\nabla a$ and $\Lambda$. More precisely, we address the following problems.

\begin{problem} \label{prb:1}
	Identify the optimal range of parameters $d\in\N^+$, $s \geq 1$, and $p\in[1,\infty]$,\linebreak such that the nonlinear operator $\psi \mapsto \nabla_{\!\circ} \psi$ is bounded from $W^{s,p}(\R^d; \C)$ into $W^{s-1,p}(\R^d; \C^d)$.
\end{problem}
\begin{problem} \label{prb:2}
	Identify the optimal range of parameters $d\in\N^+$, $s \geq 1$, and $p\in[1,\infty]$,\linebreak such that the nonlinear operator $\psi \mapsto \nabla_{\!\circ} \psi$ is continuous from $W^{s,p}(\R^d; \C)$ into $W^{s-1,p}(\R^d; \C^d)$.
\end{problem}

Establishing mapping properties of \eqref{eqn:phase_corrected} poses nontrivial technical issues: when $\psi$ vanishes on some subset of $\R^d$ (a vacuum region, in hydrodynamic terms), one cannot appeal to general results on nonlinear operators on Sobolev spaces (see the discussion in Section \ref{sec:idea}).

Our main results resolve Problems \ref{prb:1} and \ref{prb:2} conclusively in the case $p = 2$, and our method is not naturally applicable to the cases $p \neq 2$.

\begin{theorem}\label{thm:main-pre1} Let $d \in \N^+$.
	\begin{itemize}
		\item[(i)] Let $s \in [1, \frac32)$. There exists $C > 0$ such that for all $\psi \in W^{s,2}(\R^d; \C)$ we have $\nabla_{\!\circ} \psi \in W^{s-1,2}(\R^d; \C^d)$ and
		\begin{align*}
			\|\nabla_{\!\circ} \psi\|_{W^{s-1,2}} &\leq C \|\nabla \psi\|_{W^{s-1,2}}. 
		\end{align*}
		Furthermore, if $\psi_n \rightarrow \psi$ in $W^{s,2}(\R^d)$, then $\nabla_{\!\circ} \psi_n \rightarrow \nabla_{\!\circ} \psi$ in $W^{s-1,2}(\R^d)$.
		\item[(ii)] Let $s \geq \frac32$. There exists $\psi \in C^\infty(\R^d; \C)$ with $\nabla_{\!\circ} \psi \in C^\infty(\R^d; \C^d)$ such that $\psi \in W^{s,2}(\R^d)$ but $\nabla a, \Lambda \not\in W^{s-1,2}(\R^d)$.
	\end{itemize}
\end{theorem}

As a consequence, we also answer an open problem concerning the continuity of the map $\psi \mapsto |\psi|$ on Sobolev spaces.

\begin{corollary} \label{cor:modulus}
	The map $\psi\mapsto|\psi|$ is continuous on $W^{s,2}(\R^d; \C)$ for all $s\in[1,\frac32)$.
\end{corollary}

The boundedness of the map $\psi\mapsto |\psi|$ on $W^{s,2}(\R^d)$, for $1\leq s<\frac32$, is due to \cite{BourdaudMeyer} (see also \cite{Oswald} in the real-valued case). Continuity for $s=1$ is a standard result (see e.g.\ \cite[Proposition 1.2]{AntonelliMarcati2026}). By contrast, continuity for $s\in(1,\frac32)$, even in the real-valued case, has been formulated as an open problem in \cite{BourdaudSickel,MusinaNazarov} and is not addressed in \cite{BourdaudMeyer,Oswald}. Corollary \ref{cor:modulus} is particularly relevant for applications to quantum fluids, as discussed in Section \ref{sec:qhd}.

Our method leverages local regularity properties of complex-valued functions and is robust enough to treat wave functions $\psi\in W_{\loc}^{s,2}(\R^d)$ with $\nabla\psi\in \dot{W}^{s-1,2}(\R^d)$ (see Section \ref{sec:full} for the definition of this space). A prototypical example is given by functions of finite Ginzburg--Landau energy, e.g.\ solutions to Gross--Pitaevskii type equations with non-vanishing conditions at infinity. The associated energy space reads
\begin{equation*}
	\mathbf{E}(\R^d)=\left\{\psi\in L_{\loc}^2(\R^d;\C): \nabla\psi\in L^2(\R^d;\C^d) \text{ and } |\psi|^2-1\in L^2(\R^d;\R)\right\}.
\end{equation*}
Equipped with the distance function
\begin{equation*}
	d_{\mathbf{E}}(\psi_1,\psi_2)=\|\psi_1-\psi_2\|_{\mathcal{Z}^1+W^{1,2}}+\||\psi_1|^2-|\psi_2|^2\|_{L^2},
\end{equation*}
where $\mathcal{Z}^1$ denotes the Zhidkov space, $(\mathbf{E}(\R^d),d_{\mathbf{E}})$ is a complete metric space \cite{Gerard2006}.

\begin{theorem}\label{thm:main-pre3}
	Let $d \in \N^+$ and $s \in [1, \frac32)$. There exists $C > 0$ such that for all $\psi \in \mathbf{E}(\R^d)$ with $\nabla\psi\in {W}^{s-1,2}(\R^d; \C^d)$ we have $\nabla_{\!\circ} \psi \in W^{s-1,2}(\R^d; \C^d)$ and
	\begin{align*}
		\|\nabla_{\!\circ} \psi\|_{W^{s-1,2}} &\leq C \|\nabla \psi\|_{W^{s-1,2}}.
	\end{align*}
	Furthermore,
	\begin{align*}
		\begin{aligned}
			\psi_n &\rightarrow \psi \text{ in } \mathbf{E}(\R^d) \\
			\nabla \psi_n &\rightarrow \nabla \psi \text{ in } W^{s-1,2}(\R^d)
		\end{aligned}
		\quad&
		\Longrightarrow\quad \nabla_{\!\circ} \psi_n \rightarrow \nabla_{\!\circ} \psi \text{ in } W^{s-1,2}(\R^d) .
	\end{align*}
\end{theorem}

The example of functions of finite Ginzburg--Landau energy, as well as further applications in PDE, strongly motivates the development of a theory that is not restricted to wave functions decaying at infinity, such as those in $L^p(\R^d)$, but rather allows for very general behavior of $\psi$ at infinity (concerning nonlinear Schrödinger equations, for example, see \cite{Gerard2006,Gerard2008,Klaus2023,Klaus2026,KlausKunstmann2022,KotaniWangXuZhang2026}).

Theorems \ref{thm:main-pre1} and \ref{thm:main-pre3} are consequences of Theorem \ref{thm:main} below, which is the full version of our main results, together with Proposition \ref{prop:2} for the case $s=1$. 
In Theorem \ref{thm:main}, the convergence $\psi_n \rightarrow \psi$ is only required in $L^2_{\loc}(\R^d)$, rather than with respect to the full topology of the energy space $\mathbf{E}(\R^d)$.
By the fractional Poincaré inequality this is equivalent to convergence in $W^{s,2}_{\loc}(\R^d)$, given that $\nabla\psi_n \rightarrow \nabla\psi$ in $W^{s-1,2}(\R^d)$ is already assumed. For notational convenience, we use the latter formulation throughout the paper.
Furthermore, Theorem \ref{thm:main} includes refined statements (points (iv)--(vi)) at the threshold regularity $s=\frac32$, which provide a sharp characterization in Besov spaces.

\subsection{Realization of $\nabla_{\!\circ}\psi$} \label{sec:realization} We now provide the rigorous realization of $\nabla_{\!\circ}\psi$, which is a priori not even well-defined on the vacuum set due to the division by $|\psi|$. As $\psi$ is only of Sobolev regularity and not necessarily continuous, the vacuum set $\mathcal{V}$ is defined as the set of Lebesgue points at which $\psi$ vanishes, i.e.\
\begin{equation*}
	\mathcal{V}=\left\{ x\in \R^d: \lim_{r\rightarrow 0}\frac{1}{|B_r(x)|}\int_{B_r(x)}|\psi(y)| \dd y=0\right\}.
\end{equation*}
Following Antonelli and Marcati \cite{AntonelliMarcati2026}, we introduce the set of \emph{polar factors} of $\psi$ as
\begin{align*}
	\mathcal{P}(\psi) &= \{P\in L^{\infty}(\R^d;\C) : \psi = a P \text{ and } |P|=1\text{ almost everywhere on }\R^d\} .
\end{align*}
Any $P \in \mathcal{P}(\psi)$ is uniquely determined by
\begin{equation*}
	P(x)=\frac{\psi(x)}{|\psi(x)|}
\end{equation*}
provided that $\psi(x)\neq 0$, while it can be an arbitrary phase factor $\eul^{\imu\Theta}$, $\Theta\in\R$, on the vacuum. In other words, the polar factor is unique with respect to the measure $a^2 \dd x$. The construction of polar factors is local, i.e.\ not sensitive to the far field behavior of $\psi$, and can even be carried out for $\psi\in L_{\loc}^1(\R^d)$. 

For the realization of $\nabla_{\!\circ}\psi$ we require, as anticipated above, that $\psi \in W^{1,1}_{\loc}(\R^d)$, which we assume henceforth. Under this assumption, the set $\mathcal{P}(\psi)$ can also be characterized by a variational argument \cite[Section 1.5]{AntonelliMarcati2026}. We recall that $\nabla\psi=0$ almost everywhere on $\mathcal{V}$ (see e.g.\ \cite[Theorem 6.19]{LL01}). This crucial fact suggests that $\nabla_{\!\circ} \psi$ may be well-behaved even in the vicinity of $\mathcal{V}$. By writing
\begin{align*}
	\nabla_{\!\circ}\psi=\conj{P}\nabla\psi,
\end{align*}
where $P\in\mathcal{P}(\psi)$ is an arbitrary polar factor of $\psi$, we realize the phase-corrected gradient \eqref{eqn:phase_corrected} as a function in $L^1_{\loc}(\R^d; \C^d)$ which does not depend on the choice of $P$. It admits the explicit representation
\begin{equation}
	\nabla_{\!\circ}\psi = \begin{cases}
		\frac{\conj{\psi}\,\nabla\psi}{|\psi|}&\text{on }\mathcal{V}^c\\
		0&\text{a.e.\ on }\mathcal{V}.
	\end{cases}
\end{equation}
For $\psi\in W^{1,1}_{\loc}(\R^d; \C)$ we have $a \in W^{1,1}_{\loc}(\R^d; \R)$ with $|\nabla a|\leq |\nabla \psi|$ pointwise almost everywhere \cite[Theorem 4.4]{EvansGariepy}, and a direct computation (see e.g.\ the argument in \cite{AntonelliMarcati2026}) shows that
\begin{equation*}
	\nabla|\psi|=\RE\nabla_{\!\circ}\psi\in L^1_{\loc}(\R^d; \R^d) .
\end{equation*}
Accordingly, we define the vector field $\Lambda: \R^d \rightarrow \R^d$ by
\begin{equation*}
	\Lambda=\IM\nabla_{\!\circ}\psi\in L^1_{\loc}(\R^d; \R^d)
\end{equation*}
so that
\begin{equation}\label{eqn:splitting}
	\nabla_{\!\circ}\psi=\nabla a+\imu\Lambda
\end{equation}
holds in $L^1_{\loc}(\R^d)$. By construction,
\begin{align*}
	\nabla a(x) &= \RE\left[ \frac{\conj{\psi(x)} \nabla \psi(x)}{|\psi(x)|}\right]
	& \Lambda(x) &= \IM\left[ \frac{\conj{\psi(x)} \nabla \psi(x)}{|\psi(x)|}\right]
\end{align*}
if $\psi(x)\neq 0$, while
\begin{equation*}
	\nabla a=\Lambda = 0 \text{ almost everywhere on } \mathcal{V}.
\end{equation*}
Given a sequence $(\psi_n)_{n \in \N}$, we define $(a_n, P_n, \Lambda_n)_{n \in \N}$ analogously.
Owing to identity \eqref{eqn:splitting}, the regularity properties of the phase-corrected gradient are equivalent to those of the map
\begin{equation}\label{eqn:psitohydro}
	\psi\mapsto (\nabla a,\Lambda).
\end{equation}

\subsection{Threshold cases}
As mentioned above, for $s=1$ the boundedness and continuity of the phase-corrected gradient are known \cite[Proposition 1.2]{AntonelliMarcati2026} and rely on a simple identity. We state this result in the general setting of functions with arbitrary behavior at infinity.

\begin{proposition}[{analogous to \cite[Proposition 1.2]{AntonelliMarcati2026}}] \label{prop:2}
	Let $d \in \N^+$ and let $\psi\in L^2_{\loc}(\R^d;\C)$ with $\nabla\psi\in L^2(\R^d;\C^d)$. Then $\nabla_{\!\circ}\psi\in L^2(\R^d;\C^d)$ and
	\begin{equation} \label{eqn:L2}
		\|\nabla \psi\|_{L^2}^2 = \|\nabla_{\!\circ} \psi\|_{L^2}^2 = \|\nabla a\|_{L^2}^2+ \|\Lambda\|_{L^2}^2.
	\end{equation}
	Moreover,
	\begin{align*}
		\begin{aligned}
			\psi_n &\rightarrow \psi \text{ in } L^2_{\loc}(\R^d) \\
			\nabla \psi_n &\rightarrow \nabla \psi \text{ in } L^2(\R^d)
		\end{aligned}
		\quad&
		\Longrightarrow\quad \nabla_{\!\circ} \psi_n \rightarrow \nabla_{\!\circ} \psi \text{ in } L^2(\R^d) .
	\end{align*}
\end{proposition}

The statement in \cite{AntonelliMarcati2026} (see also \cite{AntonelliMarcati}, and \cite{AHMZ,H} for extensions to nontrivial far field) concerns only $W^{1,2}(\R^d)$, but the argument in our setting is analogous.
The polarization identity \eqref{eqn:L2} is a direct consequence of $|P|=1$ and the orthogonal decomposition of $\nabla_{\!\circ}\psi$ into $\nabla a$ and $\Lambda$. It immediately yields boundedness and the equality of norms $\|\nabla_{\!\circ}\psi\|_{L^2}=\|\nabla\psi\|_{L^2}$. Furthermore, as in the proof of \cite[Proposition 1.2]{AntonelliMarcati2026}, it is the key ingredient to obtain continuity.

We establish an analogue of \eqref{eqn:L2} in the fractional range $s\in(1,\frac32)$, namely the polarization identity in Proposition \ref{prop:formula} below. It plays the same twofold role: it underlies the equivalence of seminorms in Theorem \ref{thm:main}(ii), and it is the key ingredient in the proof of continuity.

Concerning the threshold $s=\frac32$, the following elementary examples already show its critical role: boundedness fails separately for $\nabla a$ and for $\Lambda$.

\begin{example} \label{example:threshold}
	Let $d=1$, and let $\chi\in C_c^\infty(\R;[0,1])$ with $\chi=1$ on $[-1,1]$ and $\supp\chi\subset[-2,2]$.
	\begin{enumerate}[(a)]
		\item Let $\psi(x)=\chi(x)\,x$. Then $\psi\in C_c^\infty(\R;\R)$, meaning $\Lambda=0$, while $\nabla a(x)=\sign(x)$ for $x\in(-1,1)$. Hence $\nabla a\in W^{s-1,2}(\R)$ only for $s<\frac32$, although $\psi\in W^{s,2}(\R)$ for every $s\geq1$.
		\item Let $\epsilon\in(0,1)$ and set
		\begin{equation*}
			\psi_\epsilon(x)=\chi(x)\big((x-\epsilon)+\imu(x+\epsilon)\big) .
		\end{equation*}
		The family $(\psi_\epsilon)_{\epsilon\in(0,1)}$ is bounded in $W^{s,2}(\R)$ for every $s\geq1$. On $(-1,1)$, the function $\psi_\epsilon$ does not vanish and a direct computation gives
		\begin{align*}
			\Lambda_\epsilon(x) &= -\frac{2\epsilon}{\sqrt{(x-\epsilon)^2+(x+\epsilon)^2}}=-g\Big(\frac{x}{\epsilon}\Big) \quad\text{ where }\quad
            g(y) = \frac{\sqrt{2}}{\sqrt{1+y^2}} .
		\end{align*}
		Now take $s\in(\frac32,2)$. Restricting the integral defining $[\Lambda_\epsilon]_{W^{s-1,2}}$ to $(-1,1)^2$ and rescaling, we obtain
		\begin{equation*}
		[\Lambda_\epsilon]_{W^{s-1,2}}^2
			\geq \epsilon^{3-2s}\int_{-1}^{1}\!\int_{-1}^{1}\frac{|g(x)-g(y)|^2}{|x-y|^{2s-1}}\dd x\dd y ,
		\end{equation*}
		where the last integral is positive since $g$ is not constant. Therefore $\|\Lambda_\epsilon\|_{W^{s-1,2}}$ diverges as $\epsilon\to0$, i.e.~the map $\psi\mapsto\Lambda$ is unbounded from $W^{s,2}(\R)$ into $W^{s-1,2}(\R)$. In fact, for any $s > \frac32$ the sequence $(\psi_\epsilon)_{\epsilon \in (0, 1)}$ is bounded in $W^{s,2}(\R)$, and since $W^{s-1,2}(\R)\hookrightarrow W^{s'-1,2}(\R)$ for $s\geq s'$ the same conclusion holds.  Note that in the limit $\epsilon \searrow 0$ the phase of $\psi_\epsilon$ attains a jump at the origin.
	\end{enumerate}
	Both examples can be transferred to any dimension $d\geq2$ by multiplication with a cutoff function in the remaining variables $(x_2,\dots,x_d)$.
\end{example}

Example \ref{example:threshold}(a) provides a smooth, compactly supported function $\psi$ such that $\nabla a\notin W^{s-1,2}(\R)$ for all $s\geq\frac32$, whereas (b) only yields the unboundedness of $\psi\mapsto\Lambda$ for $s>\frac32$. At $s=\frac32$, the lower bound is scale invariant and does not blow up. Moreover, in (a) the phase-corrected gradient $\nabla_{\!\circ}\psi=\sign(x)$ for $x \in (-1, 1)$ is discontinuous at the origin. The stronger statement of Theorem \ref{thm:main-pre1}(ii), where both $\psi$ and $\nabla_{\!\circ}\psi$ are smooth and the endpoint $s=\frac32$ is included also for $\Lambda$, requires a more refined construction (see Lemma \ref{lem:32} below).

\subsection{Full main results} \label{sec:full}

We now set up the functional framework needed to state our main result in its full generality. We first recall several equivalent ways to express the seminorm $[F]_{W^{\theta,2}}$, $\theta \in (0, 1)$, for vector fields $F: \R^d \rightarrow \C^d$.

\begin{definition}
	Let $\theta \in (0,1)$. For $F: \R^d \rightarrow \C^d$ we define the seminorm
	\begin{align} \label{eqn:50}
		[F]_{W^{\theta,2}} &= \left( \int_{\R^d} |h|^{-d-2\theta} \int_{\R^d} |\Delta_h F(x)|^2 \dd x \dd h \right)^{\frac12},
	\end{align}
where $\Delta_h F(x) = F(x) - F(x-h)$, and the spaces
	\begin{align*}
		W^{\theta,2}(\R^d;\C^d) &= \left\{F \in L^2(\R^d;\C^d) \, : \, [F]_{W^{\theta,2}}<\infty\right\}, \\
		\dot{W}^{\theta,2}(\R^d;\C^d) &= \left\{F \in L^{\frac{2d}{d-2\theta}}(\R^d;\C^d) \, : \, [F]_{W^{\theta,2}}<\infty\right\} \quad \text{ for } \theta < {\textstyle{\frac{d}{2}}}.
	\end{align*}
\end{definition}
\begin{lemma} \label{lem:characterizations}
	Let $d\in \N^+$ and $\theta\in (0,1)$.
	\begin{enumerate}[(i)]
		\item There exists $C_{d,\theta}>0$ such that
		\begin{equation*}
			[F]_{W^{\theta,2}}^2 = C_{d,\theta} \|F\|_{\dot{H}^\theta}^2
		\end{equation*}
		for every $F \in C_c^\infty(\R^d; \C^d)$, where $\|F\|_{\dot{H}^{\theta}}^2=\||\cdot|^{\theta}\hat{F}\|^2_{L^2}$.
		\item Define the alternative seminorm
		\begin{align}
			[F]_{X^{\theta,2}} &= \left( \int_{\R^d} |h|^{-d-2-2\theta} \int_{\R^d} |\Delta_h F(x) \cdot h|^2 \dd x \dd h \right)^{\frac12}.
		\end{align}
		Then the seminorms $[\,\cdot\,]_{X^{\theta,2}}$ and $[\,\cdot\,]_{W^{\theta,2}}$ are equivalent on $C_c^\infty(\R^d)$.
	\end{enumerate}
	In particular, the seminorms $\|\cdot\|_{\dot{H}^\theta}$, $[\,\cdot\,]_{W^{\theta,2}}$, and $[\,\cdot\,]_{X^{\theta,2}}$ are all equivalent on $W^{\theta,2}(\R^d)$, and if $\theta < \frac{d}{2}$ also on $\dot{W}^{\theta,2}(\R^d)$.
\end{lemma}

\begin{proof}
	\begin{enumerate}[(i)]
		\item If $F \in C_c^\infty(\R^d; \C^d)$, then a calculation with the Fourier transform
		\begin{align*}
			\ha{u}(\xi) &= \frac{1}{(2 \pi)^{\frac{d}{2}}} \int_{\R^d} \eul^{- \imu x \cdot \xi} u(x) \dd x
		\end{align*}
		shows that
		\begin{align*}
			[F]_{W^{\theta,2}}^2 &= C_{d,\theta} \int_{\R^d} |\xi|^{2\theta} |\ha{F}(\xi)|^2 \dd \xi
			\quad\text{ with }\quad C_{d,\theta} = \frac{2^{1-2\theta} \pi^{\frac{d}{2}} \Gamma(1-\theta)}{\theta \Gamma(\frac{d}{2}+\theta)} .
		\end{align*}
		\item By \cite[Theorem 1]{ScottMengesha2019} the seminorms $[\,\cdot\,]_{X^{\theta,2}}$ and $[\,\cdot\,]_{W^{\theta,2}}$ are equivalent on $C_c^\infty(\R^d)$.
		In our case $p = 2$, this can also be seen by a calculation with the Fourier transform (see \cite{QiangKun2011,ScottMengesha2019}), which for $F \in C_c^\infty(\R^d)$ yields
		\begin{align} \label{eqn:X-H}
			[F]_{X^{\theta,2}}^2 &= \frac{C_{d,\theta}}{d+2\theta}\left( \|F\|_{\dot{H}^\theta}^2 + 2 \theta \|\Pi F\|_{\dot{H}^\theta}^2 \right) ,
		\end{align}
		where $\ha{\Pi F}(\xi) = \frac{\xi \otimes \xi}{|\xi|^2} \ha{F}(\xi)$ is the projection onto gradient fields in the Helmholtz decomposition.
		In particular,
		\begin{align*}
			[F]_{X^{\theta,2}}^2 &\leq [F]_{W^{\theta,2}}^2 \leq (d + 2 \theta) [F]_{X^{\theta,2}}^2.\qedhere
		\end{align*}
	\end{enumerate}
\end{proof}

Our arguments naturally work with the alternative seminorm $[\,\cdot\,]_{X^{\theta,2}}$. If $F = \nabla f$ for some $f: \R^d \rightarrow \C$, then we have the following representations (see \eqref{eqn:9}--\eqref{eqn:9-alt} below):
\begin{align*}
	[\nabla f]_{X^{\theta,2}}^2 &= (2 + 2 \theta) \int_{\R^d} |h|^{-d-2-2\theta} \int_{\R^d} |R_h f(x)|^2 \dd x \dd h
	\\\nonumber &= \frac{2 + 2 \theta}{\theta} \int_{\R^d} |h|^{-d-2-2\theta} \int_{\R^d} \RE\left[R_h f(x) \conj{R_{-h} f(x-h)}\right] \dd x \dd h ,
	\\\nonumber &\quad\text{ where }\quad R_h f(x) = \Delta_h f(x) - \nabla f(x) \cdot h .
\end{align*}

Throughout, we consider $\psi \in W^{s,2}_{\loc}(\R^d; \C)$, complemented with one of the two assumptions
\begin{equation*}
	\nabla \psi \in \dot{W}^{s-1,2}(\R^d; \C^d) \qquad \text{ or } \qquad \nabla \psi \in W^{s-1,2}(\R^d; \C^d) .
\end{equation*}
Both assumptions only concern $\nabla\psi$, and hence allow for nontrivial behavior of $\psi$ at infinity. The second one is more restrictive since $W^{s-1,2}(\R^d)\subset\dot W^{s-1,2}(\R^d)$ by the Sobolev embedding. As we shall see, the first assumption is sufficient for boundedness, while we require the second assumption for the equivalence of seminorms and for continuity. 
We now state the full version of our main result: parts (i)--(iii) concern the range $s\in(1,\frac32)$, and parts (iv)--(vi) the range $s\geq\frac32$.

\begin{theorem} \label{thm:main}
	Let $d \in \N^+$ and $s > 1$. The following holds.
	\begin{enumerate}[(i)]
		\item Take $s < \frac32$ and let $\psi \in W^{s,2}_{\loc}(\R^d; \C)$ with $\nabla \psi \in \dot{W}^{s-1,2}(\R^d; \C^d)$.
		Then $\nabla_{\!\circ} \psi \in \dot{W}^{s-1,2}(\R^d; \C^d)$ and
		\begin{align}
			\label{eqn:est-0} (3-2s)^2 [\nabla a]_{X^{s-1,2}}^2 + (3-2s) [\Lambda]_{X^{s-1,2}}^2 &\leq 8 [\nabla \psi]_{X^{s-1,2}}^2 .
		\end{align}
		Moreover, the blow-up rates of the constants in \eqref{eqn:est-0} are sharp: there exists a sequence $(\psi_n)_{n \geq 4} \subset W^{\frac32-\frac{1}{n},2}(\R^d)$ with $[\nabla \psi_n]_{X^{\frac12-\frac{1}{n},2}}^2 \lesssim 1$ such that $[\nabla a_n]_{X^{\frac12-\frac{1}{n},2}}^2 \gtrsim n^2$ and $[\Lambda_n]_{X^{\frac12-\frac{1}{n},2}}^2 \gtrsim n$.
		
		\item Take $s < \frac32$ and let $\psi \in W^{1,2}_{\loc}(\R^d; \C)$ with $\nabla \psi \in L^2(\R^d; \C^d)$.
		Then $\nabla \psi \in W^{s-1,2}(\R^d; \C^d)$ if and only if $\nabla_{\!\circ} \psi, P \Lambda \in W^{s-1,2}(\R^d; \C^d)$.
		In this case, we have the equivalence
		\begin{align}
			\label{eqn:0-equiv-main}
			\frac{(3-2s)^2}{32} [\nabla \psi]_{X^{s-1,2}}^2 &\leq  [\nabla_{\!\circ} \psi]_{X^{s-1,2}}^2
			+ \frac38 [P \Lambda]_{X^{s-1,2}}^2
			\leq \frac{32}{(3-2s)^2} [\nabla \psi]_{X^{s-1,2}}^2 .
		\end{align}
		Lastly, there exists a sequence $(\psi_n)_{n \in \N} \subset W^{s,2}(\R^d)$ such that 
        \begin{equation}
		    [\nabla_{\!\circ} \psi_n]_{X^{s-1,2}}^2 + [P_n \nabla a_n]_{X^{s-1,2}}^2 \lesssim 1,
		\end{equation}
        while $[P_n \Lambda_n]_{X^{s-1,2}}$ and hence also $[\nabla \psi_n]_{X^{s-1,2}}$ are unbounded.
		\item Take $s < \frac32$. We have 
		\begin{align*}
			\begin{aligned}
				\psi_n &\rightarrow \psi \text{ in } W^{s,2}_{\loc}(\R^d) \\
				\nabla \psi_n &\rightarrow \nabla \psi \text{ in } W^{s-1,2}(\R^d)
			\end{aligned}
			\quad&
			\Longrightarrow\quad \nabla_{\!\circ} \psi_n \rightarrow \nabla_{\!\circ} \psi \text{ in } W^{s-1,2}(\R^d) .
		\end{align*}
		
		\item Take $s = \frac32$ and note the definition of the Besov spaces in \eqref{eqn:def-besov-1}--\eqref{eqn:def-besov-2} below. If $\psi \in B^{\frac32}_{2,1}(\R^d; \C)$, then $a \in B^{\frac32}_{2,\infty}(\R^d; \R)$ and $\Lambda \in W^{\frac12,2}(\R^d; \R^d)$. Furthermore, there exists a constant $C = C(d) > 0$ such that
		\begin{align} \label{eqn:endpoint}
			\|a\|_{B^{\frac32}_{2,\infty}}^2 + [\Lambda]_{X^{\frac12,2}}^2 &\leq C \|\psi\|_{B^{\frac32}_{2,1}}^2 .
		\end{align}
		Lastly, there exists a sequence $(\psi_n)_{n \in \N}$ such that
		$\psi_n \rightarrow \psi$ in $B^{\frac32}_{2,1}(\R^d)$, while $\nabla a_n \not\rightarrow \nabla a$ in $B^{\frac12}_{2,\infty}(\R^d)$ and $\Lambda_n \not\rightarrow \Lambda$ in $W^{\frac12,2}(\R^d)$.
		
		\item Take $s = \frac32$ and $\epsilon \in (0, \frac12)$. If $\psi \in W^{\frac32+\epsilon,2}_{\loc}(\R^d; \C)$ with $\nabla \psi \in W^{\frac12+\epsilon,2}(\R^d; \C^d)$, then $\Lambda \in W^{\frac12,2}(\R^d; \R^d)$ and
		\begin{align} \label{eqn:endpoint-epsilon}
			[\Lambda]_{X^{\frac12,2}}^2 &\leq \frac43 [\nabla \psi]_{X^{\frac12,2}}^2 + \frac{7}{\epsilon} [\nabla \psi]_{X^{\frac12+\epsilon,2}} [\nabla \psi]_{X^{\frac12-\epsilon,2}} .
		\end{align}
		
		\item Take $s \geq \frac32$. There exists $\psi \in C^\infty(\R^d; \C)$ with $\nabla_{\!\circ} \psi \in C^\infty(\R^d; \C^d)$ such that $\psi \in W^{s,2}(\R^d)$ but $\nabla a, \Lambda \not\in W^{s-1,2}(\R^d)$. Furthermore, there exists $\psi \in C^\infty(\R^d; \C)$ with $\nabla_{\!\circ} \psi \in C^\infty(\R^d; \C^d)$ such that $\psi \in W^{2,2}(\R^d)$ but $\nabla a \not\in W^{\frac12,2}(\R^d)$ and $\Lambda \not\in B^{\frac12}_{2,1}(\R^d)$.
	\end{enumerate}
\end{theorem}

Theorem \ref{thm:main} provides the first fractional regularity result for $\Lambda$, whose analysis is significantly more challenging than that of $\nabla a$ since it is highly sensitive to phase oscillations.
For $\nabla a$ one can use a spherical averaging formula (see \eqref{eqn:a-reduction} below) to reduce the problem to the case of real-valued $\psi$. One cannot expect such a formula to be available for $\Lambda$, since it does not correspond to a composition operator and is not, in general, a conservative vector field.

Let us comment on the individual statements of Theorem \ref{thm:main}.
\begin{itemize}
	\item[(i)] Boundedness holds in the whole range $s\in(1,\frac32)$ under the sole condition $\nabla\psi\in\dot W^{s-1,2}(\R^d)$, which allows for very general behavior of $\psi$ at infinity. For $\nabla a$ we use the method of \cite{BourdaudMeyer}, where boundedness of $\psi \mapsto a$ on $W^{s,2}(\R^d)$ is proven (see also \cite{Oswald}), to extend the result to general conditions at infinity. As $s\to\frac32$, the constants blow up at the rates $(3-2s)^{-2}$ for $\nabla a$ and $(3-2s)^{-1}$ for $\Lambda$, and both rates are sharp: $\Lambda$ is less singular than $\nabla a$ near the threshold.
	\item[(ii)] This is the fractional counterpart of the equality $\|\nabla_{\!\circ}\psi\|_{L^2}=\|\nabla\psi\|_{L^2}$ from Proposition \ref{prop:2}, and provides a converse to (i): the regularity of $\psi$ can be recovered from that of the hydrodynamic variables $(\nabla a,\Lambda)$, provided that $P\Lambda$ is controlled as well. In contrast to the case $s=1$, the phase-corrected gradient alone does not suffice: the last statement of (ii) shows that the control of $P\Lambda$ cannot be dropped, even if $P\nabla a$ is controlled.
	\item[(iii)] Continuity is new for both $\nabla a$ and $\Lambda$. It holds in the whole range $s\in(1,\frac32)$ and only requires local convergence of $\psi_n$, in addition to the convergence of $\nabla\psi_n$ in $W^{s-1,2}(\R^d)$. In particular, it yields Corollary \ref{cor:modulus}.
	\item[(iv)] At the threshold $s=\frac32$, boundedness from $W^{\frac32,2}$ into $W^{\frac12,2}$ fails (see (vi)), but it can be partially recovered by slightly strengthening the assumption to $\psi\in B^{\frac32}_{2,1}(\R^d)$. In this case $a\in B^{\frac32}_{2,\infty}(\R^d)$. This is due to \cite[Theorem 1.1]{Moussai2018}, where endpoint Besov estimates are proven for a class of composition operators including the modulus. On the other hand, $\Lambda$ belongs to $W^{\frac12,2}(\R^d)=B^{\frac12}_{2,2}(\R^d)$, which is strictly smaller than the space $B^{\frac12}_{2,\infty}(\R^d)$ containing $\nabla a$. This result is new and demonstrates once more that $\Lambda$ has better regularity properties than $\nabla a$ at the endpoint $s  =\frac32$. In contrast to (iii), continuity fails at the endpoint. We assume the fully inhomogeneous setting for simplicity, although we expect the assumption $\nabla \psi \in B^{\frac12}_{2,1}(\R^d)$ to suffice.
	\item[(v)] This is a variant of (iv) for $\Lambda$, which requires higher regularity but imposes weaker assumptions on the behavior at infinity, since only $\nabla\psi$, and not $\psi$, is required to be globally square integrable.
	\item[(vi)] The first statement shows that the range $s<\frac32$ in (i) is sharp for both $\nabla a$ and $\Lambda$, even for smooth $\psi$ and $\nabla_{\!\circ}\psi$: the failure of boundedness is thus not caused by a local singularity of $\nabla_{\!\circ}\psi$. The second statement shows that (iv) and (v) are sharp as well: $\nabla a\in W^{\frac12,2}(\R^d)$ cannot be expected even for $\psi\in W^{2,2}(\R^d)$, so that no analogue of (v) holds for $\nabla a$, and the regularity $\Lambda\in W^{\frac12,2}(\R^d)$ cannot be improved to $\Lambda\in B^{\frac12}_{2,1}(\R^d)$.
\end{itemize}

\begin{remark}[Regularity of $\nabla a$ and $\Lambda$ for $1\leq s<\frac32$] Observe that $\nabla a$ cannot be more regular than $\nabla\psi$ in general: for instance, $\nabla a=\nabla\psi$ whenever $\psi$ is real-valued and positive. For $\Lambda$, on the other hand, a gain of regularity is not excluded a priori, for two reasons. First, $\Lambda$ turns out to be better behaved than $\nabla a$, see Theorem \ref{thm:main}(i) and comment (iv) above, as well as Theorem \ref{thm:main-extra} below. Second, in $\Lambda=av$ the vanishing amplitude might dampen the singularities of the velocity field at the vacuum. The following examples show that neither of the two effects leads to a gain, i.e.\ $\Lambda$ has in general the same regularity as $\nabla\psi$, both away from and at the vacuum. Let
\begin{align*}
\psi_1(x)&=\chi(x)\big(x+\tfrac12\big)\eul^{\imu|x|} & \psi_2(x)&=\chi(x)\,x\,\eul^{\imu\ln|x|}
\end{align*}
with $\chi$ as in Example \ref{example:threshold}. Then $\psi_1,\psi_2\in W^{s,2}(\R)$ exactly for $s<\frac32$, and
\begin{align*}
\Lambda_1(x)&=\chi(x)\big|x+\tfrac12\big|\sign(x) & \Lambda_2(x)&=\chi(x)\sign(x) .
\end{align*}
Both $\Lambda_j$, $j = 1, 2$, display a jump at the origin, so that $\Lambda_j\in W^{s-1,2}(\R)$ exactly for $s<\frac32$. In the first case, the jump is produced by the kink of the phase $|x|$ at a point where $\psi_1$ does not vanish, its only zero in $(-2,2)$ being $x=-\frac12$. In the second case, the origin is a vacuum point at which the velocity field $v_2(x)=\frac1x$ is singular: the vanishing amplitude $a_2(x)=\chi(x)|x|$ compensates the blow-up of $v_2$ only up to a jump of $\Lambda_2=a_2v_2$, so that the vacuum has no regularizing effect.
\end{remark}

We now state the polarization identity anticipated above.
\begin{proposition} \label{prop:formula}
	Let $d \in \N^+$ and $s \in (1, \frac32)$. Suppose $\psi \in W^{s,2}_{\loc}(\R^d; \C)$ with $\nabla \psi \in W^{s-1,2}(\R^d; \C^d)$. Then $\nabla a, \Lambda, P \Lambda \in W^{s-1,2}(\R^d)$ and the following polarization identity holds:
	\begin{align} \label{eqn:formula-main}
		& (3-2s) [\nabla a]_{X^{s-1,2}}^2 + (2s-1) [\Lambda]_{X^{s-1,2}}^2 - (3-2s) [\nabla \psi]_{X^{s-1,2}}^2
		\\\nonumber &= \int_{\R^d} \!|h|^{-d-2s}\! \int_{\R^d} 2 (2s-1) \RE\Big[ \conj{P(x)} \Big( R_h \psi(x) \Delta_h \nabla a(x) \cdot h
		\\\nonumber &\hspace{50mm} - R_h a(x) \Delta_h \nabla \psi(x) \cdot h \Big) \Big] \dd x \dd h
		\\\nonumber &- \int_{\R^d} \!|h|^{-d-2s}\! \int_{\R^d} \IM\Big[ \conj{P(x)} R_h \psi(x) \Big( \big(8(s-1)-2(3-2s)P(x)\conj{P(x-h)}\big) \Delta_h \Lambda(x) \cdot h
		\\\nonumber &\hspace{34.5mm} + 2(3-2s)\conj{P(x-h)} \Delta_h(P\Lambda)(x) \cdot h \Big) \Big] \dd x \dd h .
	\end{align}
\end{proposition}

\begin{remark}[Concerning the real-valued case] \label{rmk:10}
	Consider the case $\psi: \R^d \rightarrow \R$, where $\Lambda = 0$. We write $u = \psi$ in order to distinguish it from the complex-valued case. For $u\in W^{s,2}(\R^d;\R)$, the estimate for $\nabla a$ in (i) is due to \cite{BourdaudMeyer,Oswald}, although we provide explicit constants. If $u$ changes sign, then the strict lower bound $[\nabla |u|]_{W^{s-1,2}(\R^d)} > [\nabla u]_{W^{s-1,2}(\R^d)}$ is proven to hold for $s \in (1, \frac32)$ in \cite[Theorem 1]{MusinaNazarov} (see also \cite{IgnatevNazarovNichitenkoTursunbaev2026}).
	This relies on the following polarization identity (see \cite[Equation (4)]{MusinaNazarov}):
	\begin{align*}
		& [\nabla |u|]_{W^{s-1,2}}^2 - [\nabla u]_{W^{s-1,2}}^2
		\\ &= - 16 s (d + 2 s - 2) \int_{\R^d} |h|^{-d-2s} \int_{\R^d} \frac{u(x) + |u(x)|}{2} \frac{u(x-h) - |u(x-h)|}{2} \dd x \dd h .
	\end{align*}
	Here the right-hand side has a sign, which immediately yields the lower bound. By contrast, the polarization identity \eqref{eqn:formula-main}, which is central to our results, is formulated in terms of the seminorm $[\,\cdot\,]_{X^{s-1,2}}$ and has a right-hand side without definite sign.
\end{remark}

\subsubsection*{Organization of the proofs of Theorem \ref{thm:main} and Proposition \ref{prop:formula}}
For the convenience of the reader, we indicate where the different parts of Theorem \ref{thm:main} are proved. Theorem \ref{thm:1} and Corollary \ref{cor:1} below imply (i) and (ii), while the sharpness of the blow-up rates of the constants is verified in Lemma \ref{lem:sharp}. Theorem \ref{thm:continuity} corresponds to (iii); Proposition \ref{prop:besov}, Corollary \ref{cor:2} and Lemma \ref{lem:endpoint-cont} to (iv); Lemma \ref{lem:endpoint-epsilon} to (v); and Lemmas \ref{lem:32} and \ref{lem:32-besov} to (vi). Finally, Proposition \ref{prop:formula} is part of Theorem \ref{thm:1}.

\subsubsection*{Additional results for $s \in [\frac32, 2]$}

Since the integrals in \eqref{eqn:formula-main} depend analytically on $s$, it is natural to ask whether the identity can be extended to the range $s \in [\frac32, 2)$, provided all integrals exist. Indeed, this is the case, and in fact several results of Theorem \ref{thm:main} can be extended to this range. Using the Bourgain--Brezis--Mironescu formula \cite{BourgainBrezisMironescu2001}, it is even possible to obtain a statement at $s = 2$. The following theorem is a consequence of Theorems \ref{thm:1} and \ref{thm:continuity}, as well as Corollaries \ref{cor:1-extra} and \ref{cor:2-extra}.

\begin{theorem} \label{thm:main-extra}
    Let $d \in \N^+$ and $s \in [\frac32, 2]$. The following holds.
	\begin{enumerate}[(i)]
		\item Take $s < 2$ and let $\psi \in W^{s,2}_{\loc}(\R^d; \C)$ with $\nabla \psi \in W^{s-1,2}(\R^d; \C^d)$. Assume in addition that $a \in W^{s,2}_{\loc}(\R^d)$ and $\nabla a \in W^{s-1,2}(\R^d)$.
		Then we have $\Lambda \in W^{s-1,2}(\R^d)$ and
		\begin{align} \label{eqn:est-0-extra}
			[\Lambda]_{X^{s-1,2}}^2 &\leq (4 s - 6) [\nabla a]_{X^{s-1,2}}^2 + 4 [\nabla a]_{X^{s-1,2}} [\nabla \psi]_{X^{s-1,2}} + 2 [\nabla \psi]_{X^{s-1,2}}^2 .
		\end{align}
        Furthermore, \eqref{eqn:formula-main} still holds.
        
		\item Take $s = \frac32$. We have 
		\begin{align*}
			\begin{aligned}
				\psi_n &\rightarrow \psi \text{ in } W^{\frac32,2}_{\loc}(\R^d)
				\\ \nabla \psi_n &\rightarrow \nabla \psi \text{ in } W^{\frac12,2}(\R^d)
                \\ &\sup_{n \in \N}\;[\nabla a_n]_{X^{\frac12,2}} < \infty
			\end{aligned}
			\quad&
			\Longrightarrow\quad \Lambda_n \rightarrow \Lambda \text{ in } W^{\frac12,2}(\R^d) .
		\end{align*}

        \item Take $s \in (\frac32, 2)$. We have 
		\begin{align*}
			\begin{aligned}
				\psi_n &\rightarrow \psi \text{ in } W^{s,2}_{\loc}(\R^d)
				\\ \nabla \psi_n &\rightarrow \nabla \psi \text{ in } W^{s-1,2}(\R^d)
                \\ \nabla a_n &\rightarrow \nabla a \text{ in } W^{s-1,2}(\R^d)
			\end{aligned}
			\quad&
			\Longrightarrow\quad \Lambda_n \rightarrow \Lambda \text{ in } W^{s-1,2}(\R^d) .
		\end{align*}
        
		\item Take $s = 2$ and let $\psi \in W^{2,2}_{\loc}(\R^d; \C)$ with $\nabla \psi \in W^{1,2}(\R^d; \C^d)$. Assume in addition that $a \in W^{2,2}_{\loc}(\R^d)$ and $\nabla a \in W^{1,2}(\R^d)$. Then $\Lambda \in W^{1,2}(\R^d)$ and
        \begin{align} 
            \|\nabla \Lambda\|_{L^2}^2 &\leq 6 \|\nabla^2 a\|_{L^2}^2 + 12 \|\nabla^2 a\|_{L^2} \|\nabla^2 \psi\|_{L^2} + 6 \|\nabla^2 \psi\|_{L^2}^2 .
        \end{align}
    \end{enumerate}
\end{theorem}

\subsection{Generalized phase-corrections and parabolic PDEs}

Theorem \ref{thm:main}, together with Proposition~\ref{prop:2} for $s=1$, implies boundedness and continuity for a large class of operators involving phase-dependent first-order nonlinearities, which we think of as more general \emph{phase-corrected gradients}. Precisely, given a Lipschitz function $f: \mathbb{S}^1 \rightarrow \C$ and $\psi \in W^{1,1}_{\loc}(\R^d; \C)$, we define
\begin{align*}
	\nabla_{\!\circ,f} \psi = f(P) \nabla \psi \quad\text{ for some } P \in \mathcal{P}(\psi) .
\end{align*}
As for $\nabla_{\!\circ}\psi$, which corresponds to $f(z) = \conj{z}$, this definition is independent of the choice of the polar factor $P \in \mathcal{P}(\psi)$ since $\nabla\psi = 0$ almost everywhere on $\mathcal{V}$. The following corollary is proven in Section \ref{section:general} below.

\begin{corollary} \label{cor:10}
	Let $d \in \N^+$, $s \in [1, \frac32)$, and let $f: \mathbb{S}^1 \rightarrow \C$ be Lipschitz.
	If $\psi \in W^{s,2}_{\loc}(\R^d; \C)$ with $\nabla \psi \in W^{s-1,2}(\R^d; \C^d)$, then $\nabla_{\!\circ,f} \psi \in W^{s-1,2}(\R^d; \C^d)$ and, if $s > 1$,
	\begin{align} \label{eqn:0-general}
		[\nabla_{\!\circ,f} \psi]_{X^{s-1,2}} &\leq \left( \|f\|_{L^\infty} + \frac{4}{3-2s} [f]_{\Lip} \right) [\nabla \psi]_{X^{s-1,2}} .
	\end{align}
	Furthermore,
	\begin{align*}
		\begin{aligned}
			\psi_n &\rightarrow \psi \text{ in } W^{s,2}_{\loc}(\R^d) \\
			\nabla \psi_n &\rightarrow \nabla \psi \text{ in } W^{s-1,2}(\R^d)
		\end{aligned}
		\quad&
		\Longrightarrow\quad \nabla_{\!\circ,f} \psi_n \rightarrow \nabla_{\!\circ,f} \psi \text{ in } W^{s-1,2}(\R^d) .
	\end{align*}
\end{corollary}

Nonlinear maps of this kind appear in certain parabolic PDEs, such as the mean curvature flow and anisotropic Allen--Cahn and Cahn--Hilliard equations.

\begin{example}[Normalized $p$-Laplacian and mean curvature flow] \label{example:mcf}
	For a function $u: \R^2 \rightarrow \R$ and $p \in [1, \infty)$, the normalized $p$-Laplacian is defined by
	\begin{align*}
		\Delta_p^{\mathrm{N}}(u) &= \frac{1}{p} \Delta_1^{\mathrm{N}}(u) + \frac{p-1}{p} \Delta_\infty^{\mathrm{N}}(u) ,
	\end{align*}
	where
	\begin{align*}
		\Delta_1^{\mathrm{N}}(u) &= |\nabla u| \diver\left( \frac{\nabla u}{|\nabla u|} \right) \quad\text{ and }\quad
		\Delta_\infty^{\mathrm{N}}(u) = \left(\frac{\nabla u}{|\nabla u|}\right)^{\mathrm{T}}\! \nabla^2 u \left(\frac{\nabla u}{|\nabla u|}\right) .
	\end{align*}
	Note that $\Delta_p^{\mathrm{N}}(u)$ is a priori only well-defined at points $x \in \R^2$ with $\nabla u(x) \neq 0$.
	Defining $\psi = \partial_x u + \imu \partial_y u$, we find that
	\begin{align*}
		\Delta_p^{\mathrm{N}}(u) &= \frac12 \Delta u + \frac{p-2}{p} \RE\left[\conj{P}^2 \frac{\partial_x + \imu \partial_y}{2} \psi\right] .
	\end{align*}
	We use this identity to extend the definition of $\Delta_p^{\mathrm{N}}(u)$ almost everywhere, including at critical points of $u$.
	Then Corollary \ref{cor:10} applied with $f(z) = \conj{z}^2$ yields boundedness and continuity of the nonlinear mapping
	\begin{align*}
		\Delta_p^{\mathrm{N}}: W^{2+\theta,2}(\R^2; \R) \longrightarrow W^{\theta,2}(\R^2; \R)
	\end{align*}
	for $\theta \in [0, \frac12)$. More precisely, using \eqref{eqn:0-general} and \eqref{eqn:X-H} to obtain an explicit constant, we have
	\begin{align*}
		\left(\frac12 - \frac{|p-2|}{p} C_\theta\right) \|\Delta u\|_{W^{\theta,2}}
		&\leq \|\Delta_p^{\mathrm{N}}(u)\|_{W^{\theta,2}}
		\leq \left(\frac12 + \frac{|p-2|}{p} C_\theta\right) \|\Delta u\|_{W^{\theta,2}} ,
	\end{align*}
	where $C_\theta = \sqrt{1+2\theta} (1 + \frac{8}{1-2\theta})$ for $\theta > 0$ and $C_0 = \frac{1}{\sqrt{2}}$.
	The parabolic evolution equation
	\begin{align*}
		\partial_t u &= \Delta_p^{\mathrm{N}}(u)
	\end{align*}
	has been studied, for example, in \cite{AttouchiParviainen2018,BanerjeeGarofalo2015}.
	Note that the case $p = 1$ corresponds to the level set formulation of mean curvature flow (see e.g.\ \cite{EvansSpruck1991}).
\end{example}

\begin{example}[Anisotropic Allen--Cahn and Cahn--Hilliard equations] \label{example:cahn}
	Consider $u: \R^2 \rightarrow \R$ and let $A \in C^2(\R^2 \setminus\{0\}; \R)$ be positively $2$-homogeneous. The nonlinear mapping
	\begin{align*}
		\mathbb{A}(u) &= \diver\left(D A(\nabla u)\right) = D^2 A(\nabla u) : \nabla^2 u
	\end{align*}
	appears in anisotropic Allen--Cahn and Cahn--Hilliard equations (see e.g.\ \cite{AlfaroEtAl2010,GarckeKnopfWittmann2023}).
	Defining again $\psi = \partial_x u + \imu \partial_y u$, we find that
	\begin{align*}
		\mathbb{A}(u) &= \RE\left[f_1(P) \partial_x \psi + f_2(P) \partial_y \psi\right],
	\end{align*}
	where
	\begin{align*}
		f_1 &= (\partial_{11} - \imu \partial_{12}) A \quad\text{ and }\quad
		f_2 = (\partial_{21} - \imu \partial_{22}) A .
	\end{align*}
	Note that $D^2 A$ is positively $0$-homogeneous. If in addition $D^2 A\vert_{\mathbb{S}^1}$ is Lipschitz continuous, then Corollary \ref{cor:10} applied with $f = f_j$, $j = 1, 2$, implies boundedness and continuity of
	\begin{align*}
		\mathbb{A}: W^{2+\theta,2}(\R^2; \R) \longrightarrow W^{\theta,2}(\R^2; \R)
	\end{align*}
	for $\theta \in [0, \frac12)$.
\end{example}

\subsection{Quantum hydrodynamics} \label{sec:qhd}
As anticipated in Section \ref{sec:setup}, the Madelung transform \eqref{eqn:formal_madelung}, introduced in 1926 \cite{Madelung1927}, allows one to formally link Schr\"odinger equations to hydrodynamic models. It associates to a complex-valued wave function $\psi = a \eul^{\imu \Theta}$ two quantities: an amplitude $a$ and a velocity field $v = \nabla \Theta$. The physical observables mass and momentum density amount to $\rho=|\psi|^2=a^2$ and $J=\rho v$. This analogy is of particular relevance in the context of Bose--Einstein condensates and superfluidity, where the effective dynamics are modeled by an order parameter $\psi$ at the macroscopic level that satisfies a nonlinear defocusing Schr\"odinger equation
\begin{equation}\label{eqn:NLS}
     \imu\partial_t\psi=-\frac{1}{2}\Delta \psi+f(|\psi|^2)\psi ,
\end{equation}
posed on $\R^d$, with associated total (or Hamiltonian) energy
\begin{align}\label{eqn:total energy NLS}
    \mathcal{H}(\psi)&=\int_{\R^d}\frac{1}{2}|\nabla \psi|^2+F(|\psi|^2)\dd x & F(\rho)&=\int_{\rho_\infty}^\rho f(r) \dd r .
\end{align}
Here $\rho_\infty = \lim_{|x| \rightarrow \infty} \rho(x) \geq 0$ is the background density, which for now we set to $\rho_\infty = 0$. The internal energy density $F(\rho)$ is characterized by the relation $p(\rho)=F'(\rho)\rho-F(\rho)$. A formal computation shows that $(\rho=a^2,v)$ satisfy the hydrodynamic equations
\begin{equation}\label{eqn:QHDvelocity}
    \begin{cases}
       \partial_t\rho+\diver(\rho v)=0 ,
       \\ \partial_t v+(v\cdot\nabla) v+\nabla f(\rho)=\frac{1}{2}\nabla\left(\frac{\Delta a}{a}\right)
    \end{cases} .
\end{equation}
In terms of the mass and momentum density $(\rho,J)$, one recovers the quantum hydrodynamic system (QHD) modeling quantum fluid flow at the macroscopic level:
\begin{equation}\label{eqn:QHD}
    \begin{cases}
    \partial_t \rho+\diver J=0
    \\ \partial_t J+\diver\left(\frac{J\otimes J}{\rho}\right)+\nabla p(\rho)=\frac{1}{2}\rho\nabla\left(\frac{\Delta a}{a}\right) .
    \end{cases}
\end{equation}
The associated total energy is
\begin{equation}\label{eqn:en_QHD1}
\mathcal{E}(\rho,J)=\int_{\R^d} \frac{1}{2}|\nabla a|^2+\frac12{|\Lambda|^2}+F(a^2) \dd x.
\end{equation}
For power-law type nonlinearities $f(|\psi|^2)=|\psi|^{2(\gamma-1)}$ in \eqref{eqn:NLS}, one recovers a barotropic pressure law given by $p(\rho) \sim \rho^{\gamma}$ for \eqref{eqn:QHD}. Moreover, \eqref{eqn:QHDvelocity} may be written as a Hamiltonian system provided that $v$ is irrotational, where the total energy \eqref{eqn:en_QHD1} serves as the Hamiltonian. In view of Proposition \ref{prop:2}, the quantities $\nabla a, \Lambda$ are bounded in $L^2(\R^d)$ provided that $\psi$ is of finite energy. 

Thus, on a formal level, the Madelung transform allows one to pass from \eqref{eqn:NLS} to \eqref{eqn:QHD}. 
Conversely, one may pass from the hydrodynamic formulation to Schrödinger dynamics provided that the hydrodynamic states allow for a wave function lifting (see \cite{AntonelliMarcatiZheng1, AntonelliMarcatiZheng2, Bianchini}). Both links suffer from being largely restricted to the vacuum-free scenario. To the best of our knowledge, Proposition \ref{prop:2} is the only prior result establishing good mapping properties of $\psi\mapsto (a,\Lambda)$ in the presence of vacuum. 

In the following, we discuss applications of this analogy in the study of both NLS and QHD, and highlight the relevance of the fractional regularity provided by Theorem \ref{thm:main}. We refer to \cite{AntonelliMarcati2026} and references therein for an in-depth discussion of quantum hydrodynamic models.

\subsubsection{Fractional regularity of the momentum density}
For the momentum density
\begin{equation} \label{eqn:def-J}
    J = a \Lambda = \IM\left[\conj{\psi} \nabla\psi\right],
\end{equation}
it follows from H\"older and Sobolev inequalities that, given $\psi\in H^1(\R^d;\C)$, 
\begin{equation*}
    J\in L_{\loc}^{p}(\R^d) \quad \text{for} \, \begin{cases}
        p\in [1,2] & \text{if } d=1,\\
        p\in [1,2) & \text{if } d=2,\\
        p\in \left[1,\frac{d}{d-1}\right] & \text{if } d\geq 3.\\
    \end{cases}
\end{equation*}
We provide the following fractional regularity result for $J$.
\begin{proposition} \label{prop:J}
    Let $d \geq 2$, $s \in (1, 2)$. For $j \in \{1, 2\}$ let $\psi_j \in (L^\infty \,\cap\, W^{s,2}_{\loc})(\R^d; \C)$ with $\nabla \psi_j \in \dot{W}^{s-1,2}(\R^d; \C^d)$.
    We have $J_j \in \dot{W}^{s-1,2}(\R^d)$ and
    \begin{align*} 
        [J_j]_{X^{s-1,2}}^2 &\leq 2 \|\psi_j\|_{L^\infty}^2 [\nabla \psi_j]_{X^{s-1,2}}^2,
        \\ [J_1 - J_2]_{X^{s-1,2}}^2 &\leq 4 \left( \|\psi_1\|_{L^\infty}^2 [\nabla (\psi_1 - \psi_2)]_{X^{s-1,2}}^2 + \|\psi_1 - \psi_2\|_{L^\infty}^2 [\nabla \psi_2]_{X^{s-1,2}}^2\right) \!.
    \end{align*}
    If $d = 1$ then the same results hold for $s \in (1, \frac32)$, and also for $s \in [\frac32, 2)$ when the space $\dot{W}^{s-1,2}(\R)$ is replaced with $W^{s-1,2}(\R)$.
\end{proposition}
\begin{remark}
  While the boundedness of $\psi \mapsto J$ can be proven with standard product estimates,  since $J=\IM[\conj{\psi}\nabla\psi]$ does not involve the polar factor (and is therefore not subject to the threshold $s<\frac32$), the proof we give illustrates some aspects of the method used for Theorem \ref{thm:main} and yields additional Lipschitz-type control.
\end{remark}
\subsubsection{The Cauchy problem}
The QHD system amounts to the compressible Euler equations with an additional dispersive stress tensor that depends in a highly nonlinear way on $a$ and its derivatives up to third order. As such, it belongs to the class of Euler--Korteweg systems, but with a particular choice of the capillarity tensor that closely links it to a semilinear NLS equation. This analogy can be exploited to establish a Cauchy theory for \eqref{eqn:QHD}.

A first case concerns the scenario away from vacuum. Typically, one considers non-vanishing boundary conditions such as $\rho(x)\rightarrow 1$ for $|x|\rightarrow \infty$ on $\R^d$ with hydrodynamic initial data $(a_0-1,v_0=\frac{1}{a_0}\Lambda_0)\in H^1(\R^d)\times L^{2}(\R^d)$. Under suitable regularity assumptions, one may lift $(a_0,v_0)$ to a wave function $\psi_0$. Relying on the well-posedness property of the underlying defocusing Gross--Pitaevskii type equation
\begin{equation} \label{eqn:GP}
     \begin{cases}
         \imu \partial_t \psi=-\frac{1}{2}\Delta \psi+(|\psi|^2-1)\psi\\
         |\psi(x)|^2\rightarrow 1 \quad \text{as} \quad |x|\rightarrow \infty
     \end{cases} ,
\end{equation}
the existence of a unique solution $\psi\in C(\R,\mathbf{E}(\R^d))$ follows from \cite{Gerard2008} for $d = 1$, \cite{Gerard2006} for $d \in \{2, 3\}$, and \cite{KillipOhPocovnicuVisan2012} for $d = 4$ and \cite{AHM23NLS} for general nonlinearities. In turn, this solution $\psi$ can be shown to yield existence of weak solutions to \eqref{eqn:QHD}. We refer to \cite{AntonelliMarcati,AntonelliMarcati2012,H} for $d=2,3$, \cite{AntonelliMarcatiZheng1,Mohamad2014,Robert} for $d = 1$, and \cite{AudiardHaspot} for $d \geq 3$.

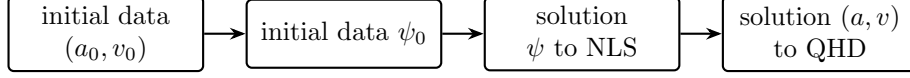
\begin{center}
\begin{tikzpicture}[
box/.style={rectangle,draw,text width=2.34cm,minimum height=.9cm,align=center,thick,rounded corners=2pt},
arrow/.style={-{Stealth[scale=.9]},thick}
]
\node[box] at (0,0) (box1) {initial data $(a_0,v_0)$};
\node[box] at (3.15,0) (box2) {initial data $\psi_0$};
\node[box] at (6.3,0) (box3) {solution $\psi$ to NLS};
\node[box] at (9.45,0) (box4) {solution $(a,v)$ to QHD};
\draw[arrow] (box1)--(box2);
\draw[arrow] (box2)--(box3);
\draw[arrow] (box3)--(box4);
\end{tikzpicture}
\end{center}

This scheme provides a class of weak solutions that we refer to as \emph{Schr\"odinger}-generated solutions, and in general it does not yield uniqueness. For $d=1$ and away from vacuum, it is shown by the third author in \cite{Robert} (see also \cite{Mohamad2014}) that the map \eqref{eqn:formal_madelung} is locally bilipschitz continuous between relevant function spaces for regularities $s > \frac12$. The well-posedness of \eqref{eqn:GP} in fractional energy spaces \cite{KochLiao1} yields a global solution $\psi\in C(\R, \mathbf{E}(\R))$ that enjoys additional fractional regularity. As a result, global well-posedness of QHD with non-vanishing far field for small data $(\rho_0,v_0)\in (1+H^s(\R)) \times H^{s-1}(\R)$ is obtained for $s \geq 1$. The proof fundamentally relies on the fact that a positive distance to vacuum is guaranteed globally in time if the initial data is sufficiently small.

Again away from vacuum, this analogy has also successfully been used to study the long-time behavior of quantum hydrodynamics \cite{AudiardHaspot,FeolaIandoliMurgante,GiulianiScandone}. 

In the presence of vacuum, the Madelung transform is no longer well-defined. On the other hand, Proposition \ref{prop:2} suggests that $(a, \Lambda)$ are the appropriate unknowns to address the Cauchy theory. This motivates the 
definition of hydrodynamic states. 
\begin{definition}[Hydrodynamic state \cite{AntonelliMarcatiZheng1}]
\label{defi:hydrostates} Given $\psi\in H_{\loc}^1(\R^d; \C)$ and $(a,\Lambda)\in H^1_{\loc}(\R^d;\R)\times L^2(\R^d;\R^d)$, we say that $(a,\Lambda)$ is a hydrodynamic state associated to $\psi$ if $a=|\psi|$ and $\Lambda=\IM\left[\conj{P}\nabla\psi\right]$ where $P \in \mathcal{P}(\psi)$.
\end{definition}
Proposition \ref{prop:2} above guarantees the existence of the hydrodynamic states for $\psi\in H_{\loc}^1(\R^d)$ with $\nabla\psi\in L^2(\R^d)$, and their physical admissibility in the sense that $\nabla a(x)=\Lambda(x)=0$ a.e.\ on the vacuum set $\mathcal{V}$. This is a key tool to construct finite energy weak solutions to \eqref{eqn:QHD} for initial data that are hydrodynamic states associated to some $\psi\in H^1(\R^d)$. Specifically, global existence on $\R^d$ with initial data $(a_0,\Lambda_0)\in H^1(\R^d)\times L^2(\R^d)$ is proven in \cite{AntonelliMarcati, AntonelliMarcati2026}, while extensions to nontrivial far field behavior and vortex-type solutions are provided in \cite{H}. Note that weak solutions to \eqref{eqn:QHD} are in general not unique \cite{DonatelliFeireislMarcati, MarkowichSierra}. Our main result yields that Schr\"odinger-generated weak solutions to \eqref{eqn:QHD} dynamically inherit fractional regularity up to the threshold $s<\frac32$.

\begin{proposition}\label{prop:reg-weaksolutions}
    Let $T>0$ and $(a,\Lambda)\in L^{\infty}([0,T];H_{\loc}^1(\R^d; \R) \times L^2(\R^d; \R^d))$ be a Schr\"odinger-generated weak solution to \eqref{eqn:QHD}, i.e.\ there exists a solution $\psi\in C([0,T];H_{\loc}^1(\R^d))$ to the underlying NLS equation \eqref{eqn:NLS} such that
    \begin{equation*}
        a=|\psi| \quad\text{ and }\quad \Lambda=\IM\left[\conj{P}\nabla\psi\right] .
    \end{equation*}
    If the solution satisfies $\nabla \psi \in L^{p}([0,T];\dot{W}^{s-1,2}(\R^d))$ for some $1<s<\frac{3}{2}$ and $1\leq p \leq \infty$, then $\nabla a, \Lambda \in L^{p}([0,T];\dot{W}^{s-1,2}(\R^d))$.
\end{proposition}
Note that this also applies to the weak solutions constructed from purely hydrodynamic initial data in \cite{AntonelliMarcatiZheng1, AntonelliMarcatiZheng2} by means of a wave function lifting. A prototype application of Proposition \ref{prop:reg-weaksolutions} enables us to transfer a potential regularity gain from Kato smoothing for (nonlinear) Schr\"odinger equations to (QHD). Typically, and at least for linear Schr\"odinger equations \cite{ConstantinSaut}, local smoothing estimates yield the local gain of half a space derivative. Our fractional estimate allows one to precisely transfer this regularity gain to the associated hydrodynamic states, up to an arbitrarily small loss.

On a related note, energy sub-critical defocusing Schr\"odinger equations can be expected to enjoy a persistence of regularity property. For nonlinearities such as $f(\rho)=\rho^{\gamma-1}$, $\gamma \geq 1$, the local well-posedness theory of \eqref{eqn:NLS} ensures that the regularity of initial data propagates dynamically \cite{Bourgain2002}. In the defocusing integrable case $d = 1$ and $f(\rho)=\rho$, conserved fractional energies have been constructed both for \eqref{eqn:NLS} and \eqref{eqn:GP} in \cite{KochLiao1,KochTataru2018}, implying propagation of fractional regularities globally in time.

One may ask whether solutions exist that \emph{dynamically saturate} the regularity estimate provided by Theorem \ref{thm:main}. We give below two types of stationary solutions to \eqref{eqn:GP} for which $(a,\Lambda)$ are only of fractional regularity, while $\psi$ is smooth. 
\begin{example}
For $d=1$, we consider the black soliton
\begin{equation}\label{eqn:black-soliton}
\psi(x)=\tanh(x)\in \mathbf{E}(\R) ,
\end{equation}
which is a stationary solution to \eqref{eqn:GP} and correspondingly to \eqref{eqn:QHD} (see \cite{BethuelGravejatSautSmets}). Here $\psi \in H_{\loc}^s(\R)$ for all $s\in \R$ while $a\in H_{\loc}^s(\R)$ only for $s < \frac32$. As \eqref{eqn:black-soliton} constitutes a real-valued solution, one has $\Lambda = 0$. By applying a Galilean boost to \eqref{eqn:black-soliton}, we obtain a family of explicit solutions (depicted in Figure \ref{fig:soliton})
\begin{equation*}
    \psi_0(x)=\eul^{\imu\mu x}\tanh(x), \quad \psi(t,x)=\eul^{\imu\mu\left(x-\frac{\mu}{2}t\right)}\tanh(x-\mu t)
\end{equation*}
for which a straightforward computation shows that
\begin{align*}
    a(t,x)&=|\tanh(x-\mu t)| & \Lambda(t,x) &= \mu a(t,x)= \mu |\tanh(x-\mu t)|.
\end{align*}
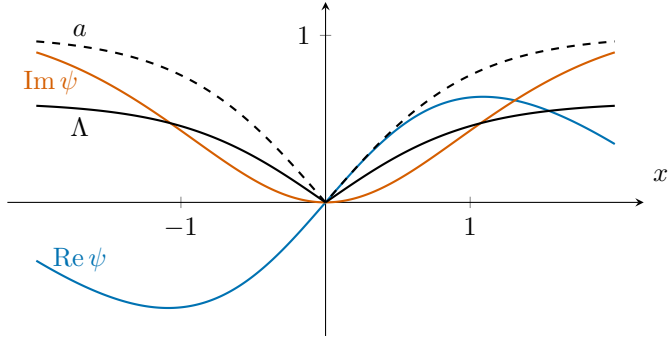
\begin{figure}
    \centering
    \def\m{0.6}
    \begin{tikzpicture}
        \begin{axis}[
            width=10cm,
            height=6cm,
            axis lines=middle,
            xmin=-2.2, xmax=2.2,
            ymin=-0.8, ymax=1.2,
            xtick={-1,1},
            xticklabels={$-1$,$1$},
            ytick={1},
            yticklabels={$1$},
            xlabel={$x$},
            xlabel style={at={(axis description cs:1,0.48)},anchor=west},
            clip=false,
        ]
        
        \addplot[
            plotblue,
            thick,
            domain=-2:2,
            samples=300
        ]
        {cos(deg(\m*x))*tanh(x)};
        \node[plotblue, above] at (axis cs:-1.7,-0.47) {$\RE\psi$};
        
        \addplot[
            plotvermillion,
            thick,
            domain=-2:2,
            samples=300
        ]
        {sin(deg(\m*x))*tanh(x)};
        \node[plotvermillion, below] at (axis cs:-1.9,0.845) {$\IM\psi$};
        
        \addplot[
            black,
            dashed,
            thick,
            domain=-2:2,
            samples=300
        ]
        {abs(tanh(x))};
        \node[above] at (axis cs:-1.7,0.94) {$a$};
        
        \addplot[
            black,
            solid,
            thick,
            domain=-2:2,
            samples=300
        ]
        {\m*abs(tanh(x))};
        \node[below] at (axis cs:-1.7,0.56) {$\Lambda$};    
        \end{axis}
    \end{tikzpicture}
    \caption{A Galilean-boosted black soliton
    $\psi(x)=\eul^{\imu\mu x}\tanh(x)$.}
    \label{fig:soliton}
\end{figure}
We observe that if $\mu \neq 0$, then $a,\Lambda \in H_{\loc}^s(\R)$ only for $s < \frac32$. While this does not dynamically saturate the regularity threshold for $\Lambda$ identified in Theorem \ref{thm:main}, the example shows a loss of regularity of the hydrodynamic states compared to the underlying wave function. In particular, it highlights a need for fractional regularity estimates. 
\end{example}

\begin{example}
For $d=2$, we discuss the example of vortex solutions to the Gross--Pitaevskii equation. For each $n\in \Z$ a stationary vortex solution of degree $n$ is given by
\begin{equation}\label{eqn: d vortex }
    \psi_n(x)=a_n(r)\eul^{\imu n\theta}
\end{equation}
where $x=(r\cos\theta,r\sin\theta)$ are the polar coordinates. The amplitude profiles $a_n(r)$ are radial, non-negative, non-decreasing, and smooth. They have the asymptotic behaviors
\begin{equation*}
    a_n(r)\leq C r^{|n|} \text{ for } r<1 \quad \text{ and } \quad 1-\frac{n^2}{r^2} \leq a_n(r)\leq 1 \text{ for } r \gg 1 .
\end{equation*}
Here
\begin{equation*}
    \Lambda_n(x)=na_n(|x|)\frac{x^{\perp}}{|x|^2} . 
\end{equation*}
For $n=1$ one computes that $\psi_1 \in C^{\infty}(\R^2)$ while $\Lambda_1 \in H_{\loc}^s(\R^2)$ only for $s<1$. Indeed, we have $\Lambda_1 \approx \frac{x^{\perp}}{|x|}$ close to the origin. In terms of integrability, it is well-known that $\Lambda_n\notin L^2(\R^2)$ for $n \neq 0$ only due to the logarithmic divergence at infinity while the multiplication with the amplitude desingularizes the behavior around the origin.
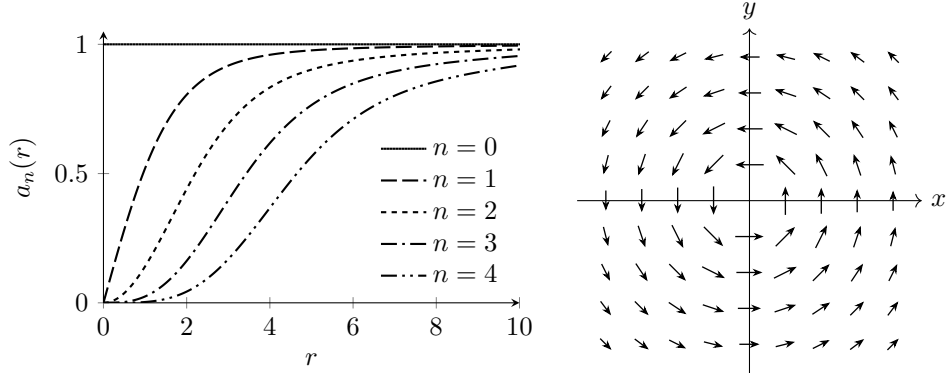
\begin{figure}[t]
\centering

\begin{minipage}[t]{0.56\textwidth}
\centering
\begin{tikzpicture}
\begin{axis}[
    width=\textwidth,
    height=0.73\textwidth,
    axis lines=left,
    xlabel={$r$},
    ylabel={$a_n(r)$},
    xmin=0, xmax=10,
    ymin=0, ymax=1.05,
    legend style={draw=none, fill=none, at={(0.98,0.02)}, anchor=south east},
]
\addplot[thick,dash pattern=on 1pt off 0pt] table[x=r,y=f0] {vortex_profiles.dat};
\addlegendentry{$n=0$}
\addplot[thick,dash pattern=on 6pt off 2pt] table[x=r,y=f1] {vortex_profiles.dat};
\addlegendentry{$n=1$}
\addplot[thick,dash pattern=on 2pt off 2pt] table[x=r,y=f2] {vortex_profiles.dat};
\addlegendentry{$n=2$}
\addplot[thick,dash pattern=on 6pt off 2pt on 1pt off 2pt] table[x=r,y=f3] {vortex_profiles.dat};
\addlegendentry{$n=3$}
\addplot[thick,dash pattern=on 6pt off 2pt on 1pt off 2pt on 1pt off 2pt] table[x=r,y=f4] {vortex_profiles.dat};
\addlegendentry{$n=4$}
\end{axis}
\end{tikzpicture}
\end{minipage}
\hfill
\begin{minipage}[t]{0.40\textwidth}
\centering
\begingroup
\def\vortexdegree{1}
\begin{tikzpicture}[x=0.95cm,y=0.95cm]
\def\plotradius{2.4}
\def\arrowradius{2.3}
\def\gridmax{5}
\def\gridstep{0.5}
\def\lambdaarrowscale{0.7}

\draw[->] (-\plotradius,0) -- (\plotradius,0) node[right] {$x$};
\draw[->] (0,-\plotradius) -- (0,\plotradius) node[above] {$y$};

\foreach \ix in {-\gridmax,...,\gridmax} {
  \foreach \iy in {-\gridmax,...,\gridmax} {
    \pgfmathsetmacro{\xx}{\gridstep*\ix}
    \pgfmathsetmacro{\yy}{\gridstep*\iy}
    \pgfmathsetmacro{\rr}{sqrt(\xx*\xx+\yy*\yy)}
    \pgfmathtruncatemacro{\skipthis}{
      ((\ix==0 && \iy==0) || abs(\xx)>\arrowradius || abs(\yy)>\arrowradius) ? 1 : 0
    }
    \ifnum\skipthis=0
      \pgfmathtruncatemacro{\vortexrow}{min(\vortexlastrow,max(0,round(\vortexlastrow*\rr/12)))}
      \getvortexamp{\vortexrow}
      \pgfmathsetmacro{\amp}{\pgfplotsretval}
      \pgfmathsetmacro{\lamx}{-\vortexdegree*\amp*\yy/(\rr*\rr)}
      \pgfmathsetmacro{\lamy}{\vortexdegree*\amp*\xx/(\rr*\rr)}
      \pgfmathsetmacro{\sx}{\xx - 0.5*\lambdaarrowscale*\lamx}
      \pgfmathsetmacro{\sy}{\yy - 0.5*\lambdaarrowscale*\lamy}
      \pgfmathsetmacro{\tx}{\xx + 0.5*\lambdaarrowscale*\lamx}
      \pgfmathsetmacro{\ty}{\yy + 0.5*\lambdaarrowscale*\lamy}
      \draw[-{Stealth[length=1.3333mm]}, line width=0.5pt] (\sx,\sy) -- (\tx,\ty);
    \fi
  }
}
\end{tikzpicture}
\endgroup
\end{minipage}

\caption{Left: radial vortex profiles $a_n(r)$ for $n=0,1,2,3,4$. Right: the vector field $\Lambda$ for the degree-$1$ vortex $\psi(r,\theta)=a_1(r)\eul^{\imu \theta}$.}
\label{fig:vortex-combined}
\end{figure} 
Note that global well-posedness of \eqref{eqn:GP} for vortex type solutions in $d = 2$ is proven in \cite{BethuelSmets}, and the corresponding infinite energy weak solutions to \eqref{eqn:QHD} are constructed in \cite{H}. 
\end{example}

For the sake of completeness, we mention that the Madelung transform and its inverse are also exploited in the construction of blow-up solutions to energy-supercritical defocusing NLS equations in \cite{Merle}. The implosion type blow-up occurs away from vacuum.

\subsubsection{Analysis of special solutions}
The analysis of special solutions to \eqref{eqn:GP} and in particular their stability is naturally set up in the metric space $(\mathbf{E},d_{\mathbf{E}})$ due to the nontrivial far field behavior. It is convenient to rely on the hydrodynamic states, as the map
\begin{align*}
    \mathbf{E}(\R^d) &\rightarrow H^1(\R^d; \R)\times L^2(\R^d; \R^d)
    \\ \psi &\mapsto (1-a, \Lambda)
\end{align*}
maps into a Banach space.
In a scenario with uniformly positive distance to vacuum, one may prove that $\psi\mapsto (1-a,v=\frac{1}{a}\Lambda)$ is indeed bounded from $\mathbf{E}(\R^d)$ into the product space $H^1(\R^d) \times L^2(\R^d)$ with bounded inverse (modulo constant phase rotations) into a weakened topology (see e.g.\ \cite{Mohamad2014}). Consequently, one can set up the mathematical analysis in Banach spaces, rather than the metric space $(\mathbf{E}(\R^d),d_{\mathbf{E}})$ (see e.g.\ \cite{BGSstab, BGSasympstab}). 

This approach breaks down in scenarios that allow for the occurrence of vacuum. Indeed, the orbital stability analysis in \cite{BethuelGravejatSautSmets} is complicated by this fact and largely relies on the wave function dynamics, rather than the hydrodynamic unknowns. Theorem \ref{thm:main} provides dynamical fractional regularity for the associated hydrodynamic states, allowing for vacuum regions to be included, which we expect to be highly relevant in this framework.

\subsection{Comparison to literature on nonlinear operators on Sobolev spaces}\label{sec:composition}
We begin by reviewing known results on the mapping properties of composition operators
\begin{align*}
    T_f \psi(x) &= f(\psi(x)),
\end{align*}
where $f: \K \rightarrow \K$ and $\psi: \R^d \rightarrow \K$ are measurable functions for $\K \in \{\R, \C\}$. For simplicity, we focus on function spaces with $p = 2$. A well-known result (see e.g.\ \cite[Theorem 2.87]{BCD}) states that if $f$ is smooth with $f(0) = 0$ and $s > 0$, then $T_f: (H^s \cap L^\infty)(\R^d; \R) \rightarrow H^s(\R^d; \R)$ is bounded. If in addition $f'(0) = 0$, then by \cite[Corollary 2.91]{BCD} it is also continuous.
The assumptions on $f$ can be substantially weakened (see for example \cite{BrezisMironescu2001}).

Note that for composition operators of the form $T_f$, boundedness and continuity are not, in general, equivalent. On the Sobolev spaces $W^{1,p}(\R^d)$ with $1 \leq p < \infty$ it has indeed been established that boundedness implies continuity \cite{MarcusMizel1979}.
This problem is also studied in \cite{BourdaudSickel}, which we refer to for an in-depth review of the matters discussed below. Note also \cite{BorowskiCianchi2025} for a very general setting using Orlicz--Sobolev spaces over open sets.

Let us transition now to the case of non-smooth functions, the prototypical example being the modulus $f(\psi) = |\psi|$. In the low-regularity range $0 < s < 1$ the inequality $|f(\psi(x)) - f(\psi(y))| \leq \|D f\|_{L^\infty} |\psi(x) - \psi(y)|$, together with $f(0) = 0$, suffices to prove that $T_f: H^s(\R^d; \C) \rightarrow H^s(\R^d; \C)$ is bounded when $D f \in L^\infty(\R^2; \R^{2 \times 2})$ (see \cite[Theorem 6]{BourdaudSickel} for a general statement). For high regularities $s > \frac32$ we refer to \cite[Theorems 8, 9]{BourdaudSickel}, which imply for the one-dimensional case that $T_f: H^s(\R; \R) \rightarrow H^s(\R; \R)$ is bounded if and only if $f \in H^s_{\loc}(\R; \R)$ and $f(0) = 0$. This is indeed not the case for the modulus, where one can also directly demonstrate the unboundedness, in any dimension, with explicit counterexamples such as we give in Section \ref{section:unbounded} below.

The interesting range is $1 < s < \frac32$, where affirmative results for the boundedness of $T_{|\,\cdot\,|}: H^s(\R^d; \R) \rightarrow H^s(\R^d; \R)$ when $s < \frac32$ were first obtained in \cite{BourdaudMeyer,Oswald}. More general results may be found in \cite{BourdaudSickel,Kateb2000}. Boundedness has also been established on bounded domains with $C^2$-boundary in \cite[Theorem 1.1]{Li}.

The case $s = 1$ is trivial and the case $s = \frac32$ fails. Here counterexamples are constructed in \cite{Oswald} for the real-valued case, but the same approach covers the complex-valued case.
We note the positive result \cite{Moussai2018}, which states that $T_{|\,\cdot\,|}: B^{\frac32}_{2,1}(\R^d; \R) \rightarrow B^{\frac32}_{2,\infty}(\R^d; \R)$ is bounded (see also Proposition \ref{prop:besov} below). To our knowledge, there are no known results on the continuity of $T_{|\,\cdot\,|}: H^s(\R^d; \C) \rightarrow H^s(\R^d; \R)$ for $s \in (1, \frac32)$. We resolve this positively in Corollary \ref{cor:modulus} above.
We mention also \cite{Kateb2003}, where for $\mu > 1$ and $s \in (1, \mu+\frac12)$ the boundedness of $T_{|\,\cdot\,|^\mu}: (H^s \,\cap\, L^\infty)(\R^d; \R) \rightarrow H^s(\R^d)$ is proven, and a similar result \cite[Theorem 3.1]{PineauTaylor2025} on the boundedness of $T_{|\,\cdot\,|^{\mu-1} \cdot}$ from $(H^s \,\cap\, H^{\frac12+\epsilon})(\R; \C)$ to $H^s(\R)$ for small $\epsilon > 0$ and $s \in (\mu, \mu + \frac12)$.
Lastly, we mention \cite{PonceSpector2020} for related results.

The mapping $\psi \mapsto \nabla_{\!\circ} \psi$, which we are interested in, is not a composition operator, but since $\RE\left[\nabla_{\!\circ} \psi\right] = \nabla T_{|\,\cdot\,|}(\psi)$ its mapping properties are related to those of $T_{|\,\cdot\,|}$.
Regarding this real part, we use the arguments from \cite{BourdaudMeyer} in Section \ref{section:a} to prove, for the sake of completeness, that $\nabla a \in \dot{W}^{s-1,2}(\R^d)$ when $\nabla \psi \in \dot{W}^{s-1,2}(\R^d)$.
Our main object of interest is the vector field $\Lambda = \IM\left[\nabla_{\!\circ} \psi\right]$. The difficulty lies in the necessity to understand how the regularity of the function interacts with the distance to vacuum. For the case of the real part, the convenient formula \eqref{eqn:a-reduction} permits a reduction to the case of real-valued $\psi$, which greatly simplifies the possible range of behaviors near vacuum. We do not expect that such an identity and corresponding reduction to real-valued $\psi$ are available for $\Lambda$ since it is not, in general, a conservative vector field, owing to the fact that it does not correspond to a composition operator.

Since the polar factor $P$ plays a major role, we mention also \cite{BourgainBrezisMironescu2000,MironescuRussSire2020} where the boundedness of the map $[P = \eul^{\imu \Theta} \mapsto \Theta]$ is studied.

\subsection{Idea of the proof}\label{sec:idea}

We first point out that under the hypothesis of strict absence of vacuum, namely $|\psi|>c>0$ for some $c>0$, the mapping properties of $\psi\mapsto \nabla_{\!\circ}\psi$ in appropriate Sobolev spaces $H^s(B)$, where $B \subset \R^d$ is some ball, follow by classical Bony-type composition estimates, as can be found in \cite{BCD}, or alternatively by directly working with the definition of the Sobolev--Slobodeckij seminorms. In the presence of vacuum, we cannot appeal to standard techniques. For example, let us try to apply the fractional Leibniz rule to estimate $\Lambda = \IM\left[\conj{P} \nabla \psi\right]$:
\begin{equation*}
	\|\Lambda\|_{H^{s-1}} \lesssim \|\conj{P}\|_{L^\infty} \|\nabla \psi\|_{H^{s-1}} + \|\conj{P}\|_{H^{s-1}} \|\nabla \psi\|_{L^\infty}.
\end{equation*}
Note that the second term cannot be controlled. In fact, we should not expect product estimates to work for the following reason: for functions $\psi(x) \in H^s(\R)$ the amplitude $a(x)$ may vanish to arbitrarily high order as $x \rightarrow 0$, and hence the polar factor $P(x)$ can have arbitrarily irregular behavior near $x = 0$.
This suggests that the regularity of $P$ must be studied in a local-in-$x$ sense, where the size of the amplitude $a(x)$ is considered. Indeed, we can informally think of $\|\Lambda\|_{H^{s-1}(\R^d)} = \|a \nabla \log P\|_{H^{s-1}(\R^d)}$ as measuring regularity of $\log P$ with a weight $a$. 
For the above reason, approaches based on the Fourier transform face difficulties. It is not clear to us whether microlocal approaches may resolve this.

Instead, we work with the Sobolev--Slobodeckij seminorm $[\Lambda]_{W^{s-1,2}}^2$.
It is natural to try to use the nonlinear structure $\Lambda = \IM\left[ \frac{\conj{\psi} \nabla \psi}{a} \right]$, manipulating $\Delta_h \Lambda$ with the goal of making differences $\Delta_h \nabla \psi$ and $\Delta_h \nabla a$ appear.
For simplicity, let us assume $d = 1$. We decompose
\begin{align*}
    \Delta_h \Lambda(x) &= \IM\left[ \frac{\conj{\psi(x)}}{a(x)} \Delta_h \psi'(x) \right] 
    \\ &+ \IM\left[ \frac{\conj{\Delta_h \psi}(x)}{a(x)} \psi'(x-h) \right] 
    \\ &+ \IM\left[ \conj{\psi(x-h)} \psi'(x-h) \Delta_h a^{-1}(x) \right] 
    \\ &= I_1 + I_2+ I_3 .
\end{align*}
The first term $I_1$ is unproblematic. 
The second term $I_2$ presents a dilemma: if we estimate $|\conj{\Delta_h \psi(x)}| \leq a(x) + a(x-h)$ in order to remove the singularity caused by the vacuum, then the singularity $|h|^{-d-2(s-1)}$ in \eqref{eqn:50} is uncontrolled.
On the other hand, if we keep $\Delta_h \psi$ and use \eqref{eqn:50}, then we utilize only $s - 1$ regularity of $\psi$.
In order to use the full $s$ regularity of $\psi$, we need to work with the second-order difference operator
\begin{align*}
    R_h \psi(x) &= \Delta_h \psi(x) - \psi'(x) h
\end{align*}
which satisfies (see \eqref{eqn:9} below)
\begin{align} \label{eqn:9-intro}
    \int_{\R} |h|^{-d-2(s-1)} \int_{\R} |h^{-1} R_h \psi(x)|^2 \dd x \dd h &= \frac{1}{2 s} [\psi']_{W^{s-1,2}}^2 .
\end{align}
We can write
\begin{align*}
    I_2 &= \IM\left[ \frac{\conj{\Delta_h \psi(x)}}{a(x)} \frac{\Delta_h \psi(x)}{h} \right]
    + \IM\left[ \frac{\conj{\Delta_h \psi(x)}}{a(x)} h^{-1} R_{-h} \psi(x-h) \right] 
    \\ &= I_{2,1} + I_{2,2}
\end{align*}
and
\begin{align*}
    I_3 &= \IM\left[ \conj{\psi(x-h)} \frac{\Delta_h \psi(x)}{h} \right] \Delta_h a^{-1}(x)
    + \IM\left[ \conj{\psi(x-h)} h^{-1} R_{-h} \psi(x-h) \right] \Delta_h a^{-1}(x) 
    \\ &= - \IM\left[ P(x) \conj{P(x-h)} \right] h^{-1} \Delta_h a(x)
    + \IM\left[ \conj{\psi(x-h)} h^{-1} R_{-h} \psi(x-h) \right] \Delta_h a^{-1}(x)
     \\ &= I_{3,1} + I_{3,2} .
\end{align*}
Conveniently, the first term $I_{2,1}$ vanishes.
Furthermore, we have the surprising identity
\begin{align*}
    I_{2,2} + I_{3,2} &= \IM\left[\left(\conj{P(x)} - \conj{P(x-h)}\right) h^{-1} R_{-h} \psi(x-h) \right],
\end{align*}
so \eqref{eqn:9-intro} yields
\begin{align*}
    \int_{\R} |h|^{-d-2(s-1)} \int_{\R} |I_{2,2} + I_{3,2}|^2 \dd x \dd h &\lesssim [\psi']_{W^{s-1,2}}^2 \,.
\end{align*}
The final term $I_{3,1}$ is problematic. We progress by observing that
\begin{align*}
    |I_{3,1}|^2 &= \left( 2 \left(1 - \RE\left[ P(x) \conj{P(x-h)} \right]\right) - \left(1 - \RE\left[ P(x) \conj{P(x-h)} \right]\right)^2 \right) (h^{-1} \Delta_h a(x))^2 
\end{align*}
where
\begin{align*}
    \left(1 - \RE\left[P(x) \conj{P(x-h)}\right]\right) h^{-1} a(x)
    &= \RE\left[ \conj{P(x-h)} h^{-1} R_{-h} \psi(x-h) \right]
    - h^{-1} R_{-h} a(x-h) .
\end{align*}
This means we can estimate
\begin{align*}
    |I_{3,1}|^2 &\lesssim |h|^{-2} |\Delta_h \psi(x)| (|R_h \psi(x)| + |R_{-h} \psi(x-h)|) 
    \\ &+ |h|^{-2} |\Delta_h a(x)| (|R_h a(x)| + |R_{-h} a(x-h)|) .
\end{align*}
Here only one power of $R_h \psi(x)$ or $R_h a(x)$ is present, instead of the two powers necessary in order to control the term as above.
Nevertheless, with the Cauchy--Schwarz inequality and \eqref{eqn:9-intro} we can still estimate
\begin{align} \label{eqn:crude}
    \int_{\R} |h|^{-1-2(s-1)} \int_{\R} |I_{3,1}|^2 \dd x \dd h &\lesssim \|\nabla \psi\|_{L^2}^2 + [\psi']_{W^{2(s-1)+\nu,2}}^2
    \\ &+ \|\nabla a\|_{L^2}^2 + [a']_{W^{2(s-1)+\nu,2}}^2
\end{align}
for any $\nu \in (0, \frac52-2s)$. Note that by Proposition \ref{prop:1} below we can estimate $[a']_{W^{2(s-1)+\nu,2}} \lesssim [\psi']_{W^{2(s-1)+\nu,2}}$ (see \cite{BourdaudMeyer}). The problem here is that the formulation \eqref{eqn:50} of the Sobolev norm is only admissible for $s-1 \in (0, 1)$, which forces us to use $s-1+\nu$ more regularity of $\psi$ than we intended. This crude estimate still turns out to be useful in Step 3.1 of the proof of Theorem \ref{thm:1}.

Theorem \ref{thm:1} avoids this loss of regularity via the usage of a number of complicated nonlinear identities \eqref{eqn:id1}--\eqref{eqn:id9}, which are of algebraic nature.
Ultimately, these identities produce only three types of terms: those which can be controlled using \eqref{eqn:9-intro} like $I_{2,2} + I_{3,2}$ above, or have a good sign, those which take the form $\RE\left[R_h \psi(x) \conj{\Delta_h \psi(x)}\right]$ or $R_h a(x) \Delta_h a(x)$ \textit{without} absolute value, and lastly a number of further terms that can be estimated to yield a strictly smaller copy of the left-hand side. We can then subtract this copy from the left-hand side to close the argument.
The key to estimating the second kind of terms is \eqref{eqn:10} below, which states
\begin{align*}
    \int_{\R^d} |h|^{-d-2s} \int_{\R^d} \RE\left[R_h \psi(x) \conj{\Delta_h \psi(x)}\right] \dd x \dd h &= \frac{s-1}{2 s (2s-1)} [\nabla \psi]_{X^{s-1,2}}^2
\end{align*}
as long as the integrand on the left-hand side is indeed integrable.
Here a cancellation weakens the singularity. Note that this step requires using translation and scaling symmetry of the integral, which fails on domains. Caution is necessary: we must already know the integral to exist in order for this cancellation to be rigorous. This is where the crude estimate \eqref{eqn:crude} turns out to be useful, as it implies the integrability of $w_\gamma |\Delta_h \psi(x)| |R_h \psi(x)|$ over $(x, h) \in \R^{d+d}$ where $\gamma \in (0, \frac14)$ is chosen so that $s - 1 = 2 \gamma + \nu$.
This allows us to obtain a small amount of regularity for $\Lambda$. We finally gain the full regularity by repeating the argument in a weaker fashion, i.e.\ with a cutoff $\epsilon < |h| < \epsilon^{-1}$, and taking the limit $\epsilon \searrow 0$. Since the cutoff is not invariant under scaling symmetry, we use the bootstrap-type Lemma \ref{lem:bootstrap}.
Our argument yields more: once the integrability of all relevant terms and the existence of the limit $\epsilon \searrow 0$ have been established, we obtain the explicit formula \eqref{eqn:formula-main}.

Note that the obtained estimates depend only on $[\nabla \psi]_{W^{s-1,2}}$ and $[\nabla a]_{W^{s-1,2}}$. There are, however, some points where a certain control of the low-frequency part of $\nabla \psi$ is needed. Specifically, observe that \eqref{eqn:crude} requires $\nabla \psi \in L^2(\R^d)$ (see Lemma \ref{lem:20} below). Our proof strategy therefore requires $\nabla \psi \in W^{s-1,2}(\R^d)$ to be in an inhomogeneous Sobolev space. 
In Corollary \ref{cor:1} we use an approximation argument to extend the result to $\nabla \psi \in \dot{W}^{s-1,2}(\R^d)$.

This argument also works in higher dimensions, although it is fundamentally an argument along one-dimensional subspaces of $\R^d$. This means it only controls the quantity $\Delta_h \Lambda(x) \cdot h$. Fortunately, the fractional Korn's inequality (see \cite[Theorem 1]{ScottMengesha2019} and \eqref{eqn:korn} below) implies that controlling $\Delta_h \Lambda(x) \cdot h$, which corresponds to the ``symmetric part of the fractional gradient'' of $\Lambda$, is sufficient.

We expect that the argument can be transferred to the case of cylindrical domains $\T^k \times \R^d$, although this requires nontrivial changes.

As mentioned above, the argument does not work on domains due to its usage of translation and scaling symmetry. If one carries out these symmetrizations nonetheless, certain terms that require additional regularity near the boundary appear.

\subsection*{Statements and Declarations}
\begin{itemize}
    \item This work was funded by the Deutsche Forschungsgemeinschaft (DFG, German Research Foundation) – Project-ID 258734477 – SFB 1173. Raffaele Scandone acknowledges partial support to INdAM-Gnampa through the project ``Fluid dynamics and transport models with irregular velocity fields''.
    \item No datasets were generated or analyzed during the current study.
    \item The authors have no competing interests to declare that are relevant to the content of this article.
    \item SymPy \cite{sympy} was used to carry out symbolic computations for the purposes of exploration and verification of various algebraic identities.
    \item Various versions of LLMs\footnote{In particular, the OpenAI services ChatGPT and Codex.} were used for exploratory discussion, literature search, proofreading, typesetting, creation of figures, and in the development of some mathematical arguments. Beyond these aspects, they were not used for the writing of the document. \\
    Specifically concerning the development of mathematical arguments, we would like to clarify where significant ideas were contributed. 
    We note their use in providing the known construction in Lemma \ref{lem:g-n}, when prompted to establish the sharpness results in Lemma \ref{lem:sharp}, as well as in the development of downstream counterexamples in Section \ref{section:unbounded}.
    Furthermore, in the derivation of Theorem \ref{thm:continuity} from the results of Section \ref{section:lambda}, they proposed the use of a scale decomposition and pointwise convergence of the integrand when prompted to exploit the structure of \eqref{eqn:0} to obtain convergence of norms, which previously needed additional assumptions. \\
    Lastly, when prompted to investigate generalizations and applications of the main result, they provided the basic observation used to obtain Corollary \ref{cor:10} and noted its relevance for Examples \ref{example:mcf} and \ref{example:cahn}.
\end{itemize}

\section{Preliminaries}

When we do not specify the target space of a function space under discussion, the intended meaning is that any finite-dimensional vector space is admissible.

Recall the definitions of the space of test functions $\mathcal{D}(\R^d) = C_c^\infty(\R^d)$, the space of distributions $\mathcal{D}'(\R^d)$, the Schwartz space $\mathcal{S}(\R^d)$, and the corresponding space of tempered distributions $\mathcal{S}'(\R^d)$. For $f \in \mathcal{S}(\R^d)$, we recall the definition of the Fourier transform
\begin{align*}
    \ha{f}(\xi) &= \frac{1}{(2 \pi)^{\frac{d}{2}}} \int_{\R^d} \eul^{- \imu \xi \cdot x} f(x) \dd x ,
\end{align*}
and for $f \in \mathcal{S}'(\R^d)$ the corresponding distributional Fourier transform $\ha{f}(\phi) = f(\ha{\phi})$ for $\phi \in \mathcal{S}(\R^d)$.

For $\Omega \subseteq \R^d$ and $\theta \in (0, 1)$, we recall the definition of the Sobolev--Slobodeckij seminorm
\begin{align*}
    [f]_{W^{\theta,2}(\Omega)} &= \left( \int_{\Omega} \int_{\Omega} |x-y|^{-d-2\theta} |f(x) - f(y)|^2 \dd x \dd y \right)^{\frac12}
\end{align*}
and the fractional inhomogeneous Sobolev spaces
\begin{align*}
    W^{\theta,2}(\Omega) &= \{f \in L^2(\Omega): [f]_{W^{\theta,2}(\Omega)} < \infty\}
    \\ W^{1+\theta,2}(\Omega) &= \{f \in L^2(\Omega): \nabla f \in W^{\theta,2}(\Omega)\}
    \\ W^{1+\theta,2}_{\loc}(\R^d) &= \{f \in L^2_{\loc}(\R^d): f\vert_\Omega \in W^{1+\theta,2}(\Omega) \text{ for all bounded open } \Omega \subseteq \R^d \}
\end{align*}
with the standard associated norms.
Recall also for $k \in \N$ the definition of the inhomogeneous Sobolev spaces
\begin{align*}
    W^{k,2}(\Omega) &= \{f \in L^2(\Omega): \nabla^1 f, \dots, \nabla^k f \in L^2(\Omega)\} 
    \\ W^{k,2}_{\loc}(\R^d) &= \{f \in L^2_{\loc}(\R^d): f\vert_\Omega \in W^{k,2}(\Omega) \text{ for all bounded open } \Omega \subseteq \R^d \}
\end{align*}
with the standard associated norms.
When $\Omega = \R^d$, we write $[\,\cdot\,]_{W^{\theta,2}} = [\,\cdot\,]_{W^{\theta,2}(\R^d)}$.
We define the homogeneous Sobolev--Slobodeckij spaces
\begin{align*}
    \dot{W}^{\theta,2}(\R^d) &= \conj{C_c^\infty(\R^d)}^{[\,\cdot\,]_{W^{\theta,2}}}
\end{align*}
and note that they are Hilbert spaces.
For $\theta < \frac{d}{2}$, which is the only range where we use these spaces, it is known (see e.g.\ \cite[Theorem 3.1]{BrascoGomezCastroVazquez2021}) that the identification
\begin{align} \label{eqn:characterization}
    \dot{W}^{\theta,2}(\R^d) = \{f \in L^{\frac{2d}{d-2\theta}}(\R^d): [f]_{W^{\theta,2}} < \infty\}
\end{align}
is possible.
We are interested specifically in vector fields $F \in \dot{W}^{\theta,2}(\R^d; \C^d)$. Here we generally use the alternative norm
\begin{align} \label{eqn:def-korn-norm}
    [F]_{X^{\theta,2}} &= \left( \int_{\R^d} \int_{\R^d} |x-y|^{-d-2-2\theta} |(F(x) - F(y)) \cdot (x-y)|^2 \dd x \dd y \right)^{\frac12} .
\end{align}
By \cite[Theorem 1]{ScottMengesha2019} there exists a constant $C(d, \theta)$ such that
\begin{align} \label{eqn:korn}
    [F]_{X^{\theta,2}} \leq [F]_{W^{\theta,2}} &\leq C(d, \theta) [F]_{X^{\theta,2}}
\end{align}
for all $F \in C_c^\infty(\R^d)$, meaning that $[\,\cdot\,]_{X^{\theta,2}}$ and $[\,\cdot\,]_{W^{\theta,2}}$ are equivalent norms on $\dot{W}^{\theta,2}(\R^d)$. Further details are given in Lemma \ref{lem:characterizations} above.

We define for $h \in \R^d$ the first- and second-order finite difference operators
\begin{align}
    \label{eqn:def-Delta-h} \Delta_h f(x) &= f(x) - f(x-h)
    \\ \label{eqn:def-R-h}  R_h f(x) &=  \Delta_h f(x) - \nabla f(x) \cdot h .
\end{align}
For $\epsilon \geq 0$ we define the annulus
\begin{align*}
    A_\epsilon &= \begin{cases}
        \{h \in \R^d: |h| \in (\epsilon, \epsilon^{-1})\} &, \epsilon > 0
        \\ \R^d &, \epsilon = 0
    \end{cases} .
\end{align*}
We often use the following shorthand notation for an important weight:
\begin{align*}
    w_\theta &= w_\theta(h) = |h|^{-d-2-2\theta} .
\end{align*}
\subsection{Inequalities}

The following lemma summarizes several equivalent formulations of a technical assumption that is needed in Lemma \ref{lem:inequ} below.
\begin{lemma}[Integrability] \label{lem:equiv}
    Let $d \in \N^+$ and $\alpha \in (0, 1)$. Let $f \in W^{1+\alpha,2}_{\loc}(\R^d; \C)$ with $[\nabla f]_{X^{\alpha,2}} < \infty$. The following are equivalent:
    \begin{enumerate}[(i)]
        \item $w_\alpha R_h f(x) \conj{\Delta_h f(x)} \in L^1_h L^1_x(\R^{d+d}; \C)$.
        \item $w_\alpha (\Delta_h \nabla f(x) \cdot h) \conj{\nabla f(x) \cdot h} \in L^1_h L^1_x(\R^{d+d})$.
        \item $w_\alpha r^{-2\alpha} (\Delta_{rh} \nabla f(x) \cdot h) \conj{\nabla f(x) \cdot h} \in L^1_h L^1_x L^1_r(\R^{d+d} \times [0, 1])$.
    \end{enumerate}
\end{lemma}
\begin{proof}
    For (i) $\Longrightarrow$ (ii) we calculate
    \begin{align*}
        (\Delta_h \nabla f(x) \cdot h) \conj{\nabla f(x) \cdot h} 
        &= - R_h f(x) \conj{\Delta_h f(x)} 
        \\ &+ R_{-h} f(x-h) \conj{\Delta_{-h} f(x-h)} 
        \\ &- (\Delta_h \nabla f(x) \cdot h) \conj{R_h f(x)} .
    \end{align*}
    Here the second term is integrable by the transformation $(h, x) \mapsto (-h, x-h)$ and the last term is integrable by the Cauchy--Schwarz inequality and (i) of Lemma \ref{lem:inequ} below. Next, (ii) $\Longrightarrow$ (iii) follows from the transformation $h \mapsto r h$. Lastly, for (iii) $\Longrightarrow$ (i) we write
    \begin{align*}
        R_h f(x) \conj{\Delta_h f(x)} &= R_h f(x) \conj{R_h f(x)} 
        - \int_0^1 (\Delta_{rh} \nabla f(x) \cdot h) \conj{\nabla f(x) \cdot h} \dd r
    \end{align*}
    and note that (i) of Lemma \ref{lem:inequ} below implies the integrability of the first term, while (iii) implies the integrability of the second term.
\end{proof}

Our main method of estimating quantities involving $f$ and $\nabla f$ by $[\nabla f]_{X^{\alpha,2}}$ is the collection of inequalities given below.
\begin{lemma}[Higher-order differences] \label{lem:inequ}
    Let $d \in \N^+$ and $\alpha \in (0, 1)$. Let $f \in W^{1+\alpha,2}_{\loc}(\R^d; \C)$ with $[\nabla f]_{X^{\alpha,2}} < \infty $.
    \begin{enumerate}[(i)]

        \item We have $w_\alpha |R_h f(x)|^2 \in L^1_h L^1_x(\R^{d+d})$ and
        \begin{align} \label{eqn:9}
            \int_{\R^d} w_\alpha \int_{\R^d} |R_h f(x)|^2 \dd x \dd h &= \frac{1}{2+2\alpha} [\nabla f]_{X^{\alpha,2}}^2 .
        \end{align}
        Similarly, we have $w_\alpha \RE\left[R_h f(x) \conj{R_{-h} f(x-h)}\right] \in L^1_h L^1_x(\R^{d+d})$ and
        \begin{align} \label{eqn:9-alt}
            \int_{\R^d} w_\alpha \int_{\R^d} \RE\left[R_h f(x) \conj{R_{-h} f(x-h)}\right] \dd x \dd h &= \frac{\alpha}{2+2\alpha} [\nabla f]_{X^{\alpha,2}}^2 .
        \end{align}
        
        \item If $w_\nu R_h f(x) \conj{\Delta_h f(x)} \in L^1_h L^1_x(\R^d \times \R^d)$ and $[\nabla f]_{X^{\nu,2}} < \infty$ for some $\nu \in (0, 1)$, then for all $\epsilon > 0$ we have
        \begin{align} \label{eqn:10-eps}
            & \int_{A_\epsilon} \!\!w_\alpha\!\! \int_{\R^d} \RE\left[R_h f(x) \conj{\Delta_h f(x)}\right] \dd x \dd h
            \\\nonumber &= \int_{A_\epsilon} \!\!w_\alpha\!\! \int_{\R^d} |R_h f(x)|^2 \dd x \dd h
            - \frac12 \int_0^1 r^{2\alpha} \int_{r A_\epsilon} \!\!w_\alpha\!\! \int_{\R^d} |\Delta_h \nabla f(x) \cdot h|^2 \dd x \dd h \dd r .
        \end{align}

        \item If $w_\alpha R_h f(x) \conj{\Delta_h f(x)} \in L^1_h L^1_x(\R^{d+d})$, then
        \begin{align} \label{eqn:10}
            \int_{\R^d} w_\alpha \int_{\R^d} \RE\left[R_h f(x) \conj{\Delta_h f(x)}\right] \dd x \dd h
            &= \frac{\alpha}{2 (1+\alpha) (1+2\alpha)} [\nabla f]_{X^{\alpha,2}}^2 .
        \end{align}
    \end{enumerate}
\end{lemma}
\begin{proof}
    \begin{enumerate}[(i)]
        \item We decompose
        \begin{align*}
            & \int_{\R^d} w_\alpha \int_{\R^d} |R_h f(x)|^2 \dd x \dd h 
            \\ &= \int_{\R^d} w_\alpha \int_0^1 \int_0^1  \int_{\R^d} \RE\left[(\Delta_{rh} \nabla f(x) \cdot h) \conj{\Delta_{th} \nabla f(x) \cdot h}\right] \dd x \dd t \dd r \dd h
            \\ &= \frac12 \int_{\R^d} w_\alpha \int_0^1 \int_0^1  \int_{\R^d} |\Delta_{rh} \nabla f(x) \cdot h|^2 + |\Delta_{th} \nabla f(x) \cdot h|^2 \dd x \dd t \dd r \dd h
            \\ &- \frac12 \int_{\R^d} w_\alpha \int_0^1 \int_0^1  \int_{\R^d} |(\Delta_{rh} \nabla f(x) - \Delta_{th} \nabla f(x)) \cdot h|^2 \dd x \dd t \dd r \dd h
            \\ &= (I) + (II) .
        \end{align*}
        Then we use a change of variables to write
        \begin{align*}
            (I) &= \int_0^1 r^{2\alpha} \dd r \int_{\R^d} w_\alpha \int_{\R^d} |\Delta_h \nabla f(x) \cdot h|^2 \dd x \dd h
        \end{align*}
        and
        \begin{align*}
            (II) &= - \frac12 \int_{\R^d} w_\alpha \int_0^1 \int_0^1 \int_{\R^d} |\Delta_{(r-t)h} \nabla f(x) \cdot h|^2 \dd x \dd t \dd r \dd h
            \\ &= - \frac12 \int_0^1 \int_0^1 |r-t|^{2\alpha} \dd t \dd r \int_{\R^d} w_\alpha \int_{\R^d} |\Delta_{h} \nabla f(x) \cdot h|^2 \dd x \dd h .
        \end{align*}
        We obtain \eqref{eqn:9} from
        \begin{align*}
            (I) + (II) &= \frac12 \left( \frac{2}{1+2\alpha} - \frac{2}{(1+2\alpha)(2+2\alpha)} \right) [\nabla f]_{X^{\alpha,2}}^2
            \\ &= \frac{1}{2+2\alpha} [\nabla f]_{X^{\alpha,2}}^2 ,
        \end{align*}
        and then \eqref{eqn:9-alt} via the identity
        \begin{align*}
            |\Delta_h \nabla f(x) \cdot h|^2 = |R_h f(x) + R_{-h} f(x-h)|^2
            &= |R_h f(x)|^2 + |R_{-h} f(x-h)|^2 
            \\ &+ 2 \RE\left[R_h f(x) \conj{R_{-h} f(x-h)}\right] .
        \end{align*}

        \item Here we write
        \begin{align*}
        	& \int_{A_\epsilon} w_\alpha \int_{\R^d} \RE\left[R_h f(x) \conj{\Delta_h f(x)}\right] \dd x \dd h
        	\\ &= \int_{A_\epsilon} w_\alpha \int_{\R^d} \RE\left[R_h f(x) \conj{R_h f(x)}\right] \dd x \dd h
        	\\ &- \int_{A_\epsilon} w_\alpha \int_{\R^d} \int_0^1 \RE\left[(\Delta_{rh} \nabla f(x) \cdot h) \conj{\nabla f(x) \cdot h}\right] \dd r \dd x \dd h .
        \end{align*}
        For the second integral, we note that (iii) from Lemma \ref{lem:equiv}, applied with $\nu$ instead of $\alpha$, implies $(\Delta_{rh} \nabla f(x) \cdot h) \conj{\nabla f(x) \cdot h} \in L^1_x L^1_r(\R^d \times [0, 1])$ for almost all $h \in \R^d$.
        Now the second integral can be symmetrized and estimated as follows:
        \begin{align*}
        	& - \int_{A_\epsilon} w_\alpha \int_{\R^d} \int_0^1 \RE\left[(\Delta_{rh} \nabla f(x) \cdot h) \conj{\nabla f(x) \cdot h}\right] \dd r \dd x \dd h
            \\ &= - \frac12 \int_{A_\epsilon} w_\alpha \int_{\R^d} \int_0^1 |\Delta_{rh} \nabla f(x) \cdot h|^2 \dd r \dd x \dd h
            \\ &= - \frac12 \int_0^1 r^{2\alpha} \int_{r A_\epsilon} w_\alpha \int_{\R^d} |\Delta_h \nabla f(x) \cdot h|^2 \dd x \dd h \dd r .
        \end{align*}

        \item We perform the argument of (ii) with $\epsilon = 0$, observing that the required integrability is provided by Lemma \ref{lem:equiv}, and that now we can calculate the second term exactly.
    \end{enumerate}
\end{proof}

\subsection{Approximation results}

We state three approximation results that are needed in the subsequent sections.
The techniques used to prove them are standard.
\begin{lemma} \label{lem:approx-1}
    Let $d \in \N^+$ and $\theta \in (0, \frac12)$.
    If $(F_n)_{n \in \N} \subset \dot{W}^{\theta,2}(\R^d; \C^d)$ is a sequence with $F_n \rightarrow F$ in $L^2_{\loc}(\R^d)$ and $\sup_{n \in \N} [F_n]_{X^{\theta,2}} < \infty$, then $F \in \dot{W}^{\theta,2}(\R^d)$ and
    \begin{align*}
        [F]_{X^{\theta,2}} &\leq \liminf_{n \rightarrow \infty} [F_n]_{X^{\theta,2}} .
    \end{align*}
\end{lemma}
\begin{proof}
    There exists a subsequence $F_{n_k}$ such that $\lim_{k \rightarrow \infty} [F_{n_k}]_{X^{\theta,2}} = \liminf_{n \rightarrow \infty} [F_n]_{X^{\theta,2}}$.
    Since $(F_{n_k})_{k \in \N} \subset \dot{W}^{\theta,2}(\R^d)$ is a bounded sequence in a Hilbert space, there exists a further subsequence, which we also denote by $F_{n_k}$, that converges weakly to some limit $G \in \dot{W}^{\theta,2}(\R^d)$. Then $F_{n_k} \rightarrow F$ in $L^2_{\loc}(\R^d)$ implies $F = G$.
    Lastly, there exists another subsequence, which we still denote by $F_{n_k}$, converging sufficiently quickly on arbitrarily large balls so that, with Chebyshev's inequality and the Borel--Cantelli lemma, one can show $F_{n_k} \rightarrow F$ almost everywhere.  Then lower semicontinuity of the norm yields
    \begin{align*}
        [F]_{X^{\theta,2}} &\leq \liminf_{k \rightarrow \infty}\, [F_{n_k}]_{X^{\theta,2}}
        = \liminf_{n \rightarrow \infty} [F_n]_{X^{\theta,2}} .
    \end{align*}
    Furthermore, Fatou's lemma and the Sobolev inequality imply $F \in L^{\frac{2d}{d-2\theta}}(\R^d)$ and
    \begin{align*}
        \|F\|_{L^{\frac{2d}{d-2\theta}}(\R^d)} &\leq \liminf_{k \rightarrow \infty} \|F_{n_k}\|_{L^{\frac{2d}{d-2\theta}}(\R^d)}
        \lesssim \limsup_{n \rightarrow \infty} [F_n]_{X^{\theta,2}} .
    \end{align*}
\end{proof}

\begin{lemma} \label{lem:approx-2}
    Let $d \in \N^+$, $\alpha \in (0, \frac12)$, and $\nu \in (0, 1)$.
    Let $f \in W^{1+\alpha,2}_{\loc}(\R^d; \C)$ with $\nabla f \in \dot{W}^{\alpha,2}(\R^d; \C^d)$.
    There exist $f_n \in W^{1+\nu,2}_{\loc}(\R^d; \C)$ with $\nabla f_n \in W^{\nu,2}(\R^d; \C^d)$ such that $[\nabla f_n - \nabla f]_{W^{\alpha,2}(\R^d)} \rightarrow 0$ and $f_n \rightarrow f$ in $W^{1,2}_{\loc}(\R^d)$.
\end{lemma}
\begin{proof}
    One can calculate that there exists $c_{d,\alpha} > 0$ such that $[F]_{W^{\alpha,2}(\R^d)} = c_{d,\alpha} \||\xi|^\alpha \ha{F}\|_{L^2}$ for all $F \in C_c^\infty(\R^d; \C^d)$.
    Define $F = \nabla f$ and note that $F \in L^{\frac{2d}{d-2\alpha}}(\R^d; \C^d) \subset \mathcal{S}'(\R^d; \C^d)$ has a well-defined distributional Fourier transform $\ha{F} \in \mathcal{S}'(\R^d; \C^d)$. 
    
    Let $\eta_n \in C_c^\infty(\R^d \setminus \{0\}; [0, 1])$ with $\eta_n(\xi) \rightarrow 1$ for every $\xi \in \R^d \setminus \{0\}$ and define $\ha{F_n}(\xi) = \eta_n(\xi) \ha{F}(\xi)$. Then dominated convergence yields $[F_n - F]_{W^{\alpha,2}} \rightarrow 0$. Since $\ha{F_n}$ is supported in an annulus, we have $F_n \in W^{\nu,2}(\R^d)$. 
    
    Since $F = \nabla f$, the identity $\xi_j \ha{F}_k(\xi) = \xi_k \ha{F}_j(\xi)$ holds distributionally for every $j, k \in \{1, \dots, d\}$, and the same is true for each $F_n$. This implies $\ha{F_n}(\xi) = \xi \frac{\xi \cdot \ha{F_n}(\xi)}{|\xi|^2}$. We now define $\ha{g_n} = - \imu \frac{\xi \cdot \ha{F_n}(\xi)}{|\xi|^2}$ so that $\nabla g_n = F_n$. Since $\supp \ha{F_n}$ is an annulus, this construction is well-defined and $g_n \in W^{1+\nu,2}(\R^d)$. 
    The homogeneous Sobolev inequality now implies $\|F_n - F\|_{L^{\frac{2d}{d-2\alpha}}} \lesssim [F_n - F]_{W^{\alpha,2}} \rightarrow 0$. In particular, $F_n \rightarrow F$ in $L^2_{\loc}(\R^d)$.

    Set $f_n = g_n - \frac{1}{|B_1(0)|} \int_{B_1(0)} (g_n - f) \dd x$. By the classical Poincaré inequality (see e.g.\ \cite[Section 5.8.1, Theorem 1]{Evans2010}) we have $\|f_n - f\|_{L^2(B_1(0))} \lesssim \|\nabla f - \nabla f_n\|_{L^2(B_1(0))} \rightarrow 0$. For any larger ball $B \supseteq B_1(0)$ we again use the Poincaré inequality to estimate $\|f_n - f\|_{L^2(B)} \lesssim \|\nabla f_n - \nabla f\|_{L^2(B)} + \|f_n - f\|_{L^2(B_1(0))} \rightarrow 0$. Therefore $f_n \rightarrow f$ in $L^2_{\loc}(\R^d)$, and hence also in $W^{1,2}_{\loc}(\R^d)$.
\end{proof}
For $d = 1$ we can prove a stronger result.
\begin{lemma} \label{lem:approx-3}
    Let $\theta \in (0, \frac12)$ and $f \in W^{1+\theta,2}_{\loc}(\R; \C)$ with $f' \in \dot{W}^{\theta,2}(\R; \C)$.
    There exist $f_n \in C^\infty(\R; \C)$ with $f_n' \in C_c^\infty(\R)$ such that $[f_n' - f']_{W^{\theta,2}} \rightarrow 0$ and $f_n \rightarrow f$ in $W^{1,2}_{\loc}(\R)$.
\end{lemma}
\begin{proof}
    By definition of $\dot{W}^{\theta,2}(\R; \C)$ there exist $F_n \in C_c^\infty(\R; \C)$ such that $[F_n - f']_{W^{\theta,2}} \rightarrow 0$. We can now define $f_n(x) = f(0) + \int_0^x F_n(y) \dd y$ and argue as in Lemma \ref{lem:approx-2}, using Sobolev and Poincaré inequalities, to prove that $f_n \rightarrow f$ in $W^{1,2}_{\loc}(\R)$.
\end{proof}

\section{Boundedness}

\subsection{Regularity of $\nabla a$} \label{section:a}

In this section we slightly generalize \cite[Théorème 1]{BourdaudMeyer} to functions $\psi \in W^{1+\theta,2}_{\loc}(\R^d; \C)$ with $\nabla \psi \in \dot{W}^{\theta,2}(\R^d; \C^d)$. The proof proceeds as in \cite{BourdaudMeyer} and requires no significant new ideas.

We use the following ``slicing'' result for Sobolev--Slobodeckij seminorms as a substitute for the ``Fubini'' result that is used in \cite[Lemme 3]{BourdaudMeyer}.

\begin{proposition}[{\cite[Theorem 6.47]{Leoni2023}}] \label{prop:slicing}
    Let $d \in \N^+$, $\theta \in (0, 1)$, and $\Omega \subseteq \R^d$ be open. Write $\nu^\bot = \{x \in \R^d: x \cdot \nu = 0\}$ for $\nu \in \R^d$.
    For every $F \in L^2_{\loc}(\Omega; \C^d)$ we have
    \begin{align} \label{eqn:slice}
        [F]_{X^{\theta,2}(\Omega)}^2 &= \frac12 \int_{\mathbb{S}^{d-1}} \int_{\nu^\bot}
        [F(x + \cdot\, \nu) \cdot \nu]_{W^{\theta,2}(\Omega(x, \nu))}^2 \dd \mathcal{H}^{d-1}(x) \dd \mathcal{H}^{d-1}(\nu)
    \end{align}
    where $\Omega(x, \nu) = \{t \in \R: x + t \nu \in \Omega\}$. 
    In particular, a function $F \in L^2_{\loc}(\Omega)$ satisfies $[F]_{X^{\theta,2}(\Omega)} < \infty$ if and only if it admits a representative $\widetilde{F}$ such that $[\widetilde{F}(x + \cdot\, \nu) \cdot\, \nu]_{W^{\theta,2}(\Omega(x, \nu))} < \infty$ for $\mathcal{H}^{d-1}$-a.e.\ $\nu \in \mathbb{S}^{d-1}$ and $\mathcal{H}^{d-1}$-a.e.\ $x \in \nu^\bot$ with $\Omega(x, \nu) \neq \emptyset$, and
    \begin{align*}
        \int_{\mathbb{S}^{d-1}} \int_{\nu^\bot}
        [\widetilde{F}(x + \cdot\, \nu) \cdot \nu]_{W^{\theta,2}(\Omega(x, \nu))}^2 \dd \mathcal{H}^{d-1}(x) \dd \mathcal{H}^{d-1}(\nu) < \infty .
    \end{align*}
\end{proposition}
\begin{proof}
    This is \cite[Theorem 6.47]{Leoni2023}, except the seminorm $[\,\cdot\,]_{W^{\theta,2}(\Omega)}$ has been replaced with $[\,\cdot\,]_{X^{\theta,2}(\Omega)}$. The only change necessary to the proof is that \cite[Lemma 6.48]{Leoni2023} is now applied to $f(x, y) = \left|(F(x) - F(y)) \cdot \frac{x-y}{|x-y|}\right|$.
\end{proof}

We also need a fractional Hardy inequality.
\begin{lemma}
    Let $\theta \in (0, \frac12)$, $I \subset \R$ be an interval, and $f \in W^{\theta,2}(I; \C)$.
    If $I$ is bounded and $\int_I f(x) \dd x = 0$, or if $I$ is unbounded and $f$ has compact support, then we have
    \begin{align} \label{eqn:hardy}
        \int_I \dist(x, \R \setminus I)^{-2\theta} |f(x)|^2 \dd x &\leq C [f]_{W^{\theta,2}(I)}^2 .
    \end{align}
\end{lemma}
\begin{proof}
    For the bounded case we use an inhomogeneous fractional Hardy inequality such as \cite[Theorem 4]{DydaKijaczko2022} and estimate the $L^2$-norm using the fractional Poincaré inequality. For the unbounded case we may assume that $I = [0, \infty)$ and use the fractional Hardy inequality on the half-space \cite{BogdanDyda2011,FrankSeiringer2010}.
\end{proof}

\begin{proposition}[{Variant of \cite[Théorème 1]{BourdaudMeyer}}] \label{prop:1}
    For every $\theta \in (0, \frac12)$ there exists $C > 0$ such that the following holds. Let $d \in \N^+$ and $\psi \in W^{1+\theta,2}_{\loc}(\R^d; \C)$ with $\nabla \psi \in \dot{W}^{\theta,2}(\R^d; \C^d)$. Then $a \in W^{1+\theta,2}_{\loc}(\R^d; \R)$ with $\nabla a \in \dot{W}^{\theta,2}(\R^d; \R^d)$ and
    \begin{align} \label{eqn:00}
        [\nabla a]_{X^{\theta,2}}^2 &\leq C [\nabla \psi]_{X^{\theta,2}}^2 .
    \end{align}
\end{proposition}
\begin{proof}
    \textbf{Estimate.}
    We first note that $\nabla a \in L^{\frac{2d}{d-2\theta}}(\R^d)$ follows directly from $\nabla \psi \in L^{\frac{2d}{d-2\theta}}(\R^d)$.

    It suffices to prove the proposition in the one-dimensional case. To deduce the
    higher-dimensional result, we use Lemma \ref{lem:approx-2} with $\alpha = \nu = \theta$ to approximate $\psi$ by a sequence $\psi_n$ appropriate for the application of Lemma \ref{lem:approx-1}.
    Then by Proposition \ref{prop:slicing} for almost every $\nu\in\mathbb S^{d-1}$ and $x \in \nu^\bot$ we have $\nabla \psi_n(x + \cdot \nu) \cdot \nu \in W^{\theta,2}(\R)$. Applying the one-dimensional estimate, again Proposition \ref{prop:slicing}, and Lemma~\ref{lem:approx-1} yields \eqref{eqn:00}.

    By the formula
    \begin{align} \label{eqn:a-reduction}
        a' &= \frac14 \int_0^{2 \pi} a_r' \dd r ,
    \end{align}
    where $a_r = |\psi_r|$ and $\psi_r = \RE[\psi] \cos r + \IM[\psi] \sin r$, it suffices to consider the case of real-valued $\psi$.

    By Lemma \ref{lem:approx-3} there exist $\psi_n \in C^\infty(\R; \R)$ with $\psi_n' \in C_c^\infty(\R)$ such that $[\psi_n' - \psi']_{W^{\theta,2}(\R)} \rightarrow 0$ and $\psi_n \rightarrow \psi$ in $W^{1,2}_{\loc}(\R)$.
    Then a localized version of Proposition \ref{prop:2}, obtained by using cutoffs, implies also $a_n' \rightarrow a'$ in $L^2_{\loc}(\R)$.
    Therefore Lemma \ref{lem:approx-1} shows that it suffices to prove \eqref{eqn:00} for $\psi_n$ and $a_n$, which we call $\psi$ and $a$ from now on for simplicity.
   
    We decompose
    \begin{align*}
        \R &= \{\psi(x) > 0\} \,\cup\, \{\psi(x) = 0\} \,\cup\, \{\psi(x) < 0\} = \Omega_+ \,\cup\, \Omega_0 \,\cup\, \Omega_-
    \end{align*}
    and note that
    \begin{align*}
        a'\vert_{\Omega_\sigma} &= \sigma \psi'\vert_{\Omega_\sigma} \qquad \forall\, \sigma \in \{+, 0, -\} .
    \end{align*}
    Since $\psi'(x) = a'(x) = 0$ for almost all $x \in \Omega_0$, we have
    \begin{align*}
        \int_{\Omega_0} \int_{\Omega_0} |x-y|^{-1-2\theta} |\Delta_{x-y} a'(x)|^2 \dd y \dd x &= 0
        \\ \int_{\Omega_\pm} \int_{\Omega_0} |x-y|^{-1-2\theta} |\Delta_{x-y} a'(x)|^2 \dd y \dd x 
        &= \int_{\Omega_\pm} \int_{\Omega_0} |x-y|^{-1-2\theta} |a'(x)|^2 \dd y \dd x
        \\ &= \int_{\Omega_\pm} \int_{\Omega_0} |x-y|^{-1-2\theta} |\psi'(x)|^2 \dd y \dd x
        \\ &= \int_{\Omega_\pm} \int_{\Omega_0} |x-y|^{-1-2\theta} |\Delta_{x-y} \psi'(x)|^2 \dd y \dd x \\ \int_{\Omega_\pm} \int_{\Omega_\pm} |x-y|^{-1-2\theta} |\Delta_{x-y} a'(x)|^2 \dd y \dd x 
        &= \int_{\Omega_\pm} \int_{\Omega_\pm} |x-y|^{-1-2\theta} |\Delta_{x-y} \psi'(x)|^2 \dd y \dd x .
    \end{align*}
    It remains to estimate the integral on $\Omega_\pm \times \Omega_\mp$, or more simply $\Omega_+ \times \Omega_-$.
    First, we decompose the open sets
    \begin{align*}
        \Omega_\pm &= \bigcup_{j \in J^\pm} I^\pm_j
    \end{align*}
    into a disjoint union over countable index sets $J^\pm$, where each interval $I^\pm_j$ is a connected component of $\Omega_\pm$. 
    Now we estimate
    \begin{align*}
        & \sum_{j \in J^+} \int_{I^+_j}\! \int_{\Omega_-}\! |x-y|^{-1-2\theta} |\Delta_{x-y} a'(x)|^2 \dd y \dd x
        + \sum_{j \in J^-} \int_{I^-_j}\! \int_{\Omega_+}\! |x-y|^{-1-2\theta} |\Delta_{x-y} a'(x)|^2 \dd y \dd x
        \\ &\leq \sum_{j \in J^+} \int_{I^+_j}\! \frac{4}{\theta} \dist(x, \R \setminus I^+_j)^{- 2 \theta} |\psi'(x)|^2 \dd x
        + \sum_{j \in J^-} \int_{I^-_j} \!\frac{4}{\theta} \dist(x, \R \setminus I^-_j)^{- 2 \theta} |\psi'(x)|^2 \dd x .
    \end{align*}
    If $I^\pm_j$ is a bounded interval, then $\psi$ vanishes at its endpoints and hence we can apply the homogeneous fractional Hardy inequality \eqref{eqn:hardy}.
    If $I^\pm_j$ is an unbounded interval, then we use $\psi' \in C_c^\infty(\R)$ to apply the homogeneous fractional Hardy inequality instead.
\end{proof}

\subsection{Regularity of $\Lambda$} \label{section:lambda}

We define for $\gamma \in (0, \frac14)$ a set of approximants
\begin{align*}
    \mathcal{A}^{1+\gamma} &= \Big\{f \in W^{1+\gamma,2}_{\loc}(\R^d; \C) \text{ with } \nabla f \in W^{\gamma,2}(\R^d; \C^d) 
    \\ &\qquad\qquad \text{ and } w_\gamma R_h f(x) \conj{\Delta_h f(x)} \in L^1_h L^1_x(\R^{d+d}; \C) \Big\} .
\end{align*}
The following lemma provides a sufficient condition for membership in the set of approximants.
\begin{lemma} \label{lem:20}
    Let $\gamma \in (0, \frac14)$ and $\beta \in (0, 1 - 2 \gamma)$. If $f \in W^{1+\gamma,2}_{\loc}(\R^d; \C)$ with $\nabla f \in W^{2\gamma+\beta,2}(\R^d; \C^d)$, then $f \in \mathcal{A}^{1+\gamma}$.
\end{lemma}
\begin{proof}
    By separating the cases $|h| < 1$ and $|h| \geq 1$, we find that $\nabla f \in W^{\gamma,2}(\R^d)$.
    Since 
    \begin{align*}
        \Delta_h f(x) &= \int_0^1 \nabla f(x-rh) \cdot h \dd r
        \\ R_h f(x) &= - \int_0^1 \Delta_{rh} \nabla f(x) \cdot h \dd r
    \end{align*}
    we can estimate
    \begin{align*}
        \mathds{1}_{\{|h| \geq 1\}} w_\gamma \|R_h f \Delta_h f\|_{L^1} &\leq \mathds{1}_{\{|h| \geq 1\}} |h|^{-d-2\gamma} 4 \|\nabla f\|_{L^2}^2 \in L^1_h(\R^d) .
    \end{align*}
    On the other hand, we know from \eqref{eqn:9} that
    \begin{align*}
        |h|^{-d-2-4\gamma-2\beta} \|R_h f(x)\|_{L^2}^2 \in L^1_h(\R^d) ,
    \end{align*}
    so we are finished with Cauchy--Schwarz and the observation that
    \begin{align*}
        \mathds{1}_{\{|h| < 1\}} |h|^{-d-2+2\beta} \|\Delta_h f(x)\|_{L^2}^2 &\leq \mathds{1}_{\{|h| < 1\}} |h|^{-d+2\beta} 4 \|\nabla f\|_{L^2}^2 \in L^1_h(\R^d) .
    \end{align*}
\end{proof}

The following crucial lemma allows us to bootstrap a growth bound of the form $u(\epsilon) \lesssim \epsilon^{2\gamma-2\theta}$ to a uniform bound $u(\epsilon) \lesssim 1$.
\begin{lemma} \label{lem:bootstrap}
    Let $\theta \in (0, 1)$ and $u: (0, 1) \rightarrow [0, \infty)$ be measurable, satisfying the growth bound $u(\epsilon) \leq c \epsilon^{2\gamma-2\theta}$ for some $c > 0$ and $\gamma > 0$. 
    Define
    \begin{align*}
        U(\epsilon) &= \int_0^1 r^{2\theta} u(\epsilon r) \dd r = \epsilon^{-1-2\theta} \int_0^\epsilon r^{2\theta} u(r) \dd r
    \end{align*}
    and suppose that
    \begin{align*}
        \sup_{\epsilon \in (0, 1)} (u(\epsilon) - U(\epsilon)) &\leq C
    \end{align*}
    for some $C > 0$. Then
    \begin{align*}
        \sup_{\epsilon \in (0, 1)} u(\epsilon) &\leq \frac{1+2\theta}{2\theta} C .
    \end{align*}
\end{lemma}
\begin{proof}
    Define $B(\epsilon) = U(\epsilon) - \frac{C}{2\theta}$ and note that it is locally absolutely continuous. We have
    \begin{align*}
        \frac{\dd}{\dd \epsilon} (\epsilon^{2 \theta} B(\epsilon)) &= \epsilon^{2 \theta - 1} (u(\epsilon) - U(\epsilon) - C) \leq 0
    \end{align*}
    almost everywhere. The growth bound for $u$ implies $U(\epsilon) \leq \frac{c}{1+2\gamma} \epsilon^{2\gamma-2\theta}$, and in particular $\epsilon^{2 \theta} B(\epsilon) \rightarrow 0$ as $\epsilon \searrow 0$. We conclude that $\sup_{\epsilon \in (0, 1)} B(\epsilon) \leq 0$ and hence $\sup_{\epsilon \in (0, 1)} u(\epsilon) \leq \sup_{\epsilon \in (0, 1)} U(\epsilon) + C \leq \frac{C}{2 \theta} + C \leq \frac{1+2\theta}{2\theta} C$.
\end{proof}

\begin{theorem} \label{thm:1}
    Let $d \in \N^+$ and $\theta \in (0, 1)$.
    Suppose that $\psi \in W^{1,2}_{\loc}(\R^d; \C)$ with $\nabla \psi \in L^2(\R^d; \C^d)$ and let $P \in \mathcal{P}(\psi)$ be a polar factor.
    Consider the following statements:
    \begin{enumerate}[(i)]
        \item $\psi \in W^{1+\theta,2}_{\loc}(\R^d)$ and $\nabla \psi \in W^{\theta,2}(\R^d)$.
        \item $a \in W^{1+\theta,2}_{\loc}(\R^d)$ and $\nabla a \in W^{\theta,2}(\R^d)$.
        \item $\Lambda, P \Lambda \in W^{\theta,2}(\R^d; \C^d)$.
    \end{enumerate}
    If $\theta < \frac12$ then (i) $\Longleftrightarrow$ ((ii) {\normalfont\texttt{and}} (iii)).
    Next, if $\theta \geq \frac12$ then ((i) {\normalfont\texttt{and}} (ii)) $\Longrightarrow$ (iii).
    Whenever (i), (ii), and (iii) hold, we have the polarization identity
    \begin{align} \label{eqn:formula}
        & (1 - 2 \theta) [\nabla a]_{X^{\theta,2}}^2 + (1 + 2 \theta) [\Lambda]_{X^{\theta,2}}^2 - (1 - 2 \theta) [\nabla \psi]_{X^{\theta,2}}^2
        \\\nonumber &= \int_{\R^d} w_\theta \int_{\R^d} (2 + 4 \theta) \RE\Big[ \conj{P(x)} \Big( R_h \psi(x) \Delta_h \nabla a(x) \cdot h \\\nonumber &\hspace{45mm} - R_h a(x) \Delta_h \nabla \psi(x) \cdot h \Big) \Big] \dd x \dd h
        \\\nonumber &- \int_{\R^d} w_\theta \int_{\R^d} \IM\Big[ \conj{P(x)} R_h \psi(x) \Big( \big(8 \theta - (2 - 4 \theta) P(x) \conj{P(x-h)}\big) \Delta_h \Lambda(x) \cdot h 
        \\\nonumber &\hspace{45mm} + (2 - 4 \theta) \conj{P(x-h)} \Delta_h (P \Lambda)(x) \cdot h \Big) \Big] \dd x \dd h .
    \end{align}
\end{theorem}
\begin{proof}
    \textbf{Direction ``$\Longrightarrow$''}.
    Since $a, \psi \in W^{1+\theta,2}_{\loc}({\R^d})$, using Proposition \ref{prop:1} when $\theta < \frac12$, we can use the fundamental theorem of calculus to write
    \begin{align}
        \label{eqn:ftc-a} a(x) - a(x-h) &= \int_0^1 \nabla a(x-rh) \cdot h \dd r 
        \\ \label{eqn:ftc-q} \psi(x) - \psi(x-h) &= \int_0^1 \nabla \psi(x-rh) \cdot h \dd r .
    \end{align}
    We use this in the estimates \eqref{eqn:9}, \eqref{eqn:10}, and \eqref{eqn:10-eps} above, as well as in \eqref{eqn:70} below. 
    These are the only places where the fact that $\nabla a$ and $\nabla \psi$ are the weak derivatives of $a$ and $\psi$ is used.

    In the proof below we also make use of the crucial nonlinear identities \eqref{eqn:id1}--\eqref{eqn:id14}. Except for \eqref{eqn:id6}, these follow \textit{purely algebraically} from
    \begin{align}   
        \label{eqn:polar} P \conj{P} &= 1
        \\ \label{eqn:def-a} a &= \conj{P} \psi
        \\ \label{eqn:def-nabla-a} \nabla a &= \RE\left[ \conj{P} \nabla \psi \right]
        \\ \label{eqn:def-Lambda} \Lambda &= \IM\left[ \conj{P} \nabla \psi \right] .
    \end{align}
    We have verified these identities via symbolic calculations using the Python library ``SymPy'' \cite{sympy}.

    \textbf{Step 1: essential nonlinear identities.}
    All identities below must be understood in the sense of holding for almost all $(x, h) \in \R^{d+d}$.

    We introduce the auxiliary variables
    \begin{align*}
        Q_h(x) &= P(x) \conj{P(x-h)} ,
    \end{align*}
    as well as the shorthand notations
    \begin{align*}
        B_h^\pm \psi(x) &= R_h \psi(x) \conj{P(x)} \pm \conj{R_{-h} \psi(x-h)} P(x-h)
        \\ B_h^\pm a(x) &= R_h a(x) \pm R_{-h} a(x-h) .
    \end{align*}
    We also define, for general functions $L = L(x, h)$ and $r \in \R$, the operators
    \begin{align*}
        S_r L(x, h) &= \frac{L(x, h) + L(x - r h, - h)}{2} 
        \\ T_r L(x, h) &= L(x, r h) .
    \end{align*}
    
    Define
    \begin{align*}
        L_1(x, h) &= (\Delta_h \Lambda(x) \cdot h)^2 .
    \end{align*}
	The following identity is crucial:
	\begin{align} \label{eqn:id1}
		\tag{\texttt{id1}} (R_h \psi(x)) \conj{P(x)} &= R_h a(x) + a(x-h) (1 - \RE[Q_h(x)])
		\\\nonumber &+ \imu (a(x-h) \IM[Q_h(x)] - \Lambda(x) \cdot h) .
	\end{align}    
	It implies
	\begin{align*}
		\Delta_h \Lambda(x) \cdot h &= - \IM\left[B_h^- \psi(x)\right] - \Delta_h a(x) \IM[Q_h(x)] .
	\end{align*}
    We note also that
	\begin{align*}
		2 \IM[Q_h(x)]^2 &= 4 (1 - \RE[Q_h(x)]) - 2 (1 - \RE[Q_h(x)])^2
	\end{align*}
    and
	\begin{align} \label{eqn:id2} 
        \tag{\texttt{id2}} a(x) (1 - \RE[Q_h(x)]) &= \RE\left[(R_{-h} \psi(x-h)) \conj{P(x-h)}\right] - R_{-h} a(x-h)
        \\\nonumber \Delta_h a(x) (1 - \RE[Q_h(x)]) &= - \RE\left[B_h^- \psi(x)\right] + B_h^- a(x) 
    \end{align}
    Combined with \eqref{eqn:9}, these identities allow us to control the terms denoted by $G_j$, $j \in \{1, 3, 5, 6, 7, 8\}$ below. On the other hand, \eqref{eqn:10} and \eqref{eqn:10-eps} allow us to control $G_2$ and $G_4$ below. As such, the terms denoted by $G_j$ below are well-controlled.

    Observe also that
	\begin{align} \label{eqn:id3}
		\tag{\texttt{id3}}
		\RE\left[\Delta_h \psi(x) \conj{P(x)}\right] &= \Delta_h a(x) + a(x-h) (1 - \RE[Q_h(x)]) .
	\end{align}    
	Using this, one can calculate that
	\begin{align} \label{eqn:id4}
        \tag{\texttt{id4}}
		L_1 &= - J_1 + G_1 + G_2 + S_1 L_2
	\end{align}
    where
    \begin{align*}
        J_1(x, h) &= \left( \Delta_h \Lambda(x) \cdot h + 2 \IM\left[B_h^- \psi(x)\right] \right)^2
        \\ G_1(x, h) &= 2 \IM\left[B_h^- \psi(x)\right]^2 
        - 2 (\Delta_h a(x))^2 (1 - \RE[Q_h(x)])^2
        \\ G_2(x, h) &= 4 B_h^- a(x) \Delta_h a(x)
        \\ L_2(x, h) &= 8 \RE\left[- (R_h \psi(x)) \conj{P(x)}\right] \Delta_h a(x) .
    \end{align*}    
    Here and below the terms denoted by $L_j$ are key terms that remain to be estimated.
    
	We calculate
    \begin{align} \label{eqn:id5}
        \tag{\texttt{id5}}
        L_2 &= G_3 + G_4 + L_3
    \end{align}
    where
    \begin{align*}
        G_3(x, h) &= 4 |R_h \psi(x)|^2 - 8 \RE\left[(R_h \psi(x)) \conj{P(x)}\right] R_h a(x)
        \\ G_4(x, h) &= - 4 \RE\left[(R_h \psi(x)) \conj{\Delta_h \psi(x)}\right] 
        \\ L_3(x, h) &= - 4 \RE\left[(R_h \psi(x)) (\nabla \psi(x) \cdot h) \conj{P(x)}^2\right] .
    \end{align*}    
    Now we write
    \begin{align}  
    	\label{eqn:70} L_3(x, h) &= 4 \int_0^1 \RE\left[((\nabla \psi(x) - \nabla \psi(x-rh)) \cdot h) (\nabla \psi(x) \cdot h) \conj{P(x)}^2\right] \dd r
    	\\\nonumber &= 2 \int_0^1 \RE\left[((\nabla \psi(x) - \nabla \psi(x-rh)) \cdot h)^2 \conj{P(x)}^2\right] \dd r
    	\\\nonumber &+ 2 \int_0^1 \RE\left[((\nabla \psi(x) \cdot h)^2 - (\nabla \psi(x-rh) \cdot h)^2) \conj{P(x)}^2\right] \dd r
        \\\label{eqn:id6} \tag{\texttt{id6}} &= \int_0^1 r^{-2} T_r L_4(x, h) \dd r  
    \end{align}
    where
    \begin{align*}
        L_4(x, h) &= 2 \RE\left[((\nabla \psi(x) \cdot h)^2 - (\nabla \psi(x-h) \cdot h)^2) \conj{P(x)}^2\right] 
        \\ &+ 2 \RE\left[(R_h \psi(x) + R_{-h} \psi(x-h))^2 \conj{P(x)}^2\right] .
    \end{align*}

    A remarkable identity states
    \begin{align} \label{eqn:id7}
        \tag{\texttt{id7}}
        S_1 L_4 &= G_5 + J_1 - \RE\left[Q_h(x)\right] L_1
    \end{align}
    where
    \begin{align*}
        G_5(x, h) &= 2 \RE\left[(B_h^- \psi(x))^2\right] - \left|B_h^+ \psi(x) - B_h^+ a(x)\right|^2
        \\ &+ ((B_h^+ a(x))^2 - (B_h^- a(x))^2) (1 + \RE\left[Q_h(x)\right])
        - \RE\left[(B_h^+ \psi(x) - B_h^+ a(x))^2 Q_h(x)\right] .
    \end{align*}
    
    In summary, we have
    \begin{align}
        \label{eqn:final-1} L_1 &= - J_1 + G_1 + G_2 + S_1 L_2
        \\\label{eqn:final-2}  L_2 &= G_3 + G_4 + L_3
        \\\label{eqn:final-3}  L_3 &= \int_0^1 r^{-2} T_r L_4 \dd r
        \\\label{eqn:final-4}  S_1 L_4 &= G_5 + J_1 - \RE\left[Q_h(x)\right] L_1 .
    \end{align}

    We now list a number of further, less essential identities.
    Since later we apply \eqref{eqn:10-eps} to $G_2$ and $G_4$, let us define   
    \begin{align*}
        G_{2,1}(x, h) &= 4 R_h a(x)^2 + 4 R_{-h} a(x-h)^2
        \\ G_{2,2}(x, h) &= - 4 (\Delta_h \nabla a(x) \cdot h)^2
        \\ G_{4,1}(x, h) &= - 2 |R_h \psi(x)|^2 - 2 |R_{-h} \psi(x-h)|^2
        \\ G_{4,2}(x, h) &= 2 |\Delta_h \nabla \psi(x) \cdot h|^2 .
    \end{align*}
    We have the identities
    \begin{align*}
        \label{eqn:id8} \tag{\texttt{id8}} 
        G_1 + G_{2,1} + S_1 G_3 + G_{4,1} + G_{2,2} + G_{4,2} + G_5 &= G_6
        \\ \label{eqn:id9} \tag{\texttt{id9}}
        G_1 + G_{2,1} + S_1 G_3 + G_{4,1} &= G_7
    \end{align*}
    where
    \begin{align*}
        G_6(x, h) 
        &= \frac12 \left|\Delta_h \nabla \psi(x) \cdot h - \conj{\Delta_h \nabla \psi(x) \cdot h} P(x) P(x-h)\right|^2
        \\ &+ 2 B_h^+ a(x) \RE\left[B_h^+ \psi(x) Q_h(x)\right]
        + 2 B_h^- a(x) \RE\left[B_h^- \psi(x)\right]
        \\ &- 2 (B_h^+ a(x))^2 - (1 + \RE\left[Q_h(x)\right]) (B_h^- a(x))^2
        \\ G_7(x, h) &= - 2 \RE\left[(B_h^- \psi(x))^2\right] + 2 B_h^- a(x) \RE\left[B_h^- \psi(x)\right] - 2 B_h^+ a(x) \RE\left[B_h^+ \psi(x)\right]
        \\ &+ 2 (B_h^+ a(x))^2  .
    \end{align*}
    Finally, we write
    \begin{align*}
        G_6 + 2 \theta G_7 - 8 \theta \IM\left[B_h^- \psi(x)\right]^2 + (1 - \RE\left[Q_h(x)\right]) L_1 &= S_1 G_8
    \end{align*}
    where
    \begin{align*}
        G_8(x, h) &= (2 - 4 \theta) |R_h \psi(x) \conj{P(x)} - \conj{R_{-h} \psi(x-h)} P(x-h)|^2 - (2 - 4 \theta) (\Delta_h \nabla a(x) \cdot h)^2
        \\ &+ 8 R_{-h} a(x-h) \RE\left[R_h \psi(x) \conj{P(x-h) - 2 \theta P(x)}\right] .
    \end{align*}
    Here we have used again \eqref{eqn:id1}--\eqref{eqn:id2} to write $(1 - \RE\left[Q_h(x)\right]) L_1$ in terms of $a$ and $\psi$.
    
    \textbf{Step 2: formally rewriting the integral.}
    We shall now perform a series of formal calculations with integrals. 
    Afterwards, we justify them by verifying that each term is indeed integrable. 
    Recall that for $\epsilon \geq 0$ we write $A_\epsilon = \{h \in \R^d: |h| \in (\epsilon, \epsilon^{-1})\}$ if $\epsilon > 0$ and $A_\epsilon = \R^d$ if $\epsilon = 0$. Recall furthermore the shorthand notation 
    \begin{align*}
        w_\theta &= w_\theta(h) = |h|^{-d-2-2\theta} .
    \end{align*}
    
    We have
    \begin{align*}
        0 &= \int_{A_\epsilon} w_\theta \int_{\R^d} - L_1 - J_1 + G_1 + G_2 + S_1 G_3 + S_1 G_4 
        \\ &\qquad\qquad\qquad\qquad\qquad + S_1 \int_0^1 r^{-2} T_r L_4 \dd r \dd x \dd h .
    \end{align*}
    Using that $S_1$ and $S_r$ are symmetries of the integral and the identity $S_r T_r = T_r S_1$, we can write
    \begin{align*}
        & \int_{A_\epsilon} w_\theta \int_{\R^d} L_1 + J_1 \dd x \dd h
        \\ &= \int_0^1 \int_{A_\epsilon} w_\theta \int_{\R^d} G_1 + G_2 + S_1 G_3 + S_1 G_4
        + r^{-2} T_r L_4 \dd x \dd h \dd r
        \\ &= \int_0^1 \int_{A_\epsilon} w_\theta \int_{\R^d} G_1 + S_1 G_2 + S_1 G_3 + S_1 G_4
        + r^{-2} S_r T_r L_4 \dd x \dd h \dd r 
        \\ &= \int_0^1 \int_{A_\epsilon} w_\theta \int_{\R^d} G_1 + S_1 G_2 + S_1 G_3 + S_1 G_4
        + r^{-2} T_r S_1 L_4 \dd x \dd h \dd r 
        \\ &= \int_{A_\epsilon} w_\theta \int_{\R^d} G_1 + S_1 G_2 + S_1 G_3 + S_1 G_4 \dd x \dd h
        + \int_0^1 \int_{A_\epsilon} w_\theta \int_{\R^d} r^{-2} T_r G_5 \dd x \dd h \dd r
        \\ &+ \int_0^1 \int_{A_\epsilon} w_\theta \int_{\R^d} r^{-2} T_r \left(J_1 - \RE\left[Q_h(x)\right] L_1\right) \dd x \dd h \dd r.
    \end{align*}
    For any function $L: \R^{d+d} \rightarrow \C$ we note the change of variables $h \mapsto \frac{h}{r}$:
    \begin{align} \label{eqn:Tr}
        \int_{A_\epsilon} w_\theta \int_{\R^d} r^{-2} T_r L(x, h) \dd x \dd h &= \int_{r A_\epsilon} r^{2\theta} w_\theta \int_{\R^d} L(x, h) \dd x \dd h .
    \end{align}
    We also note that applying \eqref{eqn:10-eps} yields
    \begin{align*}
        \int_{A_\epsilon} w_\theta \int_{\R^d} S_1 G_2 \dd x \dd h
        &= \int_{A_\epsilon} w_\theta \int_{\R^d} G_{2,1} \dd x \dd h
        + \int_0^1 r^{2 \theta} \int_{r A_\epsilon} w_\theta \int_{\R^d} G_{2,2} \dd x \dd h \dd r ,
    \end{align*}
    and
    \begin{align*}
        \int_{A_\epsilon} w_\theta \int_{\R^d} S_1 G_4 \dd x \dd h
        &= \int_{A_\epsilon} w_\theta \int_{\R^d} G_{4,1} \dd x \dd h
        + \int_0^1 r^{2 \theta} \int_{r A_\epsilon} w_\theta \int_{\R^d} G_{4,2} \dd x \dd h \dd r .
    \end{align*}
    Therefore
    \begin{align} \label{eqn:99}
        & \int_{A_\epsilon} w_\theta \int_{\R^d} L_1 + J_1 \dd x \dd h
        - \int_0^1 r^{2 \theta} \int_{r A_\epsilon} w_\theta \int_{\R^d} - \RE\left[Q_h(x)\right] L_1 + J_1 \dd x \dd h \dd r
        \\\nonumber &= \int_{A_\epsilon} w_\theta \int_{\R^d} G_1 + G_{2,1} + S_1 G_3 + G_{4,1} \dd x \dd h
        \\\nonumber &+ \int_0^1 r^{2 \theta} \int_{r A_\epsilon} w_\theta \int_{\R^d} G_{2,2} + G_{4,2} + G_5 \dd x \dd h \dd r .
    \end{align}
    Formally taking the limit $\epsilon \searrow 0$ and multiplying both sides by $(1 + 2 \theta)$ yields
    \begin{align} \label{eqn:99-limit}
        & \int_{\R^d} w_\theta \int_{\R^d} (1 + 2 \theta + \RE\left[Q_h(x)\right]) L_1 + 2 \theta J_1 \dd x \dd h
        \\\nonumber &= \int_{\R^d} w_\theta \int_{\R^d} (1 + 2 \theta) (G_1 + G_{2,1} + S_1 G_3 + G_{4,1}) + G_{2,2} + G_{4,2} + G_5 \dd x \dd h 
        \\\nonumber &= \int_{\R^d} w_\theta \int_{\R^d} G_6 + 2 \theta G_7 \dd x \dd h .
    \end{align}
    Substituting the definition of $J_1$ lets us write
    \begin{align*}
        & \int_{\R^d} w_\theta \int_{\R^d} (2 + 4 \theta) L_1 + 8 \theta (\Delta_h \Lambda(x) \cdot h) \IM\left[B_h^- \psi(x)\right] \dd x \dd h
        \\ &= \int_{\R^d} w_\theta \int_{\R^d} G_6 + 2 \theta G_7 - 8 \theta \IM\left[B_h^- \psi(x)\right]^2 + (1 - \RE\left[Q_h(x)\right]) L_1 \dd x \dd h
        \\ &= \int_{\R^d} w_\theta \int_{\R^d} S_1 G_8 \dd x \dd h .
    \end{align*}
    In order to obtain a cleaner formula, we use 
    \begin{align*}
        8 \theta (\Delta_h \Lambda(x) \cdot h) \IM\left[B_h^- \psi(x)\right] &= S_1\left(8 (\Delta_h \Lambda(x) \cdot h) \IM\left[R_h \psi(x) \conj{P(x)} 2 \theta\right]\right) .
    \end{align*}
    We have proved that
    \begin{align} \label{eqn:formula-pre}
        & (1 + 2 \theta) [\Lambda]_{X^{\theta,2}}^2 + (1 - 2 \theta) [\nabla a]_{X^{\theta,2}}^2
        \\\nonumber &= (1 - 2 \theta) \int_{\R^d} w_\theta \int_{\R^d} \left|R_h \psi(x) \conj{P(x)} - \conj{R_{-h} \psi(x-h)} P(x-h)\right|^2 \dd x \dd h
        \\\nonumber &- 4 \int_{\R^d} w_\theta \int_{\R^d} R_{-h} a(x-h) \RE\left[(2 \theta - Q_h(x)) R_h \psi(x) \conj{P(x)}\right] \dd x \dd h
        \\\nonumber &- 4 \int_{\R^d} w_\theta \int_{\R^d} (\Delta_h \Lambda(x) \cdot h) \IM\left[2 \theta R_h \psi(x) \conj{P(x)} \right] \dd x \dd h .
    \end{align}
    This is equivalent to 
    \begin{align} \label{eqn:formula-pre-2}
        & (1 - 2 \theta) [\nabla a]_{X^{\theta,2}}^2 + (1 + 2 \theta) [\Lambda]_{X^{\theta,2}}^2
        \\\nonumber &= (1 - 2 \theta) \int_{\R^d} w_\theta \int_{\R^d} |B_h^- \psi(x)|^2 \dd x \dd h
        \\\nonumber &+ \int_{\R^d} w_\theta \int_{\R^d} \RE\left[ (Q_h(x) - 2 \theta) (B_h^+ a(x) B_h^+ \psi(x) - B_h^- a(x) B_h^- \psi(x)) \right] \dd x \dd h
        \\\nonumber &- 4 \theta \int_{\R^d} w_\theta \int_{\R^d}  \Delta_h \Lambda(x) \cdot h \IM\left[B_h^- \psi(x)\right] \dd x \dd h .
    \end{align}
    After subtracting $(1 - 2 \theta) [\nabla \psi]_{X^{\theta,2}}^2$ from both sides and using \eqref{eqn:9}--\eqref{eqn:9-alt}, we have
    \begin{align} \label{eqn:formula-pre-3}
        (1 - 2 \theta) [\nabla a]_{X^{\theta,2}}^2 + (1 + 2 \theta) [\Lambda]_{X^{\theta,2}}^2 - (1 - 2 \theta) [\nabla \psi]_{X^{\theta,2}}^2
        &= \int_{\R^d} w_\theta \int_{\R^d} K_1(x, h) \dd x \dd h
    \end{align}
    where
    \begin{align*}
        K_1(x, h) &= - (2 - 4 \theta) \RE\left[ \conj{P(x)} R_h \psi(x) \left( R_{-h} \psi(x-h) \conj{P(x-h)} + \conj{R_{-h} \psi(x-h)} P(x) \right) \right]
        \\\nonumber &+ \RE\left[ (Q_h(x) - 2 \theta) (B_h^+ a(x) B_h^+ \psi(x) - B_h^- a(x) B_h^- \psi(x)) \right]
        \\\nonumber &- 4 \theta \Delta_h \Lambda(x) \cdot h \IM\left[B_h^- \psi(x)\right] .
    \end{align*}
    Note that
    \begin{align}  \label{eqn:id10} \tag{\texttt{id10}} 
        \imu (1 - Q_h(x)) \Lambda(x-h) \cdot h &= R_{-h} \psi(x-h) \conj{P(x-h)} + \conj{R_{-h} \psi(x-h)} P(x) 
        \\\nonumber &- (1 + Q_h(x)) R_{-h} a(x-h)
    \end{align}
    is a term which can be controlled using \eqref{eqn:9}.
    We conclude the proof of \eqref{eqn:formula} with the identity
    \begin{align}  \label{eqn:id11} \tag{\texttt{id11}} 
        K_1(x, h) &= S_1 K_2 ,
    \end{align}
    where
    \begin{align*}
        K_2(x, h) &= (2+4\theta) \RE\left[\conj{P(x)} \left( R_h \psi(x) \Delta_h \nabla a(x) \cdot h - R_h a(x) \Delta_h \nabla \psi(x) \cdot h \right) \right]
        \\ &- \IM\left[ \conj{P(x)} R_h \psi(x) \left( 8 \theta \Delta_h \Lambda(x) \cdot h - (2 - 4 \theta) (1 - Q_h(x)) \Lambda(x-h) \cdot h \right) \right] .
    \end{align*}

    \textbf{Step 3: integrability.}
    Note that \eqref{eqn:99} implies
    \begin{align} \label{eqn:100}
        & \int_{A_\epsilon} w_\theta \int_{\R^d} L_1 + J_1 \dd x \dd h
        - \int_0^1 r^{2\theta} \int_{A_{r \epsilon}} w_\theta \int_{\R^d} L_1 + J_1 \dd x \dd h \dd r
        \\\nonumber &\leq \int_{A_\epsilon} w_\theta \int_{\R^d} L_1 + J_1 \dd x \dd h
        - \int_0^1 r^{2\theta} \int_{r A_\epsilon} w_\theta \int_{\R^d} L_1 + J_1 \dd x \dd h \dd r
        \\\nonumber &\leq \int_{A_\epsilon} w_\theta \int_{\R^d} L_1 + J_1 \dd x \dd h
        - \int_0^1 r^{2\theta} \int_{r A_\epsilon} w_\theta \int_{\R^d} - \RE\left[Q_h(x)\right] L_1 + J_1 \dd x \dd h \dd r
        \\\nonumber &= \int_{A_\epsilon} w_\theta \int_{\R^d} G_7 \dd x \dd h
        + \int_0^1 r^{2 \theta} \int_{r A_\epsilon} w_\theta \int_{\R^d} G_6 - G_7 \dd x \dd h \dd r .
    \end{align}
    We now justify the existence of these integrals, and thereby the validity of the calculations above.
    Once we have proven that $\Lambda \in W^{\theta,2}(\R^d)$, taking the limit $\epsilon \searrow 0$ to go from \eqref{eqn:99} to \eqref{eqn:99-limit} is justified by dominated convergence.

    Let $\gamma \in (0, \frac{\theta}{4})$ and $\beta \in (0, \frac{\theta}{2} - 2 \gamma)$. Then $\psi \in (W^{1+2\gamma+\beta,2}_{\loc} \,\cap\, W^{1+\theta,2}_{\loc})(\R^d)$ with $\nabla \psi \in (W^{2\gamma+\beta,2} \,\cap\, W^{\theta,2})(\R^d)$.
    Since $2 \gamma + \beta < \frac12$, Proposition \ref{prop:1} implies for $\theta < \frac12$ that $a, \psi \in (W^{1+2\gamma+\beta,2}_{\loc} \,\cap\, W^{1+\theta,2}_{\loc})(\R^d)$ and $\nabla a, \nabla \psi \in (W^{2\gamma+\beta,2} \cap W^{\theta,2})(\R^d)$. If $\theta \geq \frac12$, then this additional regularity has been assumed.
    The conditions of Lemma \ref{lem:20} are satisfied, which implies $a, \psi \in \mathcal{A}^{1+\gamma}$. 
    This yields
    \begin{align}
        \label{eqn:ass-9-1} w_\gamma R_h \psi(x) \conj{\Delta_h \psi(x)} &\in L^1_h L^1_x(\R^{d+d})
        \\\label{eqn:ass-9-2} w_\gamma R_h a(x) \Delta_h a(x) &\in L^1_h L^1_x(\R^{d+d}) .
    \end{align}
    In particular, both $\psi$ and $a$ admit the estimates \eqref{eqn:9} with $\alpha \in \{\gamma, \theta\}$, \eqref{eqn:10} with $\alpha = \gamma$, and \eqref{eqn:10-eps} with $\alpha = \theta$ and $\nu = \gamma$.
    Note also that $\nabla a, \Lambda \in L^2(\R^d)$ follows directly from $\nabla \psi \in L^2(\R^d)$. We therefore aim to prove
    \begin{align*}
        [\Lambda]_{X^{\theta,2}}^2 &= \sup_{\epsilon \in (0, 1)} \int_{A_\epsilon} w_\theta \int_{\R^d} L_1 \dd x \dd h < \infty .
    \end{align*}
    
    \textbf{Step 3.1: the simple argument for $\Lambda \in W^{\gamma,2}(\R^d)$.}
    Here we use the calculations from Step 2 leading up to \eqref{eqn:100} with $\theta$ replaced by $\gamma$, and $\epsilon = 0$ so that $A_\epsilon = \R^d$.

    Note that by \eqref{eqn:id2} and \eqref{eqn:9} the functions
    \begin{align*}
        w_\gamma G \text{ for } G \in \{G_1, G_{2,1}, G_{2,2}, G_3, G_{4,1}, G_{4,2}, G_5, G_6, G_7, G_8\} 
    \end{align*}
    are integrable on $\R^{d+d}$.
    For $w_\gamma G_2$ and $w_\gamma G_4$ \eqref{eqn:ass-9-1}--\eqref{eqn:ass-9-2} yield the necessary integrability.
    
    For any $L: \R^{d+d} \rightarrow \C$ the function $w_\gamma r^{-2} T_r L$ is integrable over $(x, h, r) \in \R^{d+d} \times [0, 1]$ if and only if $w_\gamma L$ is integrable on $\R^{d+d}$. We use this for $w_\gamma G_5$.
    
    Furthermore, \eqref{eqn:ass-9-1}--\eqref{eqn:ass-9-2} with (i) $\Longrightarrow$ (iii) of Lemma \ref{lem:equiv} imply that $w_\gamma r^{-2} T_r L_4$, written using the right-hand side of \eqref{eqn:70}, is integrable on $\R^{d+d} \times [0, 1]$. In particular $w_\gamma L_4$ is integrable on $\R^{d+d}$, and then also $w_\gamma L_3$, $w_\gamma L_2$, and $w_\gamma L_1 + w_\gamma J_1$ by the identities \eqref{eqn:final-1}--\eqref{eqn:final-3}.
    But $w_\gamma L_1$ and $w_\gamma J_1$ are nonnegative, so both of them are integrable on $\R^{d+d}$, and furthermore $w_\gamma r^{-2} T_r L_1$ and $w_\gamma r^{-2} T_r J_1$ are each integrable on $\R^{d+d} \times [0, 1]$.
    
    Now the above calculations are well-defined, and \eqref{eqn:100} together with \eqref{eqn:id2}, \eqref{eqn:9}, and \eqref{eqn:10} yield an estimate of the form $[\Lambda]_{X^{\gamma,2}}^2 \leq C_1(\gamma) [\nabla \psi]_{X^{\gamma,2}}^2 + C_2(\gamma) [\nabla a]_{X^{\gamma,2}}^2$, and in particular $\Lambda \in W^{\gamma,2}(\R^d)$.

    \textbf{Step 3.2: the delicate argument for $\Lambda \in W^{\theta,2}(\R^d)$.}
    Now we use the above calculations from Step 2 leading up to \eqref{eqn:100} with $\epsilon > 0$ instead of $\epsilon = 0$.
    The integrability of all terms with weight $w_\theta$ instead of $w_\gamma$ still holds, as $w_\theta(h) \leq \epsilon^{2(\gamma-\theta)} w_\gamma(h)$ for $h \in A_\epsilon$.
    An inspection of $L_1 + J_1$, applying \eqref{eqn:9} and Step 3.1, reveals that
    \begin{align*}
        u(\epsilon) &= \int_{A_\epsilon} |h|^{2\gamma-2\theta} w_\gamma \int_{\R^d} L_1 + J_1 \dd x \dd h 
    \end{align*}
    satisfies the estimate
    \begin{align*}
        u(\epsilon) &\leq \epsilon^{2\gamma-2\theta} \left( 3 [\Lambda]_{X^{\gamma,2}}^2 + \frac{16}{1+\gamma} [\nabla \psi]_{X^{\gamma,2}}^2 \right) .
    \end{align*}
    In particular, by using \eqref{eqn:100}, \eqref{eqn:9}, and \eqref{eqn:10-eps}, we see that the function $u$ satisfies the assumptions of Lemma \ref{lem:bootstrap} with a certain constant $C = C_1(\theta) [\nabla \psi]_{X^{\theta,2}}^2 + C_2(\theta) [\nabla a]_{X^{\theta,2}}^2$.
    Then Lemma \ref{lem:bootstrap} yields
    \begin{align*}
        \sup_{\epsilon \in (0, 1)} \int_{A_\epsilon} w_\theta \int_{\R^d} (L_1 + J_1) \dd x \dd h &\leq \frac{1 + 2 \theta}{2 \theta} C .
    \end{align*}
    We conclude, by taking the supremum over $\epsilon \in (0, 1)$, that $[\Lambda]_{X^{\theta,2}}^2 \leq \frac{1 + 2 \theta}{2 \theta} C$. Now dominated convergence allows us to take the limit from \eqref{eqn:99} to \eqref{eqn:99-limit}.
    The additional statement $P \Lambda \in W^{\theta,2}(\R^d)$ follows from the identity
    \begin{align}  \label{eqn:id12} \tag{\texttt{id12}}
    \Delta_h (P\Lambda)(x) \cdot h
        &= P(x) \Delta_h \Lambda(x) \cdot h
        + \imu P(x-h) \Big(
            R_{-h} \psi(x-h) \conj{P(x-h)}
        \\\nonumber &+ \conj{R_{-h} \psi(x-h)} P(x)
            - (1+Q_h(x)) R_{-h} a(x-h)
        \Big) .
    \end{align}
    
    \textbf{Direction ``$\Longleftarrow$''}.
    Here we assume $\theta < \frac12$.
    If $\Lambda \in W^{\theta,2}(\R^d; \R^d)$ then the assumptions $P \Lambda \in W^{\theta,2}(\R^d; \C^d)$ and $\mathcal{Z}(P, \Lambda) < \infty$, where
    \begin{align*}
        \mathcal{Z}(P, \Lambda) &= \int_{\R^d} |h|^{-d-2-2\theta} \int_{\R^d} |\Delta_h P(x) \Lambda(x) \cdot h|^2 \dd x \dd h ,
    \end{align*}
    are equivalent.
    We start as in the previous step, working with the enlarged annulus $W_\epsilon = \{h \in \R^d: |h| > \epsilon\}$ instead of $A_\epsilon$. This has the advantage that $r W_\epsilon = W_{r \epsilon}$ for $r \in (0, 1)$.
    Then the analogue of \eqref{eqn:99} is
    \begin{align} \label{eqn:99-reverse}
        \int_{W_\epsilon} w_\theta \int_{\R^d} L_1 + J_1 - G_7 \dd x \dd h
        &= \int_0^1 r^{2 \theta} \int_{r W_\epsilon} w_\theta \int_{\R^d} L_1 + J_1 - G_7 
        \\ &\hspace{30mm} + G_6 - (1 + \RE\left[Q_h\right]) L_1 \dd x \dd h \dd r .
    \end{align}
    We define
    \begin{align*}
        D(x, h) &= |\Delta_h \nabla \psi(x) \cdot h|^2
        \\ Z(x, h) &= |\Delta_h P(x) \Lambda(x) \cdot h|^2
    \end{align*}
    and note that 
    \begin{align} \label{eqn:id13} \tag{\texttt{id13}} 
        B_h^- \psi(x) &= - \conj{P(x)}\Delta_h \nabla \psi(x) \cdot h - (1 + \conj{P(x)} P(x-h)) R_{-h} a(x-h) 
        \\\nonumber &+ \imu \conj{P(x)} \Delta_h P(x) \Lambda(x) \cdot h - \imu (1 - \conj{Q_h(x)}) \Delta_h \Lambda(x) \cdot h .
    \end{align}
    Then
    \begin{align*}
        L_1 + J_1 - G_7 &= 2 D + H_1
    \end{align*}
    where
    \begin{align*}
        H_1(x,h) &= 2 \Big|(P(x) + P(x-h)) R_{-h} a(x-h) - \imu \Delta_h P(x) \Lambda(x-h) \cdot h\Big|^2
        \\ &+ 4 \RE\left[(\Delta_h \nabla \psi(x) \cdot h) \big((\conj{P(x)} + \conj{P(x-h)}) R_{-h} a(x-h) + \imu \conj{\Delta_h P(x)} \Lambda(x-h) \cdot h \big) \right]
        \\ &+ 2 L_1(x,h) + 4 (\Delta_h \Lambda(x) \cdot h) \IM\left[B_h^- \psi(x)\right]
        \\ &- 2 B_h^- a(x) \RE\left[B_h^- \psi(x)\right] + 2 B_h^+ a(x) \RE\left[B_h^+ \psi(x)\right] - 2 (B_h^+ a(x))^2 .
    \end{align*}
    On the other hand, we have the identity
    \begin{align} \label{eqn:id14} \tag{\texttt{id14}} 
        \Delta_h \nabla \psi(x) \cdot h + P(x) P(x-h) \conj{\Delta_h \nabla \psi(x) \cdot h}
        &= (P(x) + P(x-h)) \Delta_h \nabla a(x) \cdot h
        \\\nonumber &+ \imu \Delta_h P(x) (\Lambda(x) + \Lambda(x-h)) \cdot h .
    \end{align}
    This implies
    \begin{align*}
        G_6 - (1 + \RE\left[Q_h\right]) L_1 &= 2 D + H_2
    \end{align*}
    where
    \begin{align*}
        H_2(x, h) &= - \frac12 \Big|(P(x) + P(x-h)) \Delta_h \nabla a(x) \cdot h + \imu \Delta_h P(x) (\Lambda(x) + \Lambda(x-h)) \cdot h\Big|^2
        \\ &+ 2 B_h^+ a(x) \RE\left[B_h^+ \psi(x) Q_h(x)\right] + 2 B_h^- a(x) \RE\left[B_h^- \psi(x)\right]
        \\ &- 2 (B_h^+ a(x))^2 - (1 + \RE\left[Q_h(x)\right]) (B_h^- a(x))^2 - (1 + \RE\left[Q_h(x)\right]) L_1(x,h) .
    \end{align*}
    We therefore rewrite \eqref{eqn:99-reverse} as
    \begin{align} \label{eqn:99-reverse-2}
        & \int_0^1 r^{2 \theta} \int_{r W_\epsilon} w_\theta \int_{\R^d} 4 D \dd x \dd h \dd r 
        - \int_{W_\epsilon} w_\theta \int_{\R^d} 2 D \dd x \dd h
        \\\nonumber &= \int_{W_\epsilon} w_\theta \int_{\R^d} H_1 \dd x \dd h
        - \int_0^1 r^{2 \theta} \int_{W_{r \epsilon}} w_\theta \int_{\R^d} H_1 + H_2 \dd x \dd h \dd r .
    \end{align}
    Note that for any $\delta > 0$ there exists a constant $C_\delta > 0$ such that
    \begin{align*}
        |H_1(x, h)| + |H_2(x, h)| &\leq \frac{\delta}{3} (D + |R_h \psi(x)|^2 + |R_{-h} \psi(x-h)|^2)
        \\ &+ C_\delta (|\Delta_h \nabla a(x) \cdot h|^2 + |R_h a(x)|^2 + |R_{-h} a(x-h)|^2 
        \\ &+ L_1(x, h) + Z(x, h) + Z(x-h, -h)) .
    \end{align*}
    By the same argument as in the proof of \eqref{eqn:9}, we can show that
    \begin{align} \label{eqn:9-extra}
        \int_{W_\epsilon} w_\theta \int_{\R^d} |R_h \psi(x)|^2 \dd x \dd h
        &\leq \int_0^1 r^{2 \theta} \int_{W_{r \epsilon}} w_\theta \int_{\R^d} |\Delta_h \nabla \psi(x) \cdot h|^2 \dd x \dd h \dd r .
    \end{align}
    We also have the inequality
    \begin{align} \label{eqn:9-extra-2}
        \int_{W_\epsilon} w_\theta \int_{\R^d} D \dd x \dd h
        &\leq (1+2\theta) \int_0^1 r^{2 \theta} \int_{W_{r \epsilon}} w_\theta \int_{\R^d} D \dd x \dd h \dd r ,
    \end{align}
    which uses the fact that the integral on the left-hand side increases when $\epsilon$ decreases.
    We need to justify the existence of all integrals that appear above and in the intermediate calculations once more. It is a consequence of $\nabla \psi, \nabla a, \Lambda \in L^2(\R^d)$ and the inequality
    \begin{align} \label{eqn:300}
        \|(\Delta_h \psi, \Delta_h a)\|_{L^2} + \|(R_h \psi, R_h a)\|_{L^2} 
        + \|(\Delta_h \nabla \psi \cdot h, \Delta_h \nabla a \cdot h, \Lambda \cdot h)\|_{L^2} \lesssim |h| \|\nabla \psi\|_{L^2} .
    \end{align}
    Define
    \begin{align*}
        v(\epsilon) &= \int_0^1 r^{2 \theta} \int_{W_{r \epsilon}} w_\theta \int_{\R^d} D \dd x \dd h \dd r .
    \end{align*}
    Applying \eqref{eqn:9-extra} and \eqref{eqn:9-extra-2} to \eqref{eqn:99-reverse-2} and using \eqref{eqn:9} yields
    \begin{align*}
        & (2 - 4 \theta - \delta (2 + 2 \theta)) v(\epsilon)
        - \delta \int_0^1 r^{2 \theta} v(r \epsilon) \dd r
        \leq M
        \\ \text{ where } M &= C_\delta \frac{2+2\theta}{1+2\theta} \left( \frac{2+\theta}{1+\theta} [\nabla a]_{X^{\theta,2}}^2 + [\Lambda]_{X^{\theta,2}}^2 + 2 \mathcal{Z}(P, \Lambda) \right) .
    \end{align*}
    We choose $\delta = \frac{1 - 2 \theta}{2 + \theta}$ and note that the above implies
    \begin{align*}
        v(\epsilon) - \frac23 \int_0^1 r^{2 \theta} v(r \epsilon) \dd r
        \leq \frac{M}{2 \delta} .
    \end{align*}
    Now we define $u(\epsilon) = v(\epsilon^{\frac32})$ and $\alpha = \frac14 + \frac32 \theta$ so that
    \begin{align*}
        u(\epsilon) - \int_0^1 r^{2 \alpha} u(r \epsilon) \dd r
        \leq \frac{M}{2 \delta} .
    \end{align*}
    Note that \eqref{eqn:300} implies $v(\epsilon) \lesssim \epsilon^{-2\theta} \|\nabla \psi\|_{L^2}^2$ and hence $u(\epsilon) \lesssim \epsilon^{\frac12-2\alpha}$.
    We can now apply Lemma \ref{lem:bootstrap} to obtain a uniform bound for $v(\epsilon)$. 
    Finally another application of  \eqref{eqn:9-extra-2} yields $\nabla \psi \in W^{\theta,2}(\R^d)$.
\end{proof}

\begin{corollary} \label{cor:1}
    Let $d \in \N^+$ and $\theta \in (0, \frac12)$. Let $\psi \in W^{1+\theta,2}_{\loc}(\R^d; \C)$ with $\nabla \psi \in \dot{W}^{\theta,2}(\R^d; \C^d)$ and let $P \in \mathcal{P}(\psi)$ be a polar factor.
    We have $\nabla a, \Lambda, P \Lambda \in \dot{W}^{\theta,2}(\R^d)$ and
    \begin{align}
        \label{eqn:0} (1-2\theta)^2 [\nabla a]_{X^{\theta,2}}^2 + (1+2\theta) (1-2\theta) [\Lambda]_{X^{\theta,2}}^2 &\leq 8 [\nabla \psi]_{X^{\theta,2}}^2 .
    \end{align}
    There exists an unbounded sequence $(\psi_n)_{n \in \N} \subset W^{1+\theta,2}(\R^d)$ for which $[\nabla a_n]_{X^{\theta,2}}^2 + [\Lambda_n]_{X^{\theta,2}}^2 \lesssim 1$. Furthermore, for this sequence $[P_n \nabla a_n]_{X^{\theta,2}} \lesssim 1$ while $[P_n \Lambda_n]_{X^{\theta,2}}$ is unbounded.
    Nevertheless, if $\nabla \psi \in W^{\theta,2}(\R^d)$ then we have the equivalence
    \begin{align}
        \label{eqn:0-equiv} 
        [\nabla \psi]_{X^{\theta,2}}^2 &\leq \frac{32}{(1-2\theta)^2} [\nabla a]_{X^{\theta,2}}^2 + \frac{32}{(1-2\theta)^2} [\Lambda]_{X^{\theta,2}}^2 + 12 [P \Lambda]_{X^{\theta,2}}^2
        \leq \frac{796}{(1-2\theta)^4} [\nabla \psi]_{X^{\theta,2}}^2 .
    \end{align}
\end{corollary}
\begin{proof}
    By Lemma \ref{lem:approx-2}, applied with $\alpha = \nu = \theta$, there exist $\psi_n \in W^{1+\theta,2}_{\loc}(\R^d; \C)$ with $\nabla \psi_n \in W^{\theta,2}(\R^d; \C^d)$ such that $[\nabla \psi_n - \nabla \psi]_{X^{\theta,2}} \rightarrow 0$ and $\psi_n \rightarrow \psi$ in $W^{1,2}_{\loc}(\R^d)$. 
    Since $\psi_n \rightarrow \psi$ in $W^{1,2}_{\loc}(\R^d)$, Proposition \ref{prop:2} implies $\nabla_{\!\circ} \psi_n \rightarrow \nabla_{\!\circ} \psi$ in $L^2_{\loc}(\R^d)$.
    On the other hand, $\psi_n$ is admissible for Proposition \ref{prop:1} and Theorem \ref{thm:1}, which implies $\nabla_{\!\circ} \psi_n \in \dot{W}^{\theta,2}(\R^d)$. We similarly obtain $P \Lambda \in \dot{W}^{\theta,2}(\R^d)$ from \eqref{eqn:PLambda-est} below.
    Applying \eqref{eqn:9}--\eqref{eqn:9-alt}, \eqref{eqn:id10}, and subsequently Young's inequality for products with parameters $c_1, c_2 > 0$ to the explicit formula \eqref{eqn:formula-pre}, which is a better formulation of \eqref{eqn:formula} for this estimate, yields
    \begin{align*}
        (1-2\theta) [\nabla a]_{X^{\theta,2}}^2 + (1+2\theta) [\Lambda]_{X^{\theta,2}}^2
        &\leq (1-2\theta) \frac{4}{2+2\theta} [\nabla \psi]_{X^{\theta,2}}^2
        \\ &+ 4(1+2\theta) \left( \frac{c_1}{4+4\theta} [\nabla \psi]_{X^{\theta,2}}^2 + \frac{c_1^{-1}}{4+4\theta} [\nabla a]_{X^{\theta,2}}^2 \right)
        \\ &+ 8\theta \left( \frac{c_2}{4+4\theta} [\nabla \psi]_{X^{\theta,2}}^2 + \frac{c_2^{-1}}{2} [\Lambda]_{X^{\theta,2}}^2 \right) .
    \end{align*}
    We choose
    \begin{align*}
        c_1 &= \frac{2(1+2\theta)}{(1+\theta)(1-2\theta)} & c_2 &= \frac{8\theta}{1+2\theta}
    \end{align*}
    and obtain \eqref{eqn:0}.
    Finally, Lemma \ref{lem:approx-1} extends the inequality to $\psi$.

    For the counterexample to reverse estimates analogous to \eqref{eqn:0}, choose $a_n = n^{-1} g$ for some $0 \neq g \in C_c^\infty(\R^d; [0, 1])$, $P_n = \eul^{\imu n x_1}$ and set $\psi_n = a_n P_n$. Clearly $[\nabla a_n]_{X^{\theta,2}}^2 + [\Lambda_n]_{X^{\theta,2}}^2 \lesssim 1$. Using the Fourier characterization of the seminorm, one calculates that $\|\nabla \psi_n\|_{\dot{H}^\theta}^2$ and $\|P_n \Lambda_n\|_{\dot{H}^\theta}^2$ are unbounded and $\|P_n \nabla a_n\|_{\dot{H}^\theta}^2$ is bounded. 
    
    The reverse inequality in \eqref{eqn:0-equiv} follows similarly from \eqref{eqn:formula}, \eqref{eqn:9}, and applications of Young's inequality for products with parameters $c_3, \dots, c_6 > 0$. Specifically, we have
    \begin{align*}
        (1-2\theta) [\nabla \psi]_{X^{\theta,2}}^2
        \!&\leq (1-2\theta) [\nabla a]_{X^{\theta,2}}^2 + (1+2\theta) [\Lambda]_{X^{\theta,2}}^2
        \\ &+ (2+4\theta) \left( \left( \frac{c_3}{4+4\theta} + \frac{c_4}{2} \right) [\nabla \psi]_{X^{\theta,2}}^2 + \left( \frac{c_3^{-1}}{2} + \frac{c_4^{-1}}{4+4\theta} \right)[\nabla a]_{X^{\theta,2}}^2 \right)
        \\ &+ (2+4\theta) \left( \frac{c_5}{4+4\theta} [\nabla \psi]_{X^{\theta,2}}^2 + \frac{c_5^{-1}}{2} [\Lambda]_{X^{\theta,2}}^2 \right)
        \\ &+ (2-4\theta) \left( \frac{c_6}{4+4\theta} [\nabla \psi]_{X^{\theta,2}}^2 + \frac{c_6^{-1}}{2} [P\Lambda]_{X^{\theta,2}}^2 \right) .
    \end{align*}
    The choices
    \begin{align*}
        c_3 &= \frac{(1-2\theta)(1+\theta)}{3(1+2\theta)}
        & c_4 &= \frac{1-2\theta}{6(1+2\theta)}
        & c_5 &= \frac{(1-2\theta)(1+\theta)}{6(1+2\theta)}
        & c_6 &= \frac{1+\theta}{6}
    \end{align*}
    then yield the first inequality in \eqref{eqn:0-equiv}.
    Note that \eqref{eqn:id12} and \eqref{eqn:9} imply
    \begin{align} \label{eqn:PLambda-est}
        [P \Lambda]_{X^{\theta,2}}^2 &\leq 5 \left( [\nabla\psi]_{X^{\theta,2}}^2 + [\nabla a]_{X^{\theta,2}}^2 + [\Lambda]_{X^{\theta,2}}^2 \right).
    \end{align}
    Combining this with \eqref{eqn:0} yields the second inequality in \eqref{eqn:0-equiv}.
\end{proof}

\begin{corollary} \label{cor:1-extra}
    Let $d \in \N^+$ and $\theta \in [\frac12, 1)$. Let $\psi \in W^{1+\theta,2}_{\loc}(\R^d; \C)$ with $\nabla \psi \in W^{\theta,2}(\R^d; \C^d)$ and let $P \in \mathcal{P}(\psi)$ be a polar factor. Assume that $a \in W^{1+\theta,2}_{\loc}(\R^d; \R)$ with $\nabla a \in W^{\theta,2}(\R^d; \R^d)$.
    We have $\Lambda \in W^{\theta,2}(\R^d)$ and
    \begin{align}
        \label{eqn:0-extra} [\Lambda]_{X^{\theta,2}}^2 &\leq (4 \theta - 2) [\nabla a]_{X^{\theta,2}}^2 + 4 [\nabla a]_{X^{\theta,2}} [\nabla \psi]_{X^{\theta,2}} + 2 [\nabla \psi]_{X^{\theta,2}}^2 .
    \end{align}
\end{corollary}
\begin{proof}
    We again use Theorem \ref{thm:1} to obtain $\Lambda, P\Lambda \in W^{\theta,2}(\R^d)$.
    We may discard the first term in the right-hand side of the explicit formula \eqref{eqn:formula-pre} since $1-2\theta \leq 0$.
    Applying \eqref{eqn:9} and Young's inequality for products with parameter $c_2 > 0$ to \eqref{eqn:formula-pre} then yields
    \begin{align*}
        (1-2\theta) [\nabla a]_{X^{\theta,2}}^2 + (1+2\theta) [\Lambda]_{X^{\theta,2}}^2
        &\leq \frac{2(1+2\theta)}{1+\theta} [\nabla a]_{X^{\theta,2}} [\nabla \psi]_{X^{\theta,2}}
        \\ &+ 8\theta \left( \frac{c_2}{4+4\theta} [\nabla \psi]_{X^{\theta,2}}^2 + \frac{c_2^{-1}}{2} [\Lambda]_{X^{\theta,2}}^2 \right) .
    \end{align*}
    We choose $c_2 = \frac{8\theta}{1+2\theta}$ to obtain
    \begin{align*}
        [\Lambda]_{X^{\theta,2}}^2 
        &\leq \frac{2(2\theta-1)}{1+2\theta} [\nabla a]_{X^{\theta,2}}^2
        + \frac{4}{1+\theta} [\nabla a]_{X^{\theta,2}} [\nabla \psi]_{X^{\theta,2}}
        + \frac{32\theta^2}{(1+\theta)(1+2\theta)^2} [\nabla \psi]_{X^{\theta,2}}^2.
    \end{align*}
    This implies \eqref{eqn:0-extra}.
\end{proof}

\subsubsection{Regularity at the endpoint $s = \frac32$}

Recall the definition of the Besov spaces $B^{\frac32}_{2,1}(\R^d)$ and $B^{\frac32}_{2,\infty}(\R^d)$ from \cite[Section 2.5.12]{Triebel1983}, which use the norms
\begin{align}
    \label{eqn:def-besov-1} \|f\|_{B^{\frac32}_{2,1}}^2 &= \|f\|_{W^{1,2}}^2 + [\nabla f]_{B^{\frac12}_{2,1}}^2 \text{ where }
    [\nabla f]_{B^{\frac12}_{2,1}}^2 = \left\| \mathds{1}_{\{|h| \leq 1\}} |h|^{-d-\frac12} \|\Delta_h \nabla f(x)\|_{L^2_x} \right\|_{L^1_h}^2
    \\ \label{eqn:def-besov-2} \|f\|_{B^{\frac32}_{2,\infty}}^2 &= \|f\|_{W^{1,2}}^2 + [\nabla f]_{B^{\frac12}_{2,\infty}}^2 \text{ where }
    [\nabla f]_{B^{\frac12}_{2,\infty}}^2 = \left\|\mathds{1}_{\{|h| \leq 1\}} |h|^{-\frac12} \|\Delta_h \nabla f(x)\|_{L^2_x} \right\|_{L^\infty_h}^2 .
\end{align}
Using Jensen's inequality and $R_h f(x) = - \int_0^1 \Delta_{rh} \nabla f(x) \cdot h \dd r$, one finds that
\begin{align*}
    \int_{\R^d} |h|^{-d-\frac32} \left( \int_{\R^d} |R_h f(x)|^2 \dd x \right)^{\frac12} \dd h &\leq C(d) \|f\|_{B^{\frac32}_{2,1}}
    \\ \left\||h|^{-\frac32} \left(\int_{\R^d} |R_h f(x)|^2 \dd x\right)^{\frac12} \right\|_{L^\infty_h} &\leq C(d) \|f\|_{B^{\frac32}_{2,\infty}} .
\end{align*}

\begin{proposition}[{Consequence of \cite[Theorem 1.1]{Moussai2018}}] \label{prop:besov}
    Let $d \in \N^+$. There exists a constant $C > 0$ such that for every $\psi \in B^{\frac32}_{2,1}(\R^d; \C)$ we have $a \in B^{\frac32}_{2,\infty}(\R^d; \R)$ and
    \begin{align*}
        \|a\|_{B^{\frac32}_{2,\infty}} &\leq C \|\psi\|_{B^{\frac32}_{2,1}} .
    \end{align*}
\end{proposition}
\begin{proof}
    By \eqref{eqn:a-reduction} it suffices to assume that $\psi$ is real-valued. The proof is finished by checking that $f(x) = |x|$ is an admissible function for \cite[Theorem 1.1]{Moussai2018}.
\end{proof}

\begin{corollary} \label{cor:2}
    Let $d \in \N^+$. There exists a constant $C > 0$ such that for every $\psi \in B^{\frac32}_{2,1}(\R^d; \C)$ we have $\Lambda \in W^{\frac12,2}(\R^d; \R^d)$ and 
    \begin{align*}
        [\Lambda]_{X^{\frac12,2}}^2 &\leq C \|\psi\|_{B^{\frac32}_{2,1}}^2 .
    \end{align*}
\end{corollary}
\begin{proof}
    Note that for every $\theta \in (0, \frac12)$ we are in the setting of Theorem \ref{thm:1} and the identity \eqref{eqn:formula} is available. Let us assume $\theta > \frac14$ and apply Young's inequality for products, \eqref{eqn:id10}, and \eqref{eqn:9}--\eqref{eqn:9-alt}. This yields
    \begin{align*}
        [\Lambda]_{X^{\theta,2}}^2 &\leq C [\nabla \psi]_{W^{\theta,2}}^2 + C \int_{\R^d} |h|^{-d-2-2\theta} \int_{\R^d} |R_h \psi(x)| |R_{-h} a(x-h)| \dd x \dd h
    \end{align*}
    with a constant $C = C(d) > 0$. For the region $|h| \leq 1$ we use $|h|^{-d-2-2\theta} \leq |h|^{-d-3}$ and Hölder's inequality. For $|h| > 1$ we estimate with $\|\psi\|_{W^{1,2}(\R^d)}^2$ and the embedding $B^{\frac32}_{2,1}(\R^d) \hookrightarrow W^{1,2}(\R^d)$. Then Fatou's lemma yields
    \begin{align*}
        [\Lambda]_{X^{\frac12,2}}^2 \leq \sup_{\theta \in (\frac14, \frac12)} [\Lambda]_{X^{\theta,2}}^2 &\leq C \|\psi\|_{B^{\frac32}_{2,1}}^2 + C \|a\|_{B^{\frac32}_{2,\infty}}^2
    \end{align*}
    and Proposition \ref{prop:besov} finishes the proof.
\end{proof}

\begin{lemma} \label{lem:endpoint-epsilon}
    Let $\epsilon \in (0, \frac12)$. If $\psi \in W^{\frac32+\epsilon,2}_{\loc}(\R^d; \C)$ with $\nabla \psi \in W^{\frac12+\epsilon,2}(\R^d; \C^d)$, then $\Lambda \in W^{\frac12,2}(\R^d; \R^d)$ and
    \begin{align*}
        [\Lambda]_{X^{\frac12,2}}^2 &\leq \frac43 [\nabla \psi]_{X^{\frac12,2}}^2 + \frac{7}{\epsilon} [\nabla \psi]_{X^{\frac12+\epsilon,2}} [\nabla \psi]_{X^{\frac12-\epsilon,2}} .
    \end{align*}
\end{lemma}
\begin{proof}
    This is a consequence of taking the limit $\theta \nearrow \frac12$ (involving the use of Fatou's lemma) in \eqref{eqn:formula} and using \eqref{eqn:9}, analogous to the proof of Corollary \ref{cor:2}. 
    By assuming $\frac12+\epsilon$ regularity for $\nabla \psi$, we only need to use $\frac12-\epsilon < \frac12$ regularity of $\nabla a$. Applying Young's inequality for products with parameters $c_1, c_2, c_3 > 0$ yields
    \begin{align*}
        2[\Lambda]_{X^{\frac12,2}}^2 &\leq \left( \frac{4 c_1}{2 (3+2\epsilon)} + \frac{4 c_2}{2} \right) [\nabla\psi]_{X^{\frac12+\epsilon,2}}^2 + \left( \frac{4 c_1^{-1}}{2} + \frac{4 c_2^{-1}}{2 (3-2\epsilon)} \right) [\nabla a]_{X^{\frac12-\epsilon,2}}^2 
        \\ &+ \frac{4 c_3}{6}[\nabla\psi]_{X^{\frac12,2}}^2 + \frac{4 c_3^{-1}}{2} [\Lambda]_{X^{\frac12,2}}^2 .
    \end{align*}
    With the choices
    \begin{align*}
        c_1 &= \sqrt{3+2\epsilon} \frac{[\nabla a]_{X^{\frac12-\epsilon,2}}}{[\nabla\psi]_{X^{\frac12+\epsilon,2}}}
        & c_2 &= \frac{1}{\sqrt{3-2\epsilon}} \frac{[\nabla a]_{X^{\frac12-\epsilon,2}}}{[\nabla\psi]_{X^{\frac12+\epsilon,2}}}
        & c_3&=2
    \end{align*}
    we obtain
    \begin{align*}
        [\Lambda]_{X^{\frac12,2}}^2 &\leq
        \frac43 [\nabla\psi]_{X^{\frac12,2}}^2
        + 4 \left( \frac{1}{\sqrt{3+2\epsilon}} + \frac{1}{\sqrt{3-2\epsilon}} \right) [\nabla\psi]_{X^{\frac12+\epsilon,2}} [\nabla a]_{X^{\frac12-\epsilon,2}} .
    \end{align*}
    The proof is finished by applying \eqref{eqn:0}. Note that we may assume $[\nabla \psi]_{X^{\frac12+\epsilon,2}} > 0$ because the case of constant $\psi$ is trivial.
\end{proof}

\subsubsection{Regularity at the endpoint $s = 2$}

\begin{corollary} \label{cor:2-extra}
    Let $d \in \N^+$. For every $\psi \in W^{2,2}_{\loc}(\R^d; \C)$ with $\nabla \psi \in W^{1,2}(\R^d; \C^d)$ and in addition $a \in W^{2,2}_{\loc}(\R^d; \R)$ and $\nabla a \in W^{1,2}(\R^d; \R^d)$, we have $\Lambda \in W^{1,2}(\R^d; \R^d)$ and
    \begin{align} 
        \|\nabla \Lambda\|_{L^2}^2 &\leq 6 \|\nabla^2 a\|_{L^2}^2 + 12 \|\nabla^2 a\|_{L^2} \|\nabla^2 \psi\|_{L^2} + 6 \|\nabla^2 \psi\|_{L^2}^2 .
    \end{align}
\end{corollary}
\begin{proof}
    By interpolation we can apply Corollary \ref{cor:1-extra} for every $\theta \in [\frac12, 1)$. We multiply \eqref{eqn:0-extra} with $(1 - \theta)$ on both sides and take the limit $\theta \nearrow 1$. By \eqref{eqn:X-H} we have
    \begin{align*}
        (1 - \theta) [\Lambda]_{W^{\theta,2}}^2 &\leq (1 - \theta) (d + 2 \theta) [\Lambda]_{X^{\theta,2}}^2 
        \\ &\leq (1 - \theta) (1 + 2 \theta) \left( (4 \theta - 2) [\nabla a]_{W^{\theta,2}}^2 + 4 [\nabla a]_{W^{\theta,2}} [\nabla \psi]_{W^{\theta,2}} + 2 [\nabla \psi]_{W^{\theta,2}}^2 \right) .
    \end{align*}
    The Bourgain--Brezis--Mironescu formula \cite[Theorem 2]{BourgainBrezisMironescu2001} (see also \cite[Remark 4.3]{DiNezzaPalatucciValdinoci2012}) allows us to take the limit $\theta \nearrow 1$ on both sides, which finishes the proof.
\end{proof}

\subsection{Regularity of $J$}

\begin{proof}[Proof of Proposition \ref{prop:J}]
    Set $\theta = s - 1$. We distinguish the case $d = 1$ because we do not want to consider the homogeneous fractional Sobolev space $\dot{W}^{\theta,2}(\R^d; \C^d)$ when $\theta \geq \frac{d}{2}$, as then the characterization \eqref{eqn:characterization} is not available. The proof is given below only for the case where the homogeneous space is used.
    We use the identity
    \begin{align} 
        \Delta_h J_j(x) \cdot h
        &= - \IM\left[\conj{\psi_j(x)} R_h \psi_j(x) + \conj{\psi_j(x-h)} R_{-h} \psi_j(x-h)\right] .
	\end{align}
    After squaring both sides and applying \eqref{eqn:9}, we obtain boundedness. Since $|J_j| = a_j |\Lambda_j| \leq a_j |\nabla \psi_j|$, we have also $J_j \in L^{\frac{2d}{d-2\theta}}(\R^d)$, and hence $J_j \in \dot{W}^{\theta,2}(\R^d)$.
    Since $R_h$ is a linear operator, we furthermore have 
    \begin{align} 
        \Delta_h (J_1 - J_2)(x) \cdot h
        &= - \IM\left[\conj{\psi_1(x)} R_h (\psi_1 - \psi_2)(x) + \conj{\psi_1(x-h)} R_{-h} (\psi_1 - \psi_2)(x-h)\right] 
        \\ &- \IM\left[\conj{(\psi_1 - \psi_2)(x)} R_h \psi_2(x) + \conj{(\psi_1 - \psi_2)(x-h)} R_{-h} \psi_2(x-h)\right] ,
	\end{align}
    which yields local Lipschitz continuity.
\end{proof}

\section{Continuity}

\begin{theorem} \label{thm:continuity}
    Let $d \in \N^+$ and $\theta \in (0, 1)$. Suppose that
    \begin{align*}
        \psi_n &\rightarrow \psi \text{ in } W^{1+\theta,2}_{\loc}(\R^d; \C)
        \\ \nabla \psi_n &\rightarrow \nabla \psi \text{ in } W^{\theta,2}(\R^d; \C^d) .
    \end{align*}
    Assume in addition that
    \begin{align*}
        \theta < \frac12
        \quad\;\text{ or }\quad\;
        \begin{aligned}
            &\;\;\theta = \frac12
            \\[-0.8mm] &a_n \in W^{1+\theta,2}_{\loc}(\R^d)
            \\ &\;\sup_{n \in \N} \,[\nabla a_n]_{X^{\theta,2}} < \infty
        \end{aligned}
        \quad\;\text{ or }\quad\;
        \begin{aligned}
            \theta &> \frac12
            \\[-0.8mm] \nabla a_n &\rightarrow \nabla a \text{ in } W^{\theta,2}(\R^d)
        \end{aligned} .
    \end{align*}
    Then
    \begin{align*}
        \Lambda_n \rightarrow \Lambda \text{ in } W^{\theta,2}(\R^d) 
        \quad\;\text{ and if $\theta \neq \frac12$ also }\quad\; \begin{aligned}
            a_n &\rightarrow a \text{ in } W^{1+\theta,2}_{\loc}(\R^d)
            \\ \nabla a_n &\rightarrow \nabla a \text{ in } W^{\theta,2}(\R^d)
        \end{aligned} .
    \end{align*}
\end{theorem}
\begin{proof}
    Proposition \ref{prop:2} implies $(\nabla a_n, \Lambda_n) \rightarrow (\nabla a, \Lambda)$ in $(L^2 \times L^2)(\R^d)$.
    By Proposition \ref{prop:1} and Corollaries \ref{cor:1} and \ref{cor:1-extra} the sequences $(\nabla a_n)_{n \in \N}, (\Lambda_n)_{n \in \N} \subset W^{\theta,2}(\R^d)$ are bounded, and hence the Banach--Alaoglu theorem implies, after passing to subsequences, the weak convergence $(\nabla a_n, \Lambda_n) \rightharpoonup (\nabla a, \Lambda)$ in $(W^{\theta,2} \times W^{\theta,2})(\R^d)$.
    Furthermore, in any case we have $a \in W^{1+\theta,2}_{\loc}(\R^d)$ and $\nabla a, \Lambda \in W^{\theta,2}(\R^d)$.
    We note also that $\nabla a_n \rightarrow \nabla a$ in $W^{\theta,2}(\R^d)$ together with $\psi_n \rightarrow \psi$ in $W^{1+\theta,2}_{\loc}(\R^d)$ implies $a_n \rightarrow a$ in $W^{1+\theta,2}_{\loc}(\R^d)$.
    It therefore suffices to prove that for any subsequence $(n_k)_{k \in \N}$ there exists a subsequence $(n_{k_j})_{j \in \N}$ so that $\Lambda_{n_{k_j}} \rightarrow  \Lambda$ in $W^{\theta,2}(\R^d)$ and, if $\theta < \frac12$, also $\nabla a_{n_{k_j}} \rightarrow \nabla a$ in $W^{\theta,2}(\R^d)$. 
    This justifies the passage to a subsequence we have performed above. Note that we suppress these subsequences from our notation.

    We assume for now that $\theta < \frac12$ and equip the space $(W^{\theta,2} \times W^{\theta,2})(\R^d)$ with the equivalent norm
    \begin{align*}
        \|(\nabla a, \Lambda)\|_\theta^2 &= \|(\nabla a, \Lambda)\|_{L^2}^2 + [(\nabla a, \Lambda)]_\theta^2
        \\ \text{ where } [(\nabla a, \Lambda)]_\theta^2 &= (1-2\theta) [\nabla a]_{X^{\theta,2}}^2 + (1+2\theta) [\Lambda]_{X^{\theta,2}}^2 .
    \end{align*}
    As weak convergence in $(W^{\theta,2} \times W^{\theta,2})(\R^d)$ and strong convergence in $(L^2 \times L^2)(\R^d)$ have already been established, it suffices to show convergence of the squared seminorm
    \begin{align} \label{eqn:norm-conv}
        [(\nabla a_n, \Lambda_n)]_\theta^2 \rightarrow [(\nabla a, \Lambda)]_\theta^2 .
    \end{align}
    However, if $\theta = \frac12$ then $\nabla a_n$ disappears from the seminorm and \eqref{eqn:norm-conv} still suffices.
    Lastly, if $\theta > \frac12$ then assuming strong convergence of $\nabla a_n$ in $W^{\theta,2}(\R^d)$ again implies that \eqref{eqn:norm-conv} suffices. 
    It therefore remains to prove \eqref{eqn:norm-conv}.
    We do so by taking the limit $n \rightarrow \infty$ on both sides of the crucial identity \eqref{eqn:formula}.
    Specifically, it suffices to prove
    \begin{align} \label{eqn:limit}
        \int_{\R^d} w_\theta \int_{\R^d} I_n - I \dd x \dd h \rightarrow 0
    \end{align}
    where $w_\theta(h) I_n(x, h)$ and $w_\theta(h) I(x, h)$ are the integrands on the right-hand sides of \eqref{eqn:formula} for $\psi_n$ and $\psi$.
    
    We may assume that $(\psi_n, \nabla \psi_n, a_n, \nabla a_n, \Lambda_n) \rightarrow (\psi, \nabla \psi, a, \nabla a, \Lambda)$ pointwise almost everywhere.
    To do so, we first consider an arbitrary subsequence and use the convergence of $\psi_n$ and $\nabla \psi_n$ in the respective Sobolev spaces to switch to another subsequence where these quantities converge pointwise almost everywhere.
    By elementary inequalities the same pointwise convergence then holds, on this subsequence, for $a_n, \nabla a_n$, and $\Lambda_n$.
    Proving \eqref{eqn:limit} for a specific subsequence of an arbitrary subsequence implies the statement for the original sequence.
    
    We first prove that
    \begin{align} \label{eqn:small-h}
        \lim_{\delta \searrow 0} \sup_{n \in \N} \int_{\{|h| < \delta\}} w_\theta \int_{\R^d} |I_n| + |I| \dd x \dd h &= 0 .
    \end{align}
    Using the boundedness properties of Corollaries \ref{cor:1} and \ref{cor:1-extra} and applying Hölder's inequality, \eqref{eqn:id10}, and \eqref{eqn:9} to $I_n$ and $I$ yields 
    \begin{align*}
        & \int_{\{|h| < \delta\}} w_\theta \int_{\R^d} |I_n| + |I| \dd x \dd h 
        \\ &\leq C(\theta) \sup_{n \in \N} \big([\nabla a_n]_{X^{\theta,2}} + [\nabla a]_{X^{\theta,2}} + [\nabla \psi_n]_{X^{\theta,2}} + [\nabla \psi]_{X^{\theta,2}}\big) 
        \\ &\left( \int_{\{|h| < \delta\}} \!w_\theta\! \int_{\R^d} |R_h \psi_n(x)|^2 + |R_h \psi(x)|^2 + |\Delta_h \nabla \psi_n(x) \cdot h|^2 + |\Delta_h \nabla \psi(x) \cdot h|^2 \dd x \dd h \right)^{\frac12} .
    \end{align*}
    The weighted integrands on the right-hand side above converge strongly in $L^1_h L^1_x(\R^{d+d})$ and hence, after integration in $x$, form a uniformly integrable family of functions in $h$. This implies \eqref{eqn:small-h}.
    It remains to prove
    \begin{align} \label{eqn:large-h}
        \forall\, \delta > 0 \quad \lim_{n \rightarrow \infty} \int_{\{|h| \geq \delta\}} w_\theta \int_{\R^d} I_n - I \dd x \dd h &= 0 .
    \end{align}
    Here the key observation is that
    \begin{align} \label{eqn:lim-I}
        I_n(x, h) \rightarrow I(x, h)
    \end{align}
    for almost all $(x, h) \in \R^{d+d}$, even though polar factors $P_n(x)$ and $P_n(x-h)$ do not need to converge pointwise almost everywhere.
    For almost all $x \in \R^d$ with $\psi(x) \neq 0$ we know that there exists some $N \in \N$ such that $\psi_n(x) \neq 0$ holds for $n \geq N$, and in particular $P_n(x) = \frac{\psi_n(x)}{|\psi_n(x)|} \rightarrow \frac{\psi(x)}{|\psi(x)|} = P(x)$.
    If $\psi(x-h) \neq 0$ also holds, then $P_n(x-h) \rightarrow P(x-h)$ and we can take the limit \eqref{eqn:lim-I}. If instead $\psi(x-h) = 0$, then $\Lambda_n(x-h) \rightarrow \Lambda(x-h) = 0$ and the term containing the polar factor $P_n(x-h)$ vanishes in the limit, so we also obtain \eqref{eqn:lim-I}.
    If on the other hand $\psi(x) = 0$, and in particular $a(x) = 0$ and $\nabla \psi(x) = \nabla a(x) = \Lambda(x) = 0$, then direct substitution yields $I(x, h) = 0$.
    Furthermore, we shall prove that $I_n(x, h) \rightarrow 0$. 
    Let $c_n \approx_n d_n$ denote $\lim_{n \rightarrow \infty} (c_n - d_n) = 0$. Then
    \begin{align*}
        & (2 + 4 \theta) \RE\left[ \conj{P_n(x)} \left( R_h \psi_n(x) \Delta_h \nabla a_n(x) \cdot h - R_h a_n(x) \Delta_h \nabla \psi_n(x) \cdot h \right) \right]
        \\ &- \IM\left[ \conj{P_n(x)} R_h \psi_n(x) \left( 8 \theta \Delta_h \Lambda_n(x) \cdot h - (2 - 4 \theta) \left(1 - P_n(x) \conj{P_n(x-h)}\right) \Lambda_n(x-h) \cdot h \right) \right] 
        \\ &\approx_n (2 + 4 \theta) \RE\left[ \conj{P_n(x)} \left( R_h \psi_n(x) \Delta_h \nabla a(x) \cdot h - R_h a(x) \Delta_h \nabla \psi_n(x) \cdot h \right) \right]
        \\ &- \IM\left[ \conj{P_n(x)} R_h \psi_n(x) \left( 8 \theta \Delta_h \Lambda(x) \cdot h - (2 - 4 \theta) \left(1 - P_n(x) \conj{P_n(x-h)}\right) \Lambda(x-h) \cdot h \right) \right] 
        \\ &\approx_n (2 + 4 \theta) \RE\left[ \conj{P_n(x)} a_n(x-h) \left( P_n(x-h) \nabla a(x-h) \cdot h - \nabla \psi_n(x-h) \cdot h \right) \right]
        \\ &- \IM\left[ \conj{P_n(x)} P_n(x-h) a_n(x-h) \left( 8 \theta \Lambda(x-h) \cdot h + (2 - 4 \theta) \Lambda(x-h) \cdot h \right) \right] 
        \\ &\approx_n (2 + 4 \theta) \RE\left[ \conj{P_n(x)} P_n(x-h) a_n(x-h) (- \imu \Lambda(x-h) \cdot h) \right]
        \\ &- (2 + 4 \theta) \IM\left[ \conj{P_n(x)} P_n(x-h) a_n(x-h) \Lambda(x-h) \cdot h \right]
        \\ &= 0 .
    \end{align*}
    Generalized dominated convergence now implies \eqref{eqn:large-h}, provided we can estimate $|w_\theta I_n| \leq w_\theta J_n$ for some strongly convergent majorant sequence $(w_\theta J_n)_{n \in \N} \subset L^1(\{|h| \geq \delta\} \times \R^d)$.
    We note first that
    \begin{align*}
        \int_{\{|h| \geq \delta\}} w_\theta(h) |h|^2 \dd h < C(\delta) .
    \end{align*}
    The estimates
    \begin{align*}
        \|R_h f\|_{L^2} &\leq 2 |h| \|\nabla f\|_{L^2} & \|\Delta_h F \cdot h\|_{L^2} &\leq 2 |h| \|F\|_{L^2}
    \end{align*}
    together with the limits $(\nabla \psi_n, \nabla a_n, \Lambda_n) \rightarrow (\nabla \psi, \nabla a, \Lambda)$ in $L^2(\R^d)$ are therefore sufficient to conclude that the quantities
    \begin{align*}
        w_\theta^{\frac12} X \text{ for } X \in \{R_h \psi_n, \nabla \psi_n \cdot h, R_h a_n, \nabla a_n \cdot h, \Lambda_n \cdot h\}
    \end{align*}
    are convergent in $L^2(\{|h| \geq \delta\} \times \R^d)$.
    Then the inequalities and identities that we used above imply that
    \begin{align*}
        J_n(x, h) &= C(\theta) \big( |R_h \psi_n(x)| + |R_{-h} \psi_n(x-h)| + |\nabla \psi_n(x) \cdot h| + |\nabla \psi_n(x-h) \cdot h| \big)
        \\ &\qquad\quad \times \big( |R_h a_n(x)| + |R_{-h} a_n(x-h)| + |\nabla a_n(x) \cdot h| + |\nabla a_n(x-h) \cdot h| 
        \\ &\qquad\qquad\quad + |\Lambda_n(x) \cdot h| + |\Lambda_n(x-h) \cdot h| \big)
    \end{align*}
    is an admissible majorant sequence.
\end{proof}

\section{Generalization to $\nabla_{\!\circ,f} \psi$} \label{section:general}
The generalization from $\nabla_{\!\circ} \psi$ to $\nabla_{\!\circ,f} \psi$ is surprisingly straightforward.
\begin{proof}[Proof of Corollary \ref{cor:10}]
    We give the proof for the nontrivial case $s \in (1, \frac32)$.
    Let $x, h \in \R^d$. We calculate
    \begin{align*}
        \Delta_h \nabla_{\!\circ,f} \psi(x) &= \Delta_h f(P)(x) \nabla \psi(x) + f(P)(x-h) \Delta_h \nabla \psi(x) .
    \end{align*}
    Using the Lipschitz continuity of $f$, we can estimate
    \begin{align} \label{eqn:277}
        |\Delta_h f(P)(x) \nabla \psi(x) \cdot h| &\leq [f]_{\Lip} |\Delta_h \conj{P}(x) \nabla \psi(x) \cdot h| 
        \\\nonumber &\leq [f]_{\Lip} \big( |\Delta_h \nabla_{\!\circ} \psi(x) \cdot h| + |\Delta_h \nabla \psi(x) \cdot h| \big) ,
    \end{align}
    while
    \begin{align*}
        |f(P)(x-h) \Delta_h \nabla \psi(x) \cdot h| &\leq \|f\|_{L^\infty} |\Delta_h \nabla \psi(x) \cdot h| .
    \end{align*}
    Then \eqref{eqn:0} implies
    \begin{align*}
        [\nabla_{\!\circ,f} \psi]_{X^{s-1,2}} &\leq \left( \|f\|_{L^\infty} + [f]_{\Lip} \right) [\nabla \psi]_{X^{s-1,2}} + [f]_{\Lip} [\nabla_{\!\circ} \psi]_{X^{s-1,2}}
        \\\nonumber &\leq \left( \|f\|_{L^\infty} + \frac{4}{3-2s} [f]_{\Lip} \right) [\nabla \psi]_{X^{s-1,2}} ,
    \end{align*}
    which proves \eqref{eqn:0-general}.
    For continuity we shall prove that any subsequence has a subsequence which converges to the desired limit. In particular, we may assume without loss of generality that $\psi_n \rightarrow \psi$ and $\nabla \psi_n \rightarrow \nabla \psi$ almost everywhere. Furthermore, since $\nabla \psi_n(x) \rightarrow \nabla \psi(x) = 0$ for almost all $x \in \R^d$ with $\psi(x) = 0$, while $P_n \rightarrow P$ whenever $\psi(x) \neq 0$, we have $\nabla_{\!\circ,f} \psi_n \rightarrow \nabla_{\!\circ,f} \psi$ almost everywhere.
    On the other hand, we know that $w_{s-1}^{\frac12} \Delta_h \nabla \psi_n \cdot h \rightarrow w_{s-1}^{\frac12} \Delta_h \nabla \psi \cdot h$ and $w_{s-1}^{\frac12} \Delta_h \nabla_{\!\circ} \psi_n \cdot h \rightarrow w_{s-1}^{\frac12} \Delta_h \nabla_{\!\circ} \psi \cdot h$ in $L^2(\R^{d+d}; \C)$.
    Then \eqref{eqn:277} together with the above estimate for the second term yields a majorant for the sequence $w_{s-1}^{\frac12} \Delta_h \nabla_{\!\circ,f} \psi_n \cdot h$ in $L^2(\R^{d+d})$, allowing us to close the argument with generalized dominated convergence. Convergence in $L^2$-norm can similarly be obtained by dominated convergence.
\end{proof}

\section{Unboundedness} \label{section:unbounded}

In this section we construct counterexamples exhibiting unboundedness properties.

\subsection{Sharpness of blow-up rates}

The following construction is the basis for our counterexamples at the endpoint regularity $s = \frac32$. Its existence can be understood using results such as \cite[Theorem 1.1]{BjornBjorn2023} on critical condenser capacities of balls.
\begin{lemma} \label{lem:g-n}
    There exist $C > 0$ and odd functions $g_n \in C_c^\infty((-3\eul^{2n}, 3\eul^{2n}); \R)$ for $n \in \N^+$ with $g_n(x) = x$ for $x \in [-2, 2]$ and $[g_n']_{W^{\frac12,2}}^2 \leq \frac{C}{n}$, as well as $\|g_n'\|_{L^2}^2 \leq C \eul^n$.
\end{lemma}
\begin{proof}
    We begin by constructing an even function $g_n' \in C_c^\infty(\R; \R)$ with $\int_{\R} g_n'(x) \dd x = 0$, $g_n'(x) = 1$ for $x \in [-2, 2]$, and $[g_n']_{W^{\frac12,2}}^2 \leq \frac{C}{n}$. To do so, we let $\eta \in C^\infty(\R; [0, 1])$ with $\eta(x) = 1$ if $x \leq 0$ and $\eta(x) = 0$ if $x \geq 1$. Then we define
    \begin{align*}
            f_n'(x) &= \begin{cases}
            1 &, x \in [-2, 2]
            \\ \eta\left(\frac{\log(|x|) - \log(2)}{n}\right) &, x \not\in [-2, 2]
        \end{cases} .
    \end{align*}
    We see that $f_n'(x) = 1$ if $x \in [-2, 2]$ and $\supp f_n' \subset [- 2 \eul^n, 2 \eul^n]$.
    The definition of $f_n'$ extends to a function $F: \mathbb{H} \rightarrow \R$ on the closed upper half-plane $\mathbb{H} = \{z \in \C: \IM[z] \geq 0\}$. Here we have the estimate $|\nabla F(z)| \lesssim \frac{1}{n |z|}$ for $2 < |z| < 2 \eul^n$, and hence $\|\nabla F\|_{L^2(\mathbb{H})}^2 \lesssim \frac{1}{n}$.
    The trace characterization of $\dot{H}^{\frac12}(\R)$ then yields $[f_n']_{W^{\frac12,2}(\R)}^2 \leq \frac{C}{n}$.
    Next, we note that $| \int_{\R} f_n'(y) \dd y| \lesssim \eul^n$ and define
    \begin{align*}
        g_n'(x) &= f_n'(x) - \int_{\R} f_n'(y) \dd y\, \eul^{-2n} \omega(\eul^{-2n} x) ,
    \end{align*}
    where $\omega \in C_c^\infty(\{2<|x|<3\}; [0, 1])$ is an even function with $\int_{\R} \omega(x) \dd x = 1$, so that $\int_{\R} g_n'(x) \dd x = 0$. Due to scale invariance and the presence of the small factor $\eul^{-2n}$, we still have $[g_n']_{W^{\frac12,2}(\R)}^2 \leq \frac{C}{n}$. We define $g_n(x) = \int_{-\infty}^x g_n'(y) \dd y$ and note that $\supp g_n \subset [-3\eul^{2n}, 3\eul^{2n}]$. Since $g_n'$ is even, $g_n$ is odd. One also checks that $\|g_n'\|_{L^2}^2 \leq C \eul^n$.
\end{proof}

\begin{lemma} \label{lem:sharp}
    Let $d \in \N^+$. There exists a sequence $(\psi_n)_{n \geq 4} \subset C^\infty(\R^d; \C)$ such that $[\nabla \psi_n]_{W^{\frac12-\frac{1}{n},2}}^2, [\nabla \psi_n]_{W^{\frac12,2}}^2 \lesssim 1$ while $[\nabla a_n]_{W^{\frac12-\frac{1}{n},2}}^2 \gtrsim n^2$ and $[\Lambda_n]_{W^{\frac12-\frac{1}{n},2}}^2 \gtrsim n$.
\end{lemma}
\begin{proof}
    We consider first the case $d = 1$.
    Let $(g_n)_{n \in \N^+}$ be the sequence constructed in Lemma \ref{lem:g-n} and choose $\chi \in C_c^\infty((-2, 2); [0, 1])$ with $\chi(x) = 1$ for $x \in [-1, 1]$. We define
    \begin{align*}
        \psi_n(x) &= n^{\frac12} (g_n(x) + \imu \eul^{-n} \chi(x)) .
    \end{align*}
    Clearly $[\psi_n']_{W^{\frac12,2}}^2 \lesssim 1$ and $\|\psi_n'\|_{L^2}^2 \lesssim n \eul^n$. Then interpolation yields
    \begin{align*}
        [\psi_n']_{W^{\frac12-\frac{1}{n},2}}^2 &\lesssim \|\psi_n'\|_{L^2}^{\frac{4}{n}} [\psi_n']_{W^{\frac12,2}}^{2-\frac{4}{n}} \lesssim (n \eul^n)^{\frac{2}{n}} \lesssim 1 .
    \end{align*}
    We calculate for $x \in [-1, 1]$ that
    \begin{align*}
        a_n'(x) &= \frac{\eul^n x}{\sqrt{1 + \eul^{2n} x^2}} n^{\frac12}
        & \Lambda_n(x) &= - \frac{1}{\sqrt{1 + \eul^{2n} x^2}} n^{\frac12} .
    \end{align*}
    By scaling, we can calculate
    \begin{align*}
        [a_n']_{W^{\frac12-\frac{1}{n},2}([-1, 1])}^2 &= n \eul^{-2} \left[\frac{x}{\sqrt{1+x^2}}\right]_{W^{\frac12-\frac{1}{n},2}([-\eul^n, \eul^n])}^2 \asymp n^2
        \\ [\Lambda_n]_{W^{\frac12-\frac{1}{n},2}([-1, 1])}^2 &= n \eul^{-2} \left[\frac{1}{\sqrt{1+x^2}}\right]_{W^{\frac12-\frac{1}{n},2}([-\eul^n, \eul^n])}^2 \asymp n .
    \end{align*}
    We move on to the case $d \geq 2$. Let $h \in \mathcal{S}(\R^{d-1}; [0, 1])$ with $h(y) = 1$ for $y \in [-1, 1]^{d-1}$. Define $h_n(y) = r_n^{\frac{d-1}{2}} h(r_n y)$ where $(r_n)_{n \in \N}$ is a rapidly decaying sequence with $r_1 = 1$ that we choose later. Let $\psi^{(1)}_n$ be the sequence constructed for $d = 1$ and define for $(x, y) \in \R \times \R^{d-1}$
    \begin{align*}
        \psi_n &= \psi^{(1)}_n(x) h_n(y) .
    \end{align*}
    Then
    \begin{align*}
        \nabla_{(x,y)} \psi_n(x, y) &= \begin{pmatrix}
            (\psi^{(1)}_n)'(x) h_n(y)
            \\ \psi^{(1)}_n(x) \nabla h_n(y)
        \end{pmatrix} .
    \end{align*}
    We use the Fourier-transform representation of the seminorms with the decomposition $1 = \mathds{1}_{\{|\xi| \leq |\zeta|\}} + \mathds{1}_{\{|\xi| > |\zeta|\}}$ to estimate
    \begin{align*}
        \|(\psi^{(1)}_n)'(x) h_n(y)\|_{\dot{H}^{\frac12}}^2 &\leq \|(\psi^{(1)}_n)'\|_{\dot{H}^{\frac12}}^2 \|h_n\|_{L^2}^2 + \|(\psi^{(1)}_n)'\|_{L^2}^2 \|h_n\|_{\dot{H}^{\frac12}}^2 
        \\ &= \|(\psi^{(1)}_n)'\|_{\dot{H}^{\frac12}}^2 \|h\|_{L^2}^2 + \|(\psi^{(1)}_n)'\|_{L^2}^2 \|h\|_{\dot{H}^{\frac12}}^2 r_n .
        \\ \|\psi^{(1)}_n(x) \nabla h_n(y)\|_{\dot{H}^{\frac12}}^2 &\leq \|\psi^{(1)}_n\|_{\dot{H}^{\frac12}}^2 \|\nabla h_n\|_{L^2}^2 + \|\psi^{(1)}_n\|_{L^2}^2 \|\nabla h_n\|_{\dot{H}^{\frac12}}^2 
        \\ &\leq \|\psi^{(1)}_n\|_{\dot{H}^{\frac12}}^2 \|\nabla h\|_{L^2}^2 r_n^2 + \|\psi^{(1)}_n\|_{L^2}^2 \|\nabla h\|_{\dot{H}^{\frac12}}^2 r_n^3.
    \end{align*}
    By choosing $r_n$ to decay sufficiently quickly, we can ensure that $\|\nabla \psi_n\|_{\dot{H}^{\frac12}}^2 \lesssim 1$. Similarly, we can choose $r_n$ sufficiently small to ensure that $\|\nabla \psi_n\|_{L^2}^2 \lesssim n \eul^n$ still holds. Then the same interpolation argument as above yields $\|\nabla \psi_n\|_{\dot{H}^{\frac12-\frac{1}{n}}}^2 \lesssim 1$.
    Let $\nabla a_n + \imu \Lambda_n = \nabla_{\!\circ} \psi_n$ and $(a^{(1)}_n)' + \imu \Lambda^{(1)}_n = \nabla_{\!\circ} \psi^{(1)}_n$. Conveniently, we have
    \begin{align*}
        \nabla a_n(x, y) + \imu \Lambda_n(x, y) 
        &= \begin{pmatrix}
            ((a^{(1)}_n)'(x) + \imu \Lambda^{(1)}_n(x)) h_n(y) \\ |\psi^{(1)}_n(x)| \nabla h_n(y)
        \end{pmatrix} .
    \end{align*}
    We drop the second component on the right-hand side in the lower bound below.
    Using for frequencies $(\xi, \zeta) \in \R \times \R^{d-1}$ and $\theta \in (0, 1)$ the lower bound $(|\xi|^2 + |\zeta|^2)^\theta \geq |\xi|^{2 \theta}$, we obtain
    \begin{align*}
        \|\nabla a_n\|_{\dot{H}^{\frac12-\frac1n}}^2 &\geq \|(a^{(1)}_n)'\|_{\dot{H}^{\frac12-\frac1n}}^2 \|h_n\|_{L^2}^2 = \|(a^{(1)}_n)'\|_{\dot{H}^{\frac12-\frac1n}}^2 \|h\|_{L^2}^2 \gtrsim n^2
        \\ \|\Lambda_n\|_{\dot{H}^{\frac12-\frac1n}}^2 &\geq \|\Lambda^{(1)}_n\|_{\dot{H}^{\frac12-\frac1n}}^2 \|h_n\|_{L^2}^2 = \|\Lambda^{(1)}_n\|_{\dot{H}^{\frac12-\frac1n}}^2 \|h\|_{L^2}^2 \gtrsim n .
    \end{align*}
\end{proof}

\subsection{Failure of regularity}

\begin{lemma} \label{lem:32}
    Let $d \in \N^+$ and $s \geq \frac32$. There exists $\psi \in C^\infty(\R^d; \C)$ with $\nabla a, \Lambda \in C^\infty(\R^d; \R^d)$ such that $\psi \in W^{s,2}(\R^d)$ but $\nabla a, \Lambda \not\in W^{s-1,2}(\R^d)$.
\end{lemma}
\begin{proof}
    We first assume $d = 1$ and $s = \frac32$.
    Let $(\psi_n)_{n \geq 4}$ be the sequence constructed in Lemma \ref{lem:sharp}.
    We define
    \begin{align*}
        \psi(x) &= \sum_{n=4}^\infty n^{-1} \eul^{-n^4} \psi_n\left(\frac{x - 8 n}{\eul^{-n^4}}\right) .
    \end{align*}
    Due to the disjointness of supports and scale invariance, we have
    \begin{align*}
        [\psi']_{W^{\frac12,2}}^2 &\leq 2 \sum_{n=4}^\infty \left[n^{-1} \psi_n'\left(\frac{x - 8 n}{\eul^{-n^4}}\right)\right]_{W^{\frac12,2}}^2 \leq 2 C \sum_{n=4}^\infty n^{-2} < \infty .
    \end{align*}
    By similarly estimating the lower-regularity norms, one verifies that $\psi \in W^{\frac32,2}(\R)$.

    Note that for $x \in 8 n + \eul^{-n^4} [-1, 1]$ we have
    \begin{align*}
        a'(x) &= n^{-\frac12} \frac{\eul^{n+n^4} (x-8n)}{\sqrt{\eul^{2n+2n^4} (x-8n)^2 + 1}}
        &&& \Lambda(x) &= - n^{-\frac12} \frac{1}{\sqrt{\eul^{2n+2n^4} (x-8n)^2 + 1}} .
    \end{align*}
    The scaling invariance now implies that for a fixed constant $c > 0$
    \begin{align*}
        [a']_{W^{\frac12,2}}^2 &\geq \sum_{n=4}^\infty [a']_{W^{\frac12,2}(8 n + \eul^{-n^4} [-1, 1])}^2 \geq c \sum_{n=4}^\infty \frac{1}{n} = \infty
        \\ [\Lambda]_{W^{\frac12,2}}^2 &\geq \sum_{n=4}^\infty [\Lambda]_{W^{\frac12,2}(8 n + \eul^{-n^4} [-1, 1])}^2 \geq c \sum_{n=4}^\infty \frac{1}{n} = \infty .
    \end{align*}
    This counterexample can be transferred to the cases $d \geq 2$ as in the proof of Lemma \ref{lem:sharp}, and it suffices to take $r_n = 1$ for the verification of $\psi \in W^{\frac32,2}(\R^d)$ since an inhomogeneous space is used.
    
    Now we assume $d = 1$, $s \in (\frac32, 2)$, and set $\delta = s - \frac32$. 
    We use the function $\chi$ from the proof of Lemma \ref{lem:sharp} to define
    \begin{align} \label{eqn:phi}
        \phi_n(x) &= (x + \imu \eul^{-n}) \chi(x) .
    \end{align}
    We define
    \begin{align*}
        \psi(x) &= \sum_{n=4}^\infty n^{-\frac12} \eul^{- \delta n} \phi_n(x-8n) .
    \end{align*}
    Due to the disjointness of supports, we have
    \begin{align*}
        [\psi']_{W^{s-1,2}}^2 &\leq 2 \sum_{n=4}^\infty \left[n^{-\frac12} \eul^{- \delta n} \phi_n'(x-8n)\right]_{W^{s-1,2}}^2
        \\ &= 2 \sum_{n=4}^\infty n^{-1} \eul^{- 2 \delta n} \left[ \phi_n'\right]_{W^{s-1,2}}^2 \lesssim \sum_{n=4}^\infty n^{-1} \eul^{- 2 \delta n} < \infty .
    \end{align*}
    Similarly, the lower-regularity norms are finite and we conclude $\psi \in W^{s,2}(\R)$.
    For $x \in 8 n + [-1, 1]$ we can calculate
    \begin{align*}
        a'(x) &= n^{-\frac12} \eul^{- \delta n} \frac{\eul^n (x-8n)}{\sqrt{\eul^{2n} (x-8n)^2 + 1}}
        &&& \Lambda(x) &= - n^{-\frac12} \eul^{- \delta n} \frac{1}{\sqrt{\eul^{2n} (x-8n)^2 + 1}} .
    \end{align*}
    Then
    \begin{align*}
        [a']_{W^{s-1,2}}^2 &\geq \sum_{n=4}^\infty [a']_{W^{s-1,2}(8 n + [-1, 1])}^2 \geq c(s) \sum_{n=4}^\infty n^{-1} \eul^{- 2 \delta n} \eul^{2 (s-1) n - n} = \infty
        \\ [\Lambda]_{W^{s-1,2}}^2 &\geq \sum_{n=4}^\infty [\Lambda]_{W^{s-1,2}(8 n + [-1, 1])}^2 \geq c(s) \sum_{n=4}^\infty n^{-1} \eul^{- 2 \delta n} \eul^{2 (s-1) n - n} = \infty .
    \end{align*}
    The generalization to $d \geq 2$ works as before. Furthermore, the same examples can be used for the cases $s \geq 2$.
\end{proof}

\begin{lemma} \label{lem:endpoint-cont}
    Let $d \in \N^+$. There exists a convergent sequence $\psi_n \rightarrow \psi$ in $B^{\frac32}_{2,1}(\R^d; \C)$ for which $\nabla a_n \not\rightarrow \nabla a$ in $B^{\frac12}_{2,\infty}(\R^d; \R^d)$ and $\Lambda_n \not\rightarrow 0$ in $W^{\frac12,2}(\R^d; \R^d)$.
\end{lemma}
\begin{proof}
    For $d = 1$ we choose $\psi_n = \phi_n$ from \eqref{eqn:phi}, which satisfies $\phi_n \rightarrow x \chi$ in $B^{\frac32}_{2,1}(\R; \C)$.
    Due to scaling invariance the sequences $(\nabla a_n)_{n \in \N}$ and $(\Lambda_n)_{n \in \N}$ fail to converge respectively in $B^{\frac12}_{2,\infty}(\R^d)$ and $W^{\frac12,2}(\R^d)$, even though they converge pointwise almost everywhere.
    This example can be transferred to the cases $d \geq 2$ as before. 
\end{proof}

\begin{lemma} \label{lem:32-besov}
    Let $d \in \N^+$. There exists $\psi \in C^\infty(\R^d; \C)$ with $\nabla a, \Lambda \in C^\infty(\R^d; \R^d)$ such that $\psi \in W^{2,2}(\R^d)$ but $\nabla a \not\in W^{\frac12,2}(\R^d)$ and $\Lambda \not\in B^{\frac12}_{2,1}(\R^d)$.
\end{lemma}
\begin{proof}
    We first assume $d = 1$.
    Let $(\phi_n)_{n \in \N}$ be the sequence constructed in the proof of Lemma \ref{lem:32}.
    We define
    \begin{align*}
        \psi(x) &= \sum_{n=4}^\infty n^{-1} \phi_n(x - 8 n)
    \end{align*}
    and note that, due to the disjointness of supports in the sum, the inclusion $x \chi, \chi \in C_c^\infty(\R; \R) \subset W^{2,2}(\R)$ implies $\psi \in W^{2,2}(\R)$.
    For $x \in 8 n + [-1, 1]$ we calculate
    \begin{align*}
        a'(x) &= n^{-1} \frac{\eul^n (x-8n)}{\sqrt{\eul^{2n} (x-8n)^2 + 1}}
        &&& \Lambda(x) &= - n^{-1} \frac{1}{\sqrt{\eul^{2n} (x-8n)^2 + 1}} .
    \end{align*}
    There exists a constant $c_1 > 0$ such that for $x \in 8 n + \eul^{-n} [-\frac12, 0]$, $h \in [\eul^{-n}, 2 \eul^{-n}]$, and $n \geq 4$ we have
    \begin{align*}
        |\Delta_h \Lambda(x)| &= n^{-1} \left| \frac{1}{\sqrt{\eul^{2n} (x-8n)^2 + 1}} - \frac{1}{\sqrt{\eul^{2n} (x-h-8n)^2 +1}} \right| \geq c_1 n^{-1} .
    \end{align*}
    In particular
    \begin{align*}
        \|\Delta_h \Lambda\|_{L^2} &\geq \|\Delta_h \Lambda\|_{L^2(8 n + \eul^{-n} [-\frac12,0])} \geq c n^{-1} \eul^{-\frac{n}{2}} ,
    \end{align*}
    for some constant $c > 0$.
    Then
    \begin{align*}
        [\Lambda]_{B^{\frac12}_{2,1}} &= \int_{-1}^1 |h|^{-\frac32} \|\Delta_h \Lambda\|_{L^2} \dd h
        \geq \sum_{n=4}^\infty \int_{\eul^{-n}}^{2 \eul^{-n}} |h|^{-\frac32} c n^{-1} \eul^{-\frac{n}{2}} \dd h
        \\ &\geq c \sum_{n=4}^\infty n^{-1} \eul^{-\frac{n}{2}} \int_{\eul^{-n}}^{2 \eul^{-n}} h^{-\frac32} \dd h
        \geq c \sum_{n=4}^\infty n^{-1} = \infty .
    \end{align*}
    Estimating $[a']_{W^{\frac12,2}(8n+[-1,1])}^2 \gtrsim n^{-1}$ and using again the disjointness of supports, we find that also $a' \not\in W^{\frac12,2}(\R)$.
    For the generalization to $d \geq 2$, we proceed as before, using now the Littlewood--Paley characterization of the Besov norms.
\end{proof}

\bibliographystyle{siam}
\bibliography{bibliomain}            

@Article{BourdaudMeyer,
    AUTHOR = {Bourdaud, G\'erard and Meyer, Yves},
     TITLE = {Fonctions qui op\`erent sur les espaces de {S}obolev},
   JOURNAL = {J. Funct. Anal.},
  FJOURNAL = {Journal of Functional Analysis},
    VOLUME = {97},
      YEAR = {1991},
    NUMBER = {2},
     PAGES = {351--360},
      ISSN = {0022-1236,1096-0783},
   MRCLASS = {46E35 (47B38)},
  MRNUMBER = {1111186},
MRREVIEWER = {Nicolas\ Lerner},
       DOI = {10.1016/0022-1236(91)90006-Q},
       URL = {https://doi.org/10.1016/0022-1236(91)90006-Q},
}

@incollection{BourdaudSickel,
    AUTHOR = {Bourdaud, G\'erard and Sickel, Winfried},
     TITLE = {Composition operators on function spaces with fractional order
              of smoothness},
 BOOKTITLE = {Harmonic analysis and nonlinear partial differential
              equations},
    SERIES = {RIMS K\^oky\^uroku Bessatsu},
    VOLUME = {B26},
     PAGES = {93--132},
 PUBLISHER = {Res. Inst. Math. Sci. (RIMS), Kyoto},
      YEAR = {2011},
   MRCLASS = {47B33 (35J20 35J91 46E35 47H30)},
  MRNUMBER = {2883848},
MRREVIEWER = {Pamela\ Gorkin},
}

@article {Li,
    AUTHOR = {Li, Dong},
     TITLE = {On fractional smoothness of modulus of functions},
   JOURNAL = {Ann. Appl. Math.},
  FJOURNAL = {Annals of Applied Mathematics. Yingyong Shuxue Niankan},
    VOLUME = {37},
      YEAR = {2021},
    NUMBER = {3},
     PAGES = {394--404},
      ISSN = {2096-0174},
   MRCLASS = {35R11},
  MRNUMBER = {4320610},
       DOI = {10.4208/aam.oa-2021-0006},
       URL = {https://doi.org/10.4208/aam.oa-2021-0006},
}

@book {LL01,
    AUTHOR = {Lieb, Elliott H. and Loss, Michael},
     TITLE = {Analysis},
    SERIES = {Graduate Studies in Mathematics},
    VOLUME = {14},
   EDITION = {Second},
 PUBLISHER = {American Mathematical Society, Providence, RI},
      YEAR = {2001},
     PAGES = {xxii+346},
      ISBN = {0-8218-2783-9},
   MRCLASS = {00A05 (26-01 28-01 31-01 35J10 42-01)},
  MRNUMBER = {1817225},
       DOI = {10.1090/gsm/014},
       URL = {https://doi.org/10.1090/gsm/014},
}

@article {MusinaNazarov,
    AUTHOR = {Musina, Roberta and Nazarov, Alexander I.},
     TITLE = {A note on truncations in fractional {S}obolev spaces},
   JOURNAL = {Bull. Math. Sci.},
  FJOURNAL = {Bulletin of Mathematical Sciences},
    VOLUME = {9},
      YEAR = {2019},
    NUMBER = {1},
     PAGES = {1950001, 7},
      ISSN = {1664-3607,1664-3615},
   MRCLASS = {46E35 (47H30)},
  MRNUMBER = {3949429},
MRREVIEWER = {Jos\'e\ Francisco Alves de Oliveira},
       DOI = {10.1142/S1664360719500012},
       URL = {https://doi.org/10.1142/S1664360719500012},
}

@article {Oswald,
    AUTHOR = {Oswald, Peter},
     TITLE = {On the boundedness of the mapping {$f\to|f|$} in {B}esov
              spaces},
   JOURNAL = {Comment. Math. Univ. Carolin.},
  FJOURNAL = {Commentationes Mathematicae Universitatis Carolinae},
    VOLUME = {33},
      YEAR = {1992},
    NUMBER = {1},
     PAGES = {57--66},
      ISSN = {0010-2628,1213-7243},
   MRCLASS = {46E35 (41A15 47H30)},
  MRNUMBER = {1173747},
}

@article {ScottMengesha2019,
    AUTHOR = {Scott, James and Mengesha, Tadele},
     TITLE = {A fractional {K}orn-type inequality},
   JOURNAL = {Discrete Contin. Dyn. Syst.},
  FJOURNAL = {Discrete and Continuous Dynamical Systems},
    VOLUME = {39},
      YEAR = {2019},
    NUMBER = {6},
     PAGES = {3315--3343},
      ISSN = {1078-0947,1553-5231},
   MRCLASS = {46E35 (35A23 35R11 45G15 74B99)},
  MRNUMBER = {3959431},
MRREVIEWER = {Aapo\ Kauranen},
       DOI = {10.3934/dcds.2019137},
       URL = {https://doi.org/10.3934/dcds.2019137},
}

@book {Leoni2023,
    AUTHOR = {Leoni, Giovanni},
     TITLE = {A first course in fractional {S}obolev spaces},
    SERIES = {Graduate Studies in Mathematics},
    VOLUME = {229},
 PUBLISHER = {American Mathematical Society, Providence, RI},
      YEAR = {2023},
     PAGES = {xv+586},
      ISBN = {9781470468989; 9781470472535; 9781470472528},
   MRCLASS = {46-01 (30H05 35R11 42Bxx 42C40 46E35)},
  MRNUMBER = {4567945},
MRREVIEWER = {E.\ S.\ Dubtsov},
       DOI = {10.1090/gsm/229},
       URL = {https://doi.org/10.1090/gsm/229},
}

@book {BCD,
    AUTHOR = {Bahouri, Hajer and Chemin, Jean-Yves and Danchin, Rapha\"el},
     TITLE = {Fourier analysis and nonlinear partial differential equations},
    SERIES = {Grundlehren der mathematischen Wissenschaften [Fundamental
              Principles of Mathematical Sciences]},
    VOLUME = {343},
 PUBLISHER = {Springer, Heidelberg},
      YEAR = {2011},
     PAGES = {xvi+523},
      ISBN = {978-3-642-16829-1},
   MRCLASS = {35-02 (35L72 35Q30 42-02 42B37 76B03 76D03 76N10)},
  MRNUMBER = {2768550},
MRREVIEWER = {Peter\ R.\ Massopust},
       DOI = {10.1007/978-3-642-16830-7},
       URL = {https://doi.org/10.1007/978-3-642-16830-7},
}

@article {Gerard2006,
    AUTHOR = {G\'{e}rard, Patrick},
     TITLE = {The {C}auchy problem for the {G}ross-{P}itaevskii equation},
   JOURNAL = {Ann. Inst. H. Poincar\'{e} Anal. Non Lin\'{e}aire},
  FJOURNAL = {Annales de l'Institut Henri Poincar\'{e}. Analyse Non Lin\'{e}aire},
    VOLUME = {23},
      YEAR = {2006},
    NUMBER = {5},
     PAGES = {765--779},
      ISSN = {0294-1449},
   MRCLASS = {35Q55 (35A15 35B30 35J20 81Q20)},
  MRNUMBER = {2259616},
MRREVIEWER = {Zhi-Qiang Wang},
       DOI = {10.1016/j.anihpc.2005.09.004},
       URL = {https://doi.org/10.1016/j.anihpc.2005.09.004},
}

@incollection {Gerard2008,
    AUTHOR = {G\'erard, Patrick},
     TITLE = {The {G}ross-{P}itaevskii equation in the energy space},
 BOOKTITLE = {Stationary and time dependent {G}ross-{P}itaevskii equations},
    SERIES = {Contemp. Math.},
    VOLUME = {473},
     PAGES = {129--148},
 PUBLISHER = {Amer. Math. Soc., Providence, RI},
      YEAR = {2008},
      ISBN = {978-0-8218-4357-4},
   MRCLASS = {35Q55 (35J20)},
  MRNUMBER = {2522016},
MRREVIEWER = {Ayman\ Kachmar},
       DOI = {10.1090/conm/473/09226},
       URL = {https://doi.org/10.1090/conm/473/09226},
}

@article {BethuelGravejatSautSmets,
    AUTHOR = {B\'ethuel, Fabrice and Gravejat, Philippe and Saut,
              Jean-Claude and Smets, Didier},
     TITLE = {Orbital stability of the black soliton for the
              {G}ross-{P}itaevskii equation},
   JOURNAL = {Indiana Univ. Math. J.},
  FJOURNAL = {Indiana University Mathematics Journal},
    VOLUME = {57},
      YEAR = {2008},
    NUMBER = {6},
     PAGES = {2611--2642},
      ISSN = {0022-2518,1943-5258},
   MRCLASS = {35Q55 (35B35 35C07 35C08)},
  MRNUMBER = {2482993},
MRREVIEWER = {Hiroaki\ Kikuchi},
       DOI = {10.1512/iumj.2008.57.3632},
       URL = {https://doi.org/10.1512/iumj.2008.57.3632},
}

@article {BGSasympstab,
    AUTHOR = {Bethuel, Fabrice and Gravejat, Philippe and Smets, Didier},
     TITLE = {Asymptotic stability in the energy space for dark solitons of
              the {G}ross-{P}itaevskii equation},
   JOURNAL = {Ann. Sci. \'Ec. Norm. Sup\'er. (4)},
  FJOURNAL = {Annales Scientifiques de l'\'Ecole Normale Sup\'erieure.
              Quatri\`eme S\'erie},
    VOLUME = {48},
      YEAR = {2015},
    NUMBER = {6},
     PAGES = {1327--1381},
      ISSN = {0012-9593,1873-2151},
   MRCLASS = {35Q55 (35B35 35C08 37K40)},
  MRNUMBER = {3429470},
MRREVIEWER = {Amin\ Esfahani},
       DOI = {10.24033/asens.2271},
       URL = {https://doi.org/10.24033/asens.2271},
}

@article {BGSstab,
    AUTHOR = {B\'ethuel, Fabrice and Gravejat, Philippe and Smets, Didier},
     TITLE = {Stability in the energy space for chains of solitons of the
              one-dimensional {G}ross-{P}itaevskii equation},
   JOURNAL = {Ann. Inst. Fourier (Grenoble)},
  FJOURNAL = {Universit\'e{} de Grenoble. Annales de l'Institut Fourier},
    VOLUME = {64},
      YEAR = {2014},
    NUMBER = {1},
     PAGES = {19--70},
      ISSN = {0373-0956,1777-5310},
   MRCLASS = {35Q55 (35B35 35C07 35C08)},
  MRNUMBER = {3330540},
MRREVIEWER = {Piotr\ P.\ Goldstein},
       DOI = {10.5802/aif.2838},
       URL = {https://doi.org/10.5802/aif.2838},
}

@article {AudiardHaspot,
    AUTHOR = {Audiard, Corentin and Haspot, Boris},
     TITLE = {From the {G}ross-{P}itaevskii equation to the {E}uler
              {K}orteweg system, existence of global strong solutions with
              small irrotational initial data},
   JOURNAL = {Ann. Mat. Pura Appl. (4)},
  FJOURNAL = {Annali di Matematica Pura ed Applicata. Series IV},
    VOLUME = {197},
      YEAR = {2018},
    NUMBER = {3},
     PAGES = {721--760},
      ISSN = {0373-3114,1618-1891},
   MRCLASS = {35Q35 (35A01 35A02 35B35 35B40 35P25 35Q53)},
  MRNUMBER = {3802692},
MRREVIEWER = {Nikolay\ G.\ Kuznetsov},
       DOI = {10.1007/s10231-017-0702-z},
       URL = {https://doi.org/10.1007/s10231-017-0702-z},
}

@article {Merle,
    AUTHOR = {Merle, Frank and Rapha\"el, Pierre and Rodnianski, Igor and
              Szeftel, Jeremie},
     TITLE = {On blow up for the energy super critical defocusing nonlinear
              {S}chr\"odinger equations},
   JOURNAL = {Invent. Math.},
  FJOURNAL = {Inventiones Mathematicae},
    VOLUME = {227},
      YEAR = {2022},
    NUMBER = {1},
     PAGES = {247--413},
      ISSN = {0020-9910,1432-1297},
   MRCLASS = {35Q55},
  MRNUMBER = {4359478},
       DOI = {10.1007/s00222-021-01067-9},
       URL = {https://doi.org/10.1007/s00222-021-01067-9},
}

@article {AntonelliMarcati,
    AUTHOR = {Antonelli, Paolo and Marcati, Pierangelo},
     TITLE = {On the finite energy weak solutions to a system in quantum
              fluid dynamics},
   JOURNAL = {Comm. Math. Phys.},
  FJOURNAL = {Communications in Mathematical Physics},
    VOLUME = {287},
      YEAR = {2009},
    NUMBER = {2},
     PAGES = {657--686},
      ISSN = {0010-3616,1432-0916},
   MRCLASS = {82D50 (35D30 35Q35 76Y05 82D37 82D55)},
  MRNUMBER = {2481754},
       DOI = {10.1007/s00220-008-0632-0},
       URL = {https://doi.org/10.1007/s00220-008-0632-0},
}

@article {AntonelliMarcati2012,
    AUTHOR = {Antonelli, Paolo and Marcati, Pierangelo},
     TITLE = {The quantum hydrodynamics system in two space dimensions},
   JOURNAL = {Arch. Ration. Mech. Anal.},
  FJOURNAL = {Archive for Rational Mechanics and Analysis},
    VOLUME = {203},
      YEAR = {2012},
    NUMBER = {2},
     PAGES = {499--527},
      ISSN = {0003-9527,1432-0673},
   MRCLASS = {82D15 (35Q35)},
  MRNUMBER = {2885568},
MRREVIEWER = {Mahendra\ Panthee},
       DOI = {10.1007/s00205-011-0454-7},
       URL = {https://doi.org/10.1007/s00205-011-0454-7},
}

@incollection {AntonelliMarcati2026,
    AUTHOR = {Antonelli, Paolo and Marcati, Pierangelo},
     TITLE = {An introduction to the mathematical theory of quantum fluids},
 BOOKTITLE = {New perspectives in mathematical fluid mechanics},
    SERIES = {Ne\v cas Center Ser.},
     PAGES = {3--85},
 PUBLISHER = {Birkh\"auser/Springer, Cham},
      YEAR = {2026},
      ISBN = {978-3-032-02318-6; 978-3-032-02319-3},
   MRCLASS = {76Y05},
  MRNUMBER = {5057465},
       DOI = {10.1007/978-3-032-02319-3\_1},
       URL = {https://doi.org/10.1007/978-3-032-02319-3_1},
}

@article {AntonelliMarcatiZheng2,
    AUTHOR = {Antonelli, Paolo and Marcati, Pierangelo and Zheng, Hao},
     TITLE = {An intrinsically hydrodynamic approach to multidimensional
              {QHD} systems},
   JOURNAL = {Arch. Ration. Mech. Anal.},
  FJOURNAL = {Archive for Rational Mechanics and Analysis},
    VOLUME = {247},
      YEAR = {2023},
    NUMBER = {2},
     PAGES = {Paper No. 24, 58},
      ISSN = {0003-9527,1432-0673},
   MRCLASS = {76Y05},
  MRNUMBER = {4562812},
MRREVIEWER = {Delyan\ Zhelyazov},
       DOI = {10.1007/s00205-023-01856-x},
       URL = {https://doi.org/10.1007/s00205-023-01856-x},
}

@article {AntonelliMarcatiZheng1,
    AUTHOR = {Antonelli, Paolo and Marcati, Pierangelo and Zheng, Hao},
     TITLE = {Genuine hydrodynamic analysis to the 1-{D} {QHD} system:
              existence, dispersion and stability},
   JOURNAL = {Comm. Math. Phys.},
  FJOURNAL = {Communications in Mathematical Physics},
    VOLUME = {383},
      YEAR = {2021},
    NUMBER = {3},
     PAGES = {2113--2161},
      ISSN = {0010-3616,1432-0916},
   MRCLASS = {35Q40 (76Y05 82D50)},
  MRNUMBER = {4244267},
       DOI = {10.1007/s00220-021-03998-z},
       URL = {https://doi.org/10.1007/s00220-021-03998-z},
}

@article {ConstantinSaut,
    AUTHOR = {Constantin, Peter and Saut, Jean-Claude},
     TITLE = {Local smoothing properties of dispersive equations},
   JOURNAL = {J. Amer. Math. Soc.},
  FJOURNAL = {Journal of the American Mathematical Society},
    VOLUME = {1},
      YEAR = {1988},
    NUMBER = {2},
     PAGES = {413--439},
      ISSN = {0894-0347,1088-6834},
   MRCLASS = {35Q20 (35D10)},
  MRNUMBER = {928265},
MRREVIEWER = {Ralph\ Saxton},
       DOI = {10.2307/1990923},
       URL = {https://doi.org/10.2307/1990923},
}

@incollection {FrankSeiringer2010,
    AUTHOR = {Frank, Rupert L. and Seiringer, Robert},
     TITLE = {Sharp fractional {H}ardy inequalities in half-spaces},
 BOOKTITLE = {Around the research of {V}ladimir {M}az'ya. {I}},
    SERIES = {Int. Math. Ser. (N. Y.)},
    VOLUME = {11},
     PAGES = {161--167},
 PUBLISHER = {Springer, New York},
      YEAR = {2010},
      ISBN = {978-1-4419-1340-1; 978-5-9018-7341-0},
   MRCLASS = {35A23 (35J20 46E35)},
  MRNUMBER = {2723817},
MRREVIEWER = {Rodrigo\ da Silva Rodrigues},
       DOI = {10.1007/978-1-4419-1341-8\_6},
       URL = {https://doi.org/10.1007/978-1-4419-1341-8_6},
}

@article {BrascoGomezCastroVazquez2021,
    AUTHOR = {Brasco, Lorenzo and G\'omez-Castro, David and V\'azquez, Juan
              Luis},
     TITLE = {Characterisation of homogeneous fractional {S}obolev spaces},
   JOURNAL = {Calc. Var. Partial Differential Equations},
  FJOURNAL = {Calculus of Variations and Partial Differential Equations},
    VOLUME = {60},
      YEAR = {2021},
    NUMBER = {2},
     PAGES = {Paper No. 60, 40},
      ISSN = {0944-2669,1432-0835},
   MRCLASS = {46E35},
  MRNUMBER = {4225499},
MRREVIEWER = {Artur\ Rutkowski},
       DOI = {10.1007/s00526-021-01934-6},
       URL = {https://doi.org/10.1007/s00526-021-01934-6},
}

@article {DiNezzaPalatucciValdinoci2012,
    AUTHOR = {Di Nezza, Eleonora and Palatucci, Giampiero and Valdinoci,
              Enrico},
     TITLE = {Hitchhiker's guide to the fractional {S}obolev spaces},
   JOURNAL = {Bull. Sci. Math.},
  FJOURNAL = {Bulletin des Sciences Math\'ematiques},
    VOLUME = {136},
      YEAR = {2012},
    NUMBER = {5},
     PAGES = {521--573},
      ISSN = {0007-4497,1952-4773},
   MRCLASS = {46E35 (35A23 35S05 35S30)},
  MRNUMBER = {2944369},
MRREVIEWER = {Lanzhe\ Liu},
       DOI = {10.1016/j.bulsci.2011.12.004},
       URL = {https://doi.org/10.1016/j.bulsci.2011.12.004},
}

@book {EvansGariepy,
    AUTHOR = {Evans, Lawrence C. and Gariepy, Ronald F.},
     TITLE = {Measure theory and fine properties of functions},
    SERIES = {Textbooks in Mathematics},
   EDITION = {Second},
 PUBLISHER = {CRC Press, Boca Raton, FL},
      YEAR = {2025},
     PAGES = {xi+327},
      ISBN = {978-1-032-94644-3; 978-1-003-58300-4; 978-1-032-95055-6},
   MRCLASS = {28-01},
  MRNUMBER = {4907811},
       DOI = {10.1201/9781003583004}
}

@article {BogdanDyda2011,
    AUTHOR = {Bogdan, Krzysztof and Dyda, Bart\l omiej},
     TITLE = {The best constant in a fractional {H}ardy inequality},
   JOURNAL = {Math. Nachr.},
  FJOURNAL = {Mathematische Nachrichten},
    VOLUME = {284},
      YEAR = {2011},
    NUMBER = {5-6},
     PAGES = {629--638},
      ISSN = {0025-584X,1522-2616},
   MRCLASS = {26D15 (31B25 31C25 46E35)},
  MRNUMBER = {2663757},
MRREVIEWER = {Zheng\ Liu},
       DOI = {10.1002/mana.200810109},
       URL = {https://doi.org/10.1002/mana.200810109},
}

@book {Evans2010,
    AUTHOR = {Evans, Lawrence C.},
     TITLE = {Partial differential equations},
    SERIES = {Graduate Studies in Mathematics},
    VOLUME = {19},
   EDITION = {Second},
 PUBLISHER = {American Mathematical Society, Providence, RI},
      YEAR = {2010},
     PAGES = {xxii+749},
      ISBN = {978-0-8218-4974-3},
   MRCLASS = {35-01},
  MRNUMBER = {2597943},
MRREVIEWER = {Diego\ M.\ Maldonado},
       DOI = {10.1090/gsm/019},
       URL = {https://doi.org/10.1090/gsm/019},
}

@article{sympy,
  title = {SymPy: symbolic computing in Python},
  author = {A. Meurer and C. P. Smith and M. Paprocki and O. \v{C}ert\'{i}k and S. B. Kirpichev and M. Rocklin and A. Kumar and S. Ivanov and J. K. Moore and S. Singh and T. Rathnayake and S. Vig and B. E. Granger and R. P. Muller and F. Bonazzi and H. Gupta and S. Vats and F. Johansson and F. Pedregosa and M. J. Curry and A. R. Terrel and \v{S}. Rou\v{c}ka and A. Saboo and I. Fernando and S. Kulal and R. Cimrman and A. Scopatz},
  year = {2017},
  month = {01},
  volume = {3},
  pages = {e103},
  journal = {PeerJ Comput. Sci.},
  issn = {2376-5992},
  url = {https://doi.org/10.7717/peerj-cs.103},
  doi = {10.7717/peerj-cs.103}
}

@article {Kateb2003,
    AUTHOR = {Kateb, Djalil},
     TITLE = {On the boundedness of the mapping {$f\mapsto |f|^\mu,\mu>1$}
              on {B}esov spaces},
   JOURNAL = {Math. Nachr.},
  FJOURNAL = {Mathematische Nachrichten},
    VOLUME = {248/249},
      YEAR = {2003},
     PAGES = {110--128},
      ISSN = {0025-584X,1522-2616},
   MRCLASS = {46E35 (46M35 47H99)},
  MRNUMBER = {1950719},
MRREVIEWER = {Tohru\ Ozawa},
       DOI = {10.1002/mana.200310007},
       URL = {https://doi.org/10.1002/mana.200310007},
}

@article {Kateb2000,
    AUTHOR = {Kateb, Djalil},
     TITLE = {Fonctions qui op\`erent sur les espaces de {B}esov},
   JOURNAL = {Proc. Amer. Math. Soc.},
  FJOURNAL = {Proceedings of the American Mathematical Society},
    VOLUME = {128},
      YEAR = {2000},
    NUMBER = {3},
     PAGES = {735--743},
      ISSN = {0002-9939,1088-6826},
   MRCLASS = {47H30 (46E35)},
  MRNUMBER = {1636962},
       DOI = {10.1090/S0002-9939-99-05096-0},
       URL = {https://doi.org/10.1090/S0002-9939-99-05096-0},
}

@article {Moussai2018,
    AUTHOR = {Moussai, Madani},
     TITLE = {Composition operators on {B}esov spaces in the limiting case
              {$s=1+1/p$}},
   JOURNAL = {Studia Math.},
  FJOURNAL = {Studia Mathematica},
    VOLUME = {241},
      YEAR = {2018},
    NUMBER = {1},
     PAGES = {1--15},
      ISSN = {0039-3223,1730-6337},
   MRCLASS = {46E35 (47H30)},
  MRNUMBER = {3732927},
MRREVIEWER = {Yoichi\ Miyazaki},
       DOI = {10.4064/sm8136-4-2017},
       URL = {https://doi.org/10.4064/sm8136-4-2017},
}

@article {MironescuRussSire2020,
    AUTHOR = {Mironescu, Petru and Russ, Emmanuel and Sire, Yannick},
     TITLE = {Lifting in {B}esov spaces},
   JOURNAL = {Nonlinear Anal.},
  FJOURNAL = {Nonlinear Analysis. Theory, Methods \& Applications. An
              International Multidisciplinary Journal},
    VOLUME = {193},
      YEAR = {2020},
     PAGES = {111489, 44},
      ISSN = {0362-546X,1873-5215},
   MRCLASS = {46E35},
  MRNUMBER = {4062976},
MRREVIEWER = {Oscar\ Dom\'inguez},
       DOI = {10.1016/j.na.2019.03.012},
       URL = {https://doi.org/10.1016/j.na.2019.03.012},
}

@article {BourgainBrezisMironescu2000,
    AUTHOR = {Bourgain, Jean and Brezis, Haim and Mironescu, Petru},
     TITLE = {Lifting in {S}obolev spaces},
   JOURNAL = {J. Anal. Math.},
  FJOURNAL = {Journal d'Analyse Math\'ematique},
    VOLUME = {80},
      YEAR = {2000},
     PAGES = {37--86},
      ISSN = {0021-7670,1565-8538},
   MRCLASS = {46E35 (35J20 58D15 58E20)},
  MRNUMBER = {1771523},
MRREVIEWER = {Piotr\ Haj\l asz},
       DOI = {10.1007/BF02791533},
       URL = {https://doi.org/10.1007/BF02791533},
}

@article {BrezisMironescu2001,
    AUTHOR = {Brezis, Haim and Mironescu, Petru},
     TITLE = {Composition in fractional {S}obolev spaces},
   JOURNAL = {Discrete Contin. Dynam. Systems},
  FJOURNAL = {Discrete and Continuous Dynamical Systems},
    VOLUME = {7},
      YEAR = {2001},
    NUMBER = {2},
     PAGES = {241--246},
      ISSN = {1078-0947,1553-5231},
   MRCLASS = {46E35 (47B33)},
  MRNUMBER = {1808397},
MRREVIEWER = {Miguel\ A.\ Sanz Alix},
       DOI = {10.3934/dcds.2001.7.241},
       URL = {https://doi.org/10.3934/dcds.2001.7.241},
}

@article {MarcusMizel1979,
    AUTHOR = {Marcus, Moshe and Mizel, Victor J.},
     TITLE = {Every superposition operator mapping one {S}obolev space into
              another is continuous},
   JOURNAL = {J. Functional Analysis},
  FJOURNAL = {Journal of Functional Analysis},
    VOLUME = {33},
      YEAR = {1979},
    NUMBER = {2},
     PAGES = {217--229},
      ISSN = {0022-1236},
   MRCLASS = {47B38 (46E35)},
  MRNUMBER = {546508},
MRREVIEWER = {Vo-Khac Khoan},
       DOI = {10.1016/0022-1236(79)90113-7},
       URL = {https://doi.org/10.1016/0022-1236(79)90113-7},
}

@article {PonceSpector2020,
    AUTHOR = {Ponce, Augusto C. and Spector, Daniel},
     TITLE = {A decomposition by non-negative functions in the {S}obolev
              space {$W^{k,1}$}},
   JOURNAL = {Indiana Univ. Math. J.},
  FJOURNAL = {Indiana University Mathematics Journal},
    VOLUME = {69},
      YEAR = {2020},
    NUMBER = {1},
     PAGES = {151--169},
      ISSN = {0022-2518,1943-5258},
   MRCLASS = {46E35},
  MRNUMBER = {4077159},
MRREVIEWER = {Florin\ Isaia},
       DOI = {10.1512/iumj.2020.69.8237},
       URL = {https://doi.org/10.1512/iumj.2020.69.8237},
}

@article {QiangKun2011,
    AUTHOR = {Du, Qiang and Zhou, Kun},
     TITLE = {Mathematical analysis for the peridynamic nonlocal continuum
              theory},
   JOURNAL = {ESAIM Math. Model. Numer. Anal.},
  FJOURNAL = {ESAIM. Mathematical Modelling and Numerical Analysis},
    VOLUME = {45},
      YEAR = {2011},
    NUMBER = {2},
     PAGES = {217--234},
      ISSN = {2822-7840,2804-7214},
   MRCLASS = {45A05 (35R09 74B99)},
  MRNUMBER = {2804637},
MRREVIEWER = {Rasajit\ Kumar\ Bera},
       DOI = {10.1051/m2an/2010040},
       URL = {https://doi.org/10.1051/m2an/2010040},
}

@book {Triebel1983,
    AUTHOR = {Triebel, Hans},
     TITLE = {Theory of function spaces},
    SERIES = {Monographs in Mathematics},
    VOLUME = {78},
 PUBLISHER = {Birkh\"auser Verlag, Basel},
      YEAR = {1983},
     PAGES = {284},
      ISBN = {3-7643-1381-1},
   MRCLASS = {46Exx},
  MRNUMBER = {781540},
       DOI = {10.1007/978-3-0346-0416-1},
       URL = {https://doi.org/10.1007/978-3-0346-0416-1},
}

@article {BjornBjorn2023,
    AUTHOR = {Bj\"orn, Anders and Bj\"orn, Jana},
     TITLE = {Sharp {B}esov capacity estimates for annuli in metric spaces
              with doubling measures},
   JOURNAL = {Math. Z.},
  FJOURNAL = {Mathematische Zeitschrift},
    VOLUME = {305},
      YEAR = {2023},
    NUMBER = {3},
     PAGES = {Paper No. 41, 26},
      ISSN = {0025-5874,1432-1823},
   MRCLASS = {31C45 (30L99 31C12 31C15 31E05)},
  MRNUMBER = {4651122},
MRREVIEWER = {Juha\ K.\ Kinnunen},
       DOI = {10.1007/s00209-023-03360-0},
       URL = {https://doi.org/10.1007/s00209-023-03360-0},
}

@article {Robert,
    AUTHOR = {Wegner, Robert},
     TITLE = {Global-in-time well-posedness of the one-dimensional
              hydrodynamic {G}ross-{P}itaevskii equations without vacuum},
   JOURNAL = {Z. Angew. Math. Phys.},
  FJOURNAL = {Zeitschrift f\"ur Angewandte Mathematik und Physik. ZAMP.
              Journal of Applied Mathematics and Physics. Journal de
              Math\'ematiques et de Physique Appliqu\'ees},
    VOLUME = {74},
      YEAR = {2023},
    NUMBER = {5},
     PAGES = {Paper No. 194, 29},
      ISSN = {0044-2275,1420-9039},
   MRCLASS = {35Q55 (35Q31 37K10 76N10 81Q80)},
  MRNUMBER = {4642722},
MRREVIEWER = {Yong\ Luo},
       DOI = {10.1007/s00033-023-02089-4},
       URL = {https://doi.org/10.1007/s00033-023-02089-4},
}

@article {Bianchini,
    AUTHOR = {Bianchini, Stefano},
     TITLE = {Exact integrability conditions for cotangent vector fields},
   JOURNAL = {Manuscripta Math.},
  FJOURNAL = {Manuscripta Mathematica},
    VOLUME = {173},
      YEAR = {2024},
    NUMBER = {1-2},
     PAGES = {293--340},
      ISSN = {0025-2611,1432-1785},
   MRCLASS = {35Q40 (76Y05)},
  MRNUMBER = {4684350},
MRREVIEWER = {Robert\ M.\ Owczarek},
       DOI = {10.1007/s00229-023-01461-y},
       URL = {https://doi.org/10.1007/s00229-023-01461-y},
}

@article {FeolaIandoliMurgante,
    AUTHOR = {Feola, Roberto and Iandoli, Felice and Murgante, Federico},
     TITLE = {Long-time stability of the quantum hydrodynamic system on
              irrational tori},
   JOURNAL = {Math. Eng.},
  FJOURNAL = {Mathematics in Engineering},
    VOLUME = {4},
      YEAR = {2022},
    NUMBER = {3},
     PAGES = {Paper No. 023, 24},
      ISSN = {2640-3501},
   MRCLASS = {35Q35 (76Y05)},
  MRNUMBER = {4298389},
       DOI = {10.3934/mine.2022023},
       URL = {https://doi.org/10.3934/mine.2022023},
}

@article {AHM23NLS,
    AUTHOR = {Antonelli, Paolo and Hientzsch, Lars Eric and Marcati,
              Pierangelo},
     TITLE = {Finite energy well-posedness for nonlinear {S}chr\"odinger
              equations with non-vanishing conditions at infinity},
   JOURNAL = {Rev. Mat. Iberoam.},
  FJOURNAL = {Revista Matem\'atica Iberoamericana},
    VOLUME = {42},
      YEAR = {2026},
    NUMBER = {3},
     PAGES = {1059--1116},
      ISSN = {0213-2230,2235-0616},
   MRCLASS = {35Q55 (35B30 37L50)},
  MRNUMBER = {5065254},
       DOI = {10.4171/rmi/1584},
       URL = {https://doi.org/10.4171/rmi/1584},
}

@phdthesis{H,
  author = {Hientzsch, Lars Eric},
  title  = {Nonlinear {Schr{\"o}dinger} equations and
            quantum fluids non-vanishing at infinity:
            Incompressible limit and quantum vortices},
  school = {Gran Sasso Science Institute},
  address = {L'Aquila, Italy},
  year   = {2019},
  url    = {https://hdl.handle.net/20.500.12571/26564}
}

@article {KochLiao1,
    AUTHOR = {Koch, Herbert and Liao, Xian},
     TITLE = {Conserved energies for the one dimensional
              {G}ross-{P}itaevskii equation: low regularity case},
   JOURNAL = {Adv. Math.},
  FJOURNAL = {Advances in Mathematics},
    VOLUME = {420},
      YEAR = {2023},
     PAGES = {Paper No. 108996, 61},
      ISSN = {0001-8708,1090-2082},
   MRCLASS = {35Q55 (37K10)},
  MRNUMBER = {4568059},
       DOI = {10.1016/j.aim.2023.108996},
       URL = {https://doi.org/10.1016/j.aim.2023.108996},
}

@article {GiulianiScandone,
    AUTHOR = {Giuliani, Filippo and Scandone, Raffaele},
     TITLE = {Energy cascade and {S}obolev norms inflation for the quantum
              {E}uler equations on tori},
   JOURNAL = {Adv. Math.},
  FJOURNAL = {Advances in Mathematics},
    VOLUME = {479},
      YEAR = {2025},
     PAGES = {Paper No. 110453, 59},
      ISSN = {0001-8708,1090-2082},
   MRCLASS = {76Y05 (35B34 35B35 35Q31 35Q40 37K45)},
  MRNUMBER = {4939421},
MRREVIEWER = {Fucai\ Li},
       DOI = {10.1016/j.aim.2025.110453},
       URL = {https://doi.org/10.1016/j.aim.2025.110453},
}

@article {DonatelliFeireislMarcati,
    AUTHOR = {Donatelli, Donatella and Feireisl, Eduard and Marcati,
              Pierangelo},
     TITLE = {Well/ill posedness for the {E}uler-{K}orteweg-{P}oisson system
              and related problems},
   JOURNAL = {Comm. Partial Differential Equations},
  FJOURNAL = {Communications in Partial Differential Equations},
    VOLUME = {40},
      YEAR = {2015},
    NUMBER = {7},
     PAGES = {1314--1335},
      ISSN = {0360-5302,1532-4133},
   MRCLASS = {35Q35 (35D30 35Q53 76N10 82D50)},
  MRNUMBER = {3341206},
MRREVIEWER = {Rapha\"el\ Danchin},
       DOI = {10.1080/03605302.2014.972517},
       URL = {https://doi.org/10.1080/03605302.2014.972517},
}

@article {MarkowichSierra,
    AUTHOR = {Markowich, Peter and Sierra, Jes\'us},
     TITLE = {Non-uniqueness of weak solutions of the quantum-hydrodynamic
              system},
   JOURNAL = {Kinet. Relat. Models},
  FJOURNAL = {Kinetic and Related Models},
    VOLUME = {12},
      YEAR = {2019},
    NUMBER = {2},
     PAGES = {347--356},
      ISSN = {1937-5093,1937-5077},
   MRCLASS = {76Y05 (35A02 35Q35 81Q05)},
  MRNUMBER = {3918271},
       DOI = {10.3934/krm.2019015},
       URL = {https://doi.org/10.3934/krm.2019015},
}

@article {EvansSpruck1991,
    AUTHOR = {Evans, L. C. and Spruck, J.},
     TITLE = {Motion of level sets by mean curvature. {I}},
   JOURNAL = {J. Differential Geom.},
  FJOURNAL = {Journal of Differential Geometry},
    VOLUME = {33},
      YEAR = {1991},
    NUMBER = {3},
     PAGES = {635--681},
      ISSN = {0022-040X,1945-743X},
   MRCLASS = {35K55 (35B05 35D05 53A10)},
  MRNUMBER = {1100206},
MRREVIEWER = {Friedrich\ Sauvigny},
       URL = {http://projecteuclid.org/euclid.jdg/1214446559},
}

@article {AlfaroEtAl2010,
    AUTHOR = {Alfaro, Matthieu and Garcke, Harald and Hilhorst, Danielle and
              Matano, Hiroshi and Sch\"atzle, Reiner},
     TITLE = {Motion by anisotropic mean curvature as sharp interface limit
              of an inhomogeneous and anisotropic {A}llen-{C}ahn equation},
   JOURNAL = {Proc. Roy. Soc. Edinburgh Sect. A},
  FJOURNAL = {Proceedings of the Royal Society of Edinburgh. Section A.
              Mathematics},
    VOLUME = {140},
      YEAR = {2010},
    NUMBER = {4},
     PAGES = {673--706},
      ISSN = {0308-2105,1473-7124},
   MRCLASS = {53C44 (35B25 35K57 35K59 53C60)},
  MRNUMBER = {2672065},
MRREVIEWER = {J.\ W.\ Jerome},
       DOI = {10.1017/S0308210508000541},
       URL = {https://doi.org/10.1017/S0308210508000541},
}

@article {GarckeKnopfWittmann2023,
    AUTHOR = {Garcke, Harald and Knopf, Patrik and Wittmann, Julia},
     TITLE = {The anisotropic {C}ahn-{H}illiard equation: regularity theory
              and strict separation properties},
   JOURNAL = {Discrete Contin. Dyn. Syst. Ser. S},
  FJOURNAL = {Discrete and Continuous Dynamical Systems. Series S},
    VOLUME = {16},
      YEAR = {2023},
    NUMBER = {12},
     PAGES = {3622--3660},
      ISSN = {1937-1632,1937-1179},
   MRCLASS = {35K55 (35K61 74E10)},
  MRNUMBER = {4677067},
MRREVIEWER = {Alain\ Brillard},
       DOI = {10.3934/dcdss.2023146},
       URL = {https://doi.org/10.3934/dcdss.2023146},
}

@article {BanerjeeGarofalo2015,
    AUTHOR = {Banerjee, Agnid and Garofalo, Nicola},
     TITLE = {On the {D}irichlet boundary value problem for the normalized
              {$p$}-{L}aplacian evolution},
   JOURNAL = {Commun. Pure Appl. Anal.},
  FJOURNAL = {Communications on Pure and Applied Analysis},
    VOLUME = {14},
      YEAR = {2015},
    NUMBER = {1},
     PAGES = {1--21},
      ISSN = {1534-0392,1553-5258},
   MRCLASS = {35K55 (35B40 35B65 35K20)},
  MRNUMBER = {3299022},
MRREVIEWER = {Daria\ E.\ Apushkinskaya},
       DOI = {10.3934/cpaa.2015.14.1},
       URL = {https://doi.org/10.3934/cpaa.2015.14.1},
}

@article {AttouchiParviainen2018,
    AUTHOR = {Attouchi, Amal and Parviainen, Mikko},
     TITLE = {H\"older regularity for the gradient of the inhomogeneous
              parabolic normalized {$p$}-{L}aplacian},
   JOURNAL = {Commun. Contemp. Math.},
  FJOURNAL = {Communications in Contemporary Mathematics},
    VOLUME = {20},
      YEAR = {2018},
    NUMBER = {4},
     PAGES = {1750035, 27},
      ISSN = {0219-1997,1793-6683},
   MRCLASS = {35K55 (35B65 35D40 35K92)},
  MRNUMBER = {3804259},
       DOI = {10.1142/S0219199717500353},
       URL = {https://doi.org/10.1142/S0219199717500353},
}

@article {Klaus2023,
    AUTHOR = {Klaus, Friedrich},
     TITLE = {Wellposedness of {NLS} in modulation spaces},
   JOURNAL = {J. Fourier Anal. Appl.},
  FJOURNAL = {The Journal of Fourier Analysis and Applications},
    VOLUME = {29},
      YEAR = {2023},
    NUMBER = {1},
     PAGES = {Paper No. 9, 37},
      ISSN = {1069-5869,1531-5851},
   MRCLASS = {35Q55 (42B35)},
  MRNUMBER = {4534495},
       DOI = {10.1007/s00041-022-09985-9},
       URL = {https://doi.org/10.1007/s00041-022-09985-9},
}

@article {KlausKunstmann2022,
    AUTHOR = {Klaus, Friedrich and Kunstmann, Peer},
     TITLE = {Global wellposedness of {NLS} in {$H^1(\mathbb{R})+H^s(\mathbb{T})$}},
   JOURNAL = {J. Math. Anal. Appl.},
  FJOURNAL = {Journal of Mathematical Analysis and Applications},
    VOLUME = {514},
      YEAR = {2022},
    NUMBER = {2},
     PAGES = {Paper No. 126359, 14},
      ISSN = {0022-247X,1096-0813},
   MRCLASS = {35Q55},
  MRNUMBER = {4432296},
MRREVIEWER = {Masayuki\ Hayashi},
       DOI = {10.1016/j.jmaa.2022.126359},
       URL = {https://doi.org/10.1016/j.jmaa.2022.126359},
}

@article {Mohamad2014,
    AUTHOR = {Mohamad, Haidar},
     TITLE = {Hydrodynamical form for the one-dimensional
              {G}ross-{P}itaevskii equation},
   JOURNAL = {Electron. J. Differential Equations},
  FJOURNAL = {Electronic Journal of Differential Equations},
      YEAR = {2014},
     PAGES = {No. 141, 27},
      ISSN = {1072-6691},
   MRCLASS = {35Q55 (35B30 35C07 35C08)},
  MRNUMBER = {3239384},
}

@article {Bourgain2002,
    AUTHOR = {Bourgain, Jean},
     TITLE = {On the global {C}auchy problem for the nonlinear
              {S}chr\"odinger equation},
   JOURNAL = {Proc. Natl. Acad. Sci. USA},
  FJOURNAL = {Proceedings of the National Academy of Sciences of the United
              States of America},
    VOLUME = {99},
      YEAR = {2002},
    NUMBER = {24},
     PAGES = {15262--15268},
      ISSN = {0027-8424,1091-6490},
   MRCLASS = {35Q55 (35G25)},
  MRNUMBER = {1946761},
MRREVIEWER = {Victor\ V.\ Zharinov},
       DOI = {10.1073/pnas.222494399},
       URL = {https://doi.org/10.1073/pnas.222494399},
}

@article {KochTataru2018,
    AUTHOR = {Koch, Herbert and Tataru, Daniel},
     TITLE = {Conserved energies for the cubic nonlinear {S}chr\"odinger
              equation in one dimension},
   JOURNAL = {Duke Math. J.},
  FJOURNAL = {Duke Mathematical Journal},
    VOLUME = {167},
      YEAR = {2018},
    NUMBER = {17},
     PAGES = {3207--3313},
      ISSN = {0012-7094,1547-7398},
   MRCLASS = {35Q55 (35Q53 37K05 37K10)},
  MRNUMBER = {3874652},
MRREVIEWER = {Thierry\ Cazenave},
       DOI = {10.1215/00127094-2018-0033},
       URL = {https://doi.org/10.1215/00127094-2018-0033},
}

@article {BethuelSmets,
    AUTHOR = {Bethuel, Fabrice and Smets, Didier},
     TITLE = {A remark on the {C}auchy problem for the 2{D}
              {G}ross-{P}itaevskii equation with nonzero degree at infinity},
   JOURNAL = {Differential Integral Equations},
  FJOURNAL = {Differential and Integral Equations. An International Journal
              for Theory \& Applications},
    VOLUME = {20},
      YEAR = {2007},
    NUMBER = {3},
     PAGES = {325--338},
      ISSN = {0893-4983},
   MRCLASS = {35Q40 (35A05 35Q55 47J30)},
  MRNUMBER = {2293989},
MRREVIEWER = {Michele\ Correggi},
}

@article {KillipOhPocovnicuVisan2012,
    AUTHOR = {Killip, Rowan and Oh, Tadahiro and Pocovnicu, Oana and Vi\c
              san, Monica},
     TITLE = {Global well-posedness of the {G}ross-{P}itaevskii and
              cubic-quintic nonlinear {S}chr\"odinger equations with
              non-vanishing boundary conditions},
   JOURNAL = {Math. Res. Lett.},
  FJOURNAL = {Mathematical Research Letters},
    VOLUME = {19},
      YEAR = {2012},
    NUMBER = {5},
     PAGES = {969--986},
      ISSN = {1073-2780,1945-001X},
   MRCLASS = {35Q55 (35B30)},
  MRNUMBER = {3039823},
MRREVIEWER = {Gianmaria\ Verzini},
       DOI = {10.4310/MRL.2012.v19.n5.a1},
       URL = {https://doi.org/10.4310/MRL.2012.v19.n5.a1},
}

@article {DydaKijaczko2022,
    AUTHOR = {Dyda, Bartłomiej and Kijaczko, Michael},
     TITLE = {On density of compactly supported smooth functions in
              fractional {S}obolev spaces},
   JOURNAL = {Ann. Mat. Pura Appl. (4)},
  FJOURNAL = {Annali di Matematica Pura ed Applicata. Series IV},
    VOLUME = {201},
      YEAR = {2022},
    NUMBER = {4},
     PAGES = {1855--1867},
      ISSN = {0373-3114,1618-1891},
   MRCLASS = {46E35 (26D15 35A15)},
  MRNUMBER = {4454384},
MRREVIEWER = {Aigerim\ Kalybay},
       DOI = {10.1007/s10231-021-01181-8},
       URL = {https://doi.org/10.1007/s10231-021-01181-8},
}

@misc{Klaus2026,
  author = {Klaus, Friedrich},
  title  = {Global well-posedness of defocusing cubic {NLS}
            in {$M^{\infty,1}(\mathbb{R})$}},
  year   = {2026},
  note   = {Preprint, \href{https://arxiv.org/abs/2609.29147}
                          {arXiv:2609.29147}}
}

@misc{IgnatevNazarovNichitenkoTursunbaev2026,
  author = {Ignatev, Egor and Nazarov, Alexander I. and
            Nichitenko, Pavel and Tursunbaev, Artur},
  title  = {Truncations for fractional {Laplacians}},
  year   = {2026},
  note   = {Preprint, \href{https://arxiv.org/abs/2608.05433}
                          {arXiv:2608.05433}}
}

@misc{KotaniWangXuZhang2026,
  author = {Kotani, Shinichi and Wang, Yiqian and
            Xu, Jiahao and Zhang, Shuo},
  title  = {Global flows for the cubic nonlinear
            {Schr{\"o}dinger} equation without decay assumptions},
  year   = {2026},
  note   = {Preprint, \href{https://arxiv.org/abs/2609.23324}
                          {arXiv:2609.23324}}
}

@misc{PineauTaylor2025,
  author = {Pineau, Ben and Taylor, Mitchell A.},
  title  = {On the optimal {Sobolev} threshold for evolution
            equations with rough nonlinearities},
  year   = {2025},
  note   = {Preprint, \href{https://arxiv.org/abs/2505.14966}
                          {arXiv:2505.14966}}
}

@article{BorowskiCianchi2025,
  author  = {Borowski, Micha{\l} and Cianchi, Andrea},
  title   = {Composition operators in {Orlicz}--{Sobolev} spaces},
  journal = {Commun. Contemp. Math.},
  year    = {2026},
  note    = {Published online,
             \href{https://doi.org/10.1142/S0219199726500355}
                  {doi:10.1142/S0219199726500355}}
}

@article{AHMZ,
  author  = {Antonelli, Paolo and Hientzsch, Lars Eric and
             Marcati, Pierangelo and Zheng, Hao},
  title   = {On some results for quantum hydrodynamical models},
  journal = {RIMS K{\^o}ky{\^u}roku},
  volume  = {2070},
  year    = {2018},
  pages   = {107--129}
}

@incollection{BourgainBrezisMironescu2001,
  author    = {Bourgain, Jean and Brezis, Ha{\"i}m and Mironescu, Petru},
  title     = {Another look at {Sobolev} spaces},
  booktitle = {Optimal Control and Partial Differential Equations},
  editor    = {Menaldi, Jos{\'e} Luis and Rofman, Edmundo and Sulem, Agn{\`e}s},
  publisher = {IOS Press},
  address   = {Amsterdam},
  year      = {2001},
  pages     = {439--455},
  url       = {https://sites.math.rutgers.edu/~brezis/PUBlications/167book.pdf}
}

@article {Madelung1927,
    AUTHOR = {Madelung, Erwin},
     TITLE = {Quantentheorie in hydrodynamischer {F}orm},
   JOURNAL = {Z. Phys.},
  FJOURNAL = {Zeitschrift f\"ur Physik},
    VOLUME = {40},
      YEAR = {1927},
    NUMBER = {3-4},
     PAGES = {322--326},
       DOI = {10.1007/BF01400372},
       URL = {https://doi.org/10.1007/BF01400372},
}

\end{document}